\documentclass[11pt]{article}
\usepackage[margin=1in]{geometry}
\usepackage{authblk}
\usepackage{graphicx}
\usepackage{amsfonts}
\usepackage{mathtools}
\usepackage{amsmath}
\usepackage{amssymb}
\usepackage{amsthm}
\usepackage{booktabs}
\usepackage{placeins}
\usepackage{natbib}
\usepackage[dvipsnames]{xcolor}
\usepackage[colorlinks=true,linkcolor=NavyBlue,citecolor=NavyBlue,urlcolor=NavyBlue]{hyperref}

\allowdisplaybreaks

\newcommand{\yes}{\ensuremath{\checkmark}}
\newcommand{\no}{\ensuremath{-}}
\newcommand{\prt}{\ensuremath{\sim}}

\renewcommand{\P}{\mathbb{P}}

\newcommand{\F}{\mathcal{F}}

\newcommand{\R}{\mathbb{R}}
\newcommand{\E}{\mathbb{E}}

\newcommand{\bm}{\boldsymbol}
\newcommand{\argmin}{\operatornamewithlimits{arg\,min}}

\theoremstyle{plain}
\newtheorem{theorem}{Theorem}[section]
\newtheorem{lemma}[theorem]{Lemma}
\newtheorem{corollary}[theorem]{Corollary}
\newtheorem{proposition}[theorem]{Proposition}
\theoremstyle{definition}
\newtheorem{definition}[theorem]{Definition}
\newtheorem{assumption}[theorem]{Assumption}

\theoremstyle{remark}
\newtheorem{remark}[theorem]{Remark}

\numberwithin{equation}{section}

\title{Training Neural Networks to be Locally Adaptive}
\author[1]{Bilol Banerjee}
\author[2]{Soham Sarkar}
\affil[1]{Department of Statistics and Data Science, National University of Singapore}
\affil[2]{Interdisciplinary Statistical Research Unit, Indian Statistical Institute}
\date{\today}

\begin{document}

\maketitle

\begin{abstract}
Deep neural networks fitted by global empirical risk minimization commit to a single model for the entire covariate space. In a regression problem, if the regression function $f_0$ is inhomogeneous --- varying rapidly in one region, while being nearly flat in another --- no single model can be appropriate everywhere, and the fit degrades where the function varies most rapidly. Here we study an alternative in which, instead of fitting the network globally, we fit it locally. For each query point $x$, a \emph{local neural network} is trained on the $k$ nearest-neighbours of $x$ alone, making the training sample and the fitted function query-dependent. We analyze the performance of the resulting estimator through a localization approach, establishing that the pointwise risk of the estimator with an appropriate choice of $k$ attains the rate $n^{-2\alpha/(2\alpha+d^{*})}$, up to logarithmic factors, when $f_0$ is locally $\alpha$-H\"older with $\alpha\in(0,1]$ and the covariates are supported on a $d^{*}$-dimensional manifold embedded in $\R^{d}$. A matching Le Cam lower bound shows that this rate is pointwise minimax optimal. However, the local rate-optimality requires $k$ to be calibrated to the smoothness of $f_0$, which is unknown in practice. We therefore provide a locally adaptive choice of $k$ that makes the neural network adapt to the local H\"older exponent of $f_0$ and to the intrinsic dimension simultaneously, both of which are potentially varying over the covariate domain. An extensive simulation study demonstrates the dominance of the locally adaptive neural network over a globally trained network and several state-of-the-art methods.

\medskip
\noindent\textbf{Keywords:} Deep neural networks; $k$-nearest neighbours; local adaptivity; local Gaussian complexity; minimax rate optimality.
\end{abstract}

\section{Introduction}

Consider the nonparametric regression model
\[
Y_i = f_0(X_i) + \varepsilon_i, \qquad i = 1, \ldots, n,
\]
where $X_i \in [0,1]^d$ are $d$-dimensional covariates, $f_0:[0,1]^d \to \mathbb{R}$ is the (unknown) regression function, and $\varepsilon_i$ are random errors satisfying \(\mathbb{E}(\varepsilon_i \mid X_i) = 0.\) The objective of nonparametric regression is to estimate the function $f_0$ from the observed data $\{(X_i, Y_i)\}_{i=1}^n$, which are assumed to be i.i.d.\ samples from the distribution of $(X, Y)$. Classical approaches to this problem include kernel smoothing, spline regression, local polynomial regression, and nearest-neighbour methods \citep[see][]{Hastie2009ESL}. These techniques provide flexible alternatives to parametric models but often suffer from several challenges. These include the choice of tuning parameters (e.g., bandwidth selection in kernel methods), boundary bias effects that reduce accuracy near the edges of the covariate space, and computational complexity when applied to large datasets. Among these, the curse of dimensionality is particularly severe: As the dimension $d$ increases, local neighbourhoods become sparse, leading to slow convergence rates and poor estimation accuracy \citep{Gyorfietal2002}. Together, these limitations necessitate the development of modern approaches that can better adapt to high-dimensional and complex data structures.

Deep neural networks (DNNs) have emerged as a flexible alternative to the classical approaches, being capable of approximating a rich class of multivariate functions while adapting to underlying structural properties \citep{schmidthieber2019manifold,bauer2019deep,kohler2021rate,farrell2021deep}. \cite{dnnregression} presents a detailed theoretical exposition on estimating the regression function $f_0$ by a class of fully connected feed-forward neural networks with ReLU activation function. Let $\varrho:\mathbb{R}\to\mathbb{R}$ denote the ReLU activation function, defined as \(\varrho(t)=\max\{t,0\}.\) For integers $L \geq 1$ and a width vector $\bm p = (p_1, \ldots, p_L)$, consider a fully connected feed-forward neural network of depth $L$ defined recursively as follows. For input $x \in \mathbb{R}^d$, set \(h^{(0)} = x,\) and for each $\ell = 1, \ldots, L$ define the hidden layer
\[
h^{(\ell)} = \varrho\!\left(W_\ell h^{(\ell-1)} + v_\ell\right),
\]
where $W_\ell \in \mathbb{R}^{p_\ell \times p_{\ell-1}}$ is the weight matrix and $v_\ell \in \mathbb{R}^{p_\ell}$ is the bias vector. Here, $p_0=d$ and $p_{L+1}=1$. The network output is then given by
\[
f(x) = W_{L+1} h^{(L)} + v_{L+1},
\]
with $W_{L+1}\in\mathbb R^{1\times p_L}$ and $v_{L+1}\in\mathbb R$. Write $\theta$ for the vector of all weights and biases, including the output layer, and let $\mathcal F(L,\bm p)$ be the class of networks with this architecture. Throughout the paper, we use a constrained class of DNNs:
\begin{equation}
\label{eq:network-class}
\begin{split}
\mathcal F_n=\mathcal F(L,\bm p,s,F,K)
=\bigg\{f_\theta\in\mathcal F(L,\bm p):\;&
\sum_{j=1}^{L+1}\big(\|W_j\|_0+\|v_j\|_0\big)\le s,\\
&\|\theta\|_\infty\le1,\quad
\|f_\theta\|_{\infty,[0,1]^d}\le F,\quad
\operatorname{Lip}_{[0,1]^d}(f_\theta)\le K\bigg\}.
\end{split}
\end{equation}
Here $F\ge B>0$ and $K>0$ are fixed constants,
independent of $n$. The architecture (viz., the depth $L$, the width $\bm p$) and sparsity $s$ may depend on $n$.
We assume $\|f_0\|_\infty\le B$ with $0<B\le1$~--- a general bound reduces to this case by rescaling the responses. This is the class of \cite{dnnregression} with the additional constraint $\operatorname{Lip}_{[0,1]^d}(f_\theta)\le K$, where $\operatorname{Lip}_{[0,1]^d}(f)=\sup_{u\ne v\in[0,1]^d}|f(u)-f(v)|/\|u-v\|$.
Note that the Lipschitz constraint is part of the estimator analysed here, and does not
require knowledge of the target smoothness or intrinsic dimension.
For brevity, we also write $\mathcal F(L,\bm p,s)$, with $F,K$ understood. 
Estimation of the regression function $f_0$ is performed by minimizing the empirical risk
\[
\hat f_n(\cdot) \in \argmin_{f \in \mathcal{F}(L,\bm p,s)} \frac{1}{n}\sum_{i=1}^n \big(Y_i - f(X_i)\big)^2.
\]

The results of \citet{dnnregression} are stated for regression functions with a hierarchical composition structure. This class contains the unrestricted $d$-variate H\"older functions, for which the results recover the usual $d$-dimensional minimax rate $n^{-2\alpha/(2\alpha+d)}$ up to logarithmic factors. Faster rates require additional structure, such as a composition of functions of fewer variables. Without such additional structures, a globally trained network is therefore subject to the curse of dimensionality, %
and their convergence behavior is not theoretically guaranteed to improve upon traditional nonparametric estimators. We will subsequently establish that this undesirable property can be alleviated by appropriately modifying the training procedure of the DNN. But first, we consider a motivating example.

\subsection{Global and local fits on estimating the Doppler function}
\label{sec:doppler}

We generate $n=1000$ observations from $Y_i=f_0(X_i)+\sigma W_i$ with $X_i$ uniform on $[0,1]$,
$W_i\sim N(0,1)$, $\sigma=0.1$, and $f_0$ the Doppler function
\begin{equation}
\label{eq:doppler}
f_0(x)=\sqrt{x(1-x)}\,\sin\!\Big(\frac{2.1\pi}{x+0.05}\Big),
\qquad x\in[0,1].
\end{equation}
The Doppler function is a standard test function with rapid oscillations near zero and slower variation
farther away. The absolute value of its phase derivative is bounded by $840\pi$, so the frequency
does not diverge at zero. Moreover, the sine factor vanishes at both endpoints:
\[
f_0(u)=-840\pi u^{3/2}+O(u^{5/2})\quad\text{and}\quad
f_0(1-u)=\frac{40\pi}{21}u^{3/2}+O(u^{5/2}) \quad\text{as } u \to 0.
\]
Thus $f_0$ is globally Lipschitz.

\begin{figure}[htbp]
\centering
\begin{tabular}{cc}
\includegraphics[width=0.48\textwidth]{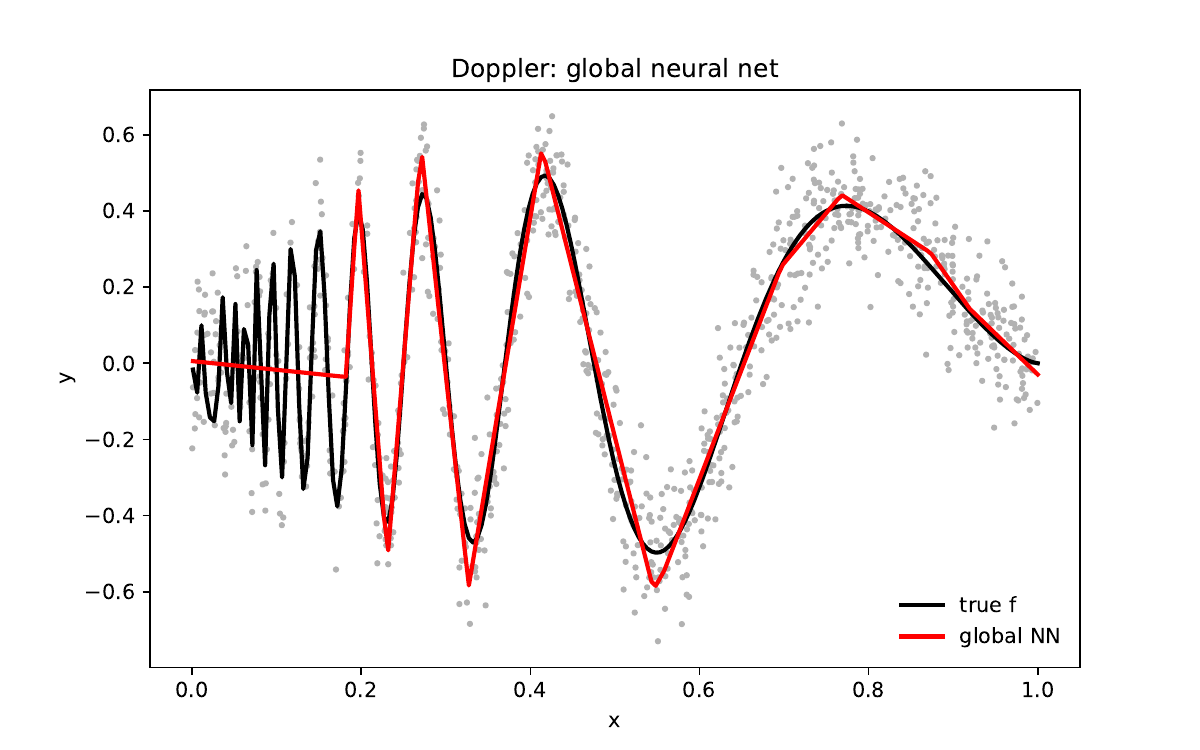} &
\includegraphics[width=0.48\textwidth]{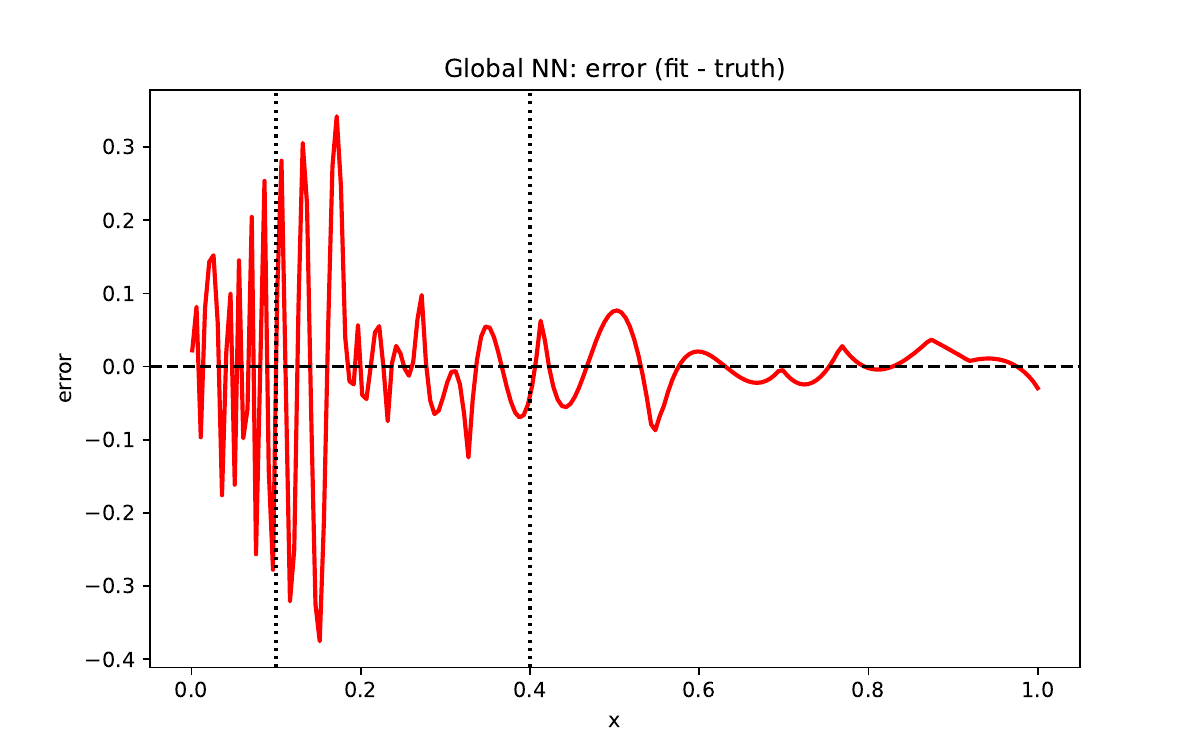}\\
(a) & (b)
\end{tabular}
\caption{Global fit to the Doppler function \eqref{eq:doppler},
$n=1000$, $\sigma=0.1$, test grid of 200 equally spaced points.
(a) Data (grey), $f_0$ (black) and a globally trained ReLU fit (red).
(b) The pointwise error $\hat f(x)-f_0(x)$ of that fit.}
\label{fig:doppler}
\end{figure}

The results from fitting DNNs to the data are shown in Figure~\ref{fig:doppler}. Panel (a) shows a fully connected ReLU network with one hidden layer containing 16 hidden units trained on the whole sample. Away from the origin the fit is accurate; near the origin it simply smooths the oscillations away. Panel (b) plots the point-wise bias of the fit against $x$. The bias is mostly concentrated on the region $x\lesssim0.2$, where the global model fails to capture the underlying variation of $f_0$. Enlarging the network does not resolve this --- a network large enough to track the oscillations near $0$ is more flexible than needed elsewhere and overfits there, because global least squares weights all $n$ residuals equally and cannot vary the weights across the domain.

The above difficulty is not merely that the guarantees of \cite{dnnregression} require a compositional structure. It is that a network fitted by \emph{global} empirical risk minimisation selects one architecture, one sparsity level and one effective degree of smoothing for the entire covariate space, whereas the difficulty of the estimation problem is typically not constant across that space. This suggests estimating $f_0$ one query point at a time, using only the observations near that point.

\subsubsection{The Proposed Estimator}

In contrast to existing approaches that train a DNN globally, we propose a localized training scheme that adapts to the query point. The central idea is to combine the flexibility of DNNs with the locality principle underlying $k$-nearest neighbour ($k$-NN) methods. Specifically, for a given query point $x \in [0,1]^d$, we first identify its $k$ nearest neighbours among the training covariates $\{X_i\}_{i=1}^n$ in the given dataset with respect to a suitable metric (e.g., Euclidean distance). Denote this neighbourhood by \(N_k(x) = \{ i : X_{i} \text{ is among the $k$ nearest neighbours of } x \}.\) Instead of fitting the neural network globally on the full sample, we restrict training to the localized dataset $N_k(x)$. The resulting estimator is then given by
\[
\widehat{f}(x) = \widehat{g}_x(x)\quad\text{where}\quad\widehat{g}_x = \argmin_{f \in \mathcal{F}(L,\bm p,s)} \frac{1}{k}\sum_{i\in N_k(x)} \big(Y_i - f(X_i)\big)^2.
\]

This approach has several advantages. First, by focusing on local neighbourhoods, the estimator adapts to heterogeneity in the regression function and can capture local smoothness properties more effectively. Second, the radius of a $k$-nearest-neighbour ball shrinks at a rate set by the intrinsic dimension $d^{*}$ of the support, not by $d$, so the estimator adapts to low-dimensional structure in the covariates (Section~\ref{sec:rate}). Third, the localized training scheme naturally incorporates ideas from nonparametric regression such as kernel and nearest-neighbour smoothing, while leveraging the expressive power of DNNs to approximate complex functional relationships.

Localized training strategies have already found applications in modern practical AI systems. For instance, $k$-NN-style neighbourhood selection has been applied in computer vision problems like handwritten digit recognition, where predictions are based on local subsets of similar feature vectors \citep{Yuchun1991}. In recommender systems, user-based and item-based $k$-NN approaches restrict training or inference to the nearest neighbours of a query user or item, thereby improving personalization and efficiency \citep{sarwar2001}. Our proposed method is also closely related to the Lazy Learning paradigm \citep{Bontempi2002}. Lazy Learning techniques avoid fitting a global model over the entire dataset. Instead, they retain the training data in memory and construct local models only when a query point is presented. The prediction is then obtained by interpolating or fitting the data in the neighbourhood of the query. This philosophy aligns naturally with our $k$-NN-inspired approach, where for each query $x$ we select its $k$ nearest neighbours and train a neural network on this local subset. 

\begin{figure}[htbp]
\centering
\begin{tabular}{cc}
\includegraphics[width=0.48\textwidth]{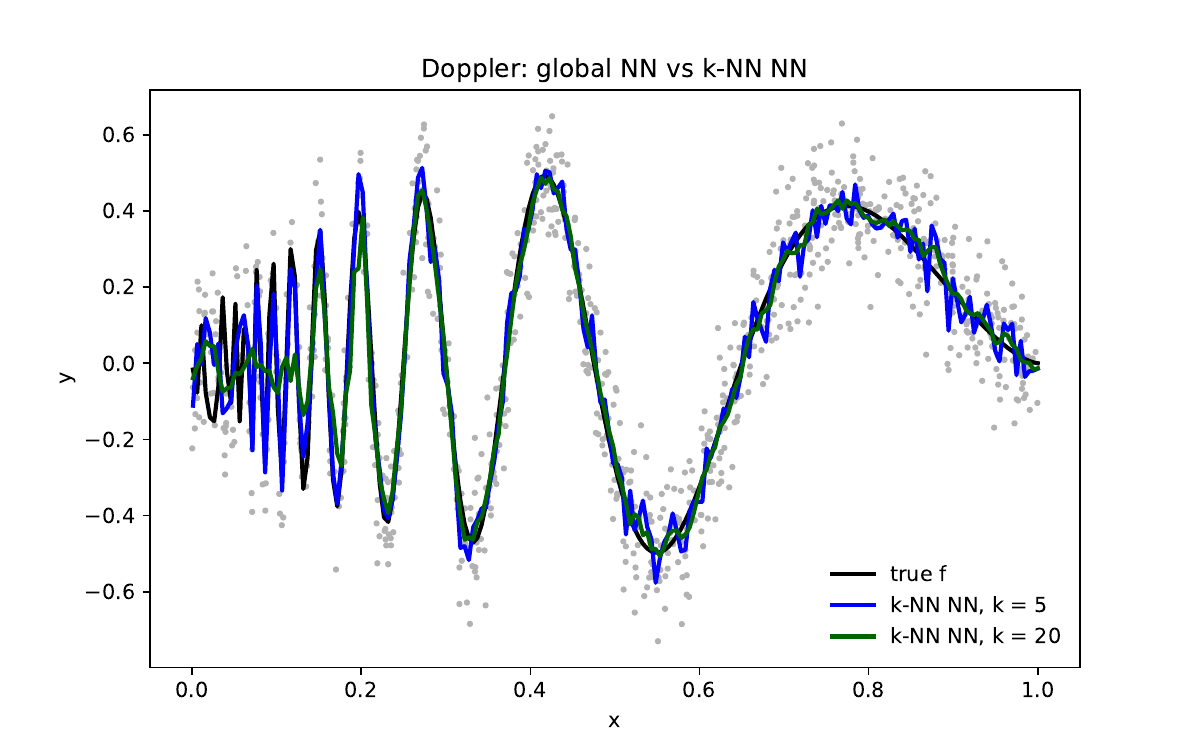} &
\includegraphics[width=0.48\textwidth]{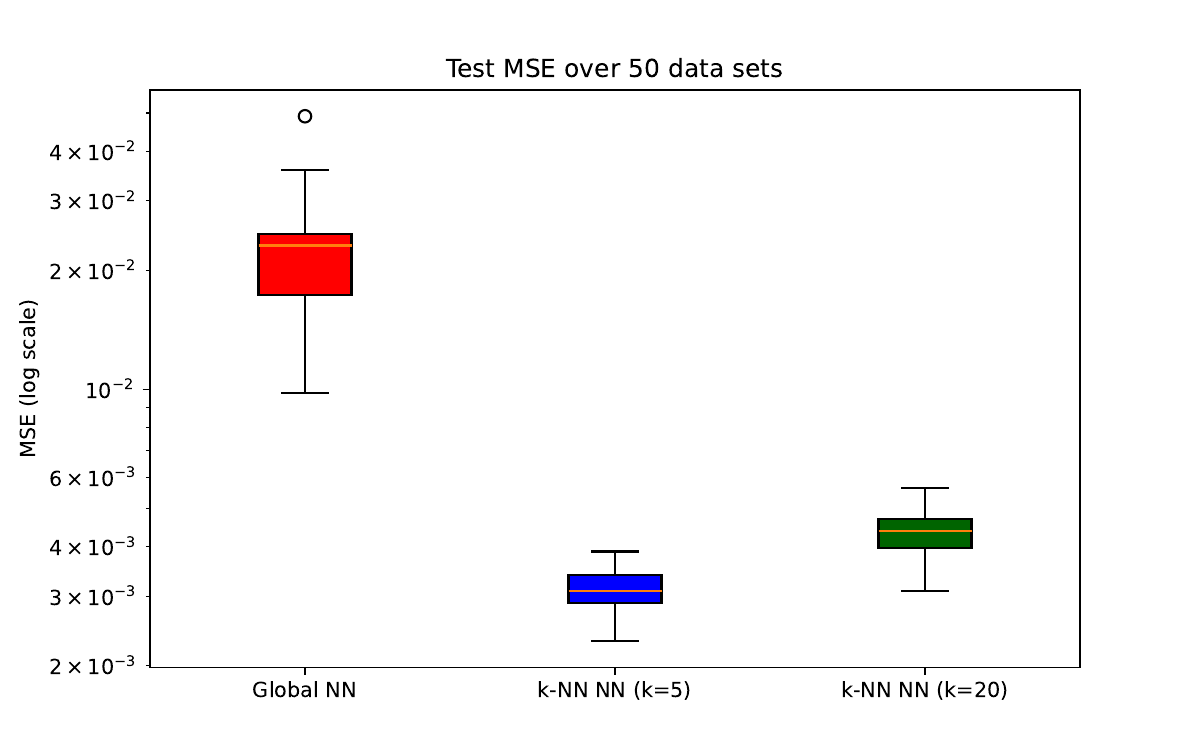}\\
(a) & (b)
\end{tabular}
\caption{Local fits to the Doppler function \eqref{eq:doppler}, same setting as Figure~\ref{fig:doppler}.
(a) Local ReLU fits with one hidden layer of width 5, using $k=5$ and $k=20$.
(b) Test MSE over 50 replicated data sets; each MSE averages squared errors over the test grid.}
\label{fig:doppler-local}
\end{figure}

Figure~\ref{fig:doppler-local}(a) shows the results for two localized fits, obtained by locally fitting a DNN with the same architecture for two different values of $k$. The model fitted on a smaller neighbourhood ($k=5$) follows the oscillations near the origin but is visibly noisy where $f_0$ is smooth. The model fitted on a larger neighbourhood ($k=20$) is accurate on the smooth part, but fails to capture the oscillations near the origin.

Figure~\ref{fig:doppler-local}(b) displays the box plot of the mean-squared errors (MSE) over $50$ replicated data sets of the global neural network estimator and the localized versions. The best localized fit has median MSE $3\times10^{-3}$ against $2.3\times10^{-2}$ for the global network; a factor of $7.6$ improvement. This shows that localization can improve the fit significantly. The best fixed choice here is $k=5$, but it may not be optimal in general. It provides visibly rough estimates when $f_0$ is smooth, whereas $k=20$ gives a much smoother estimate in that region, although its MSE is higher than the fit with $k=5$. Hence, no fixed value of $k$ can be right everywhere -- the optimal neighbourhood size should be governed by the local smoothness of $f_0$ at the query point.

Hence, two questions follow:
\begin{itemize}
    \item Does localized training actually help? What is the pointwise risk of the localized estimator, and how does it compare with the minimax benchmark?
    \item How should $k$ be chosen from the data?
\end{itemize}
These questions organize the rest of the paper, but before that we present an exposition of the current trends in the related literature.

\subsection{Related Works}
\label{sec:related}

In this section, we explicitly lay out what some of the popular contributions in the literature supply and where the present paper departs from it. Also, the word ``local'' is used in the literature in two different senses. An estimator is local \emph{in covariate space} when it uses only observations near the query point, as $k$-NN and kernel methods do and as our estimator does; an empirical process is localised \emph{in function space} when it is controlled on a small ball around $f_0$ rather than on the whole class, which is the technique of \citet{bartlett2005local} and \citet{koltchinskii2006local} that we use in Section~\ref{sec:rate}. Both appear in this article, and they are unrelated.

\citet{dnnregression} studies least squares over a sparse
ReLU class, proves an oracle inequality whose stochastic term is an entropy term
$(s+1)\log(nL\|\bm p\|_\infty)/n$, bounds the covering numbers of the class under the
normalisation $\max_\ell\|W_\ell\|_\infty\le1$, and combines these with an approximation
theorem for composite H\"older functions to obtain minimax-optimal rates.
\citet{bauer2019deep} obtain analogous rates for smooth activations under a hierarchical
composition model, and \citet{kohler2021rate} show that the sparsity constraint can be
dispensed with entirely if the network is deep enough; \citet{farrell2021deep} give
non-asymptotic bounds for unconstrained multilayer perceptrons and use them for
semiparametric inference. However, every one of these results concerns a single estimator that minimises
an empirical risk over the whole sample, and in every one the rate is governed by a single
global smoothness index. A function that is smooth on most of its domain and rough on a
small part is treated as globally rough. That is precisely the behaviour
Figure~\ref{fig:doppler} exhibits. In this article, we borrow the ideas from the above articles, but analyze the terms using the $k$-NN empirical process approach. Our procedure is significantly different from the nearest neighbour empirical process developed in \cite{portier2025nearest}, as for our analysis the summand and the conditioning variable are the same.

A substantial body of work argues that globally trained networks are already locally adaptive. \citet{suzuki2019adaptivity} shows that
functions in a Besov space $B^\alpha_{p,q}$ with $p<2$ --- the classical home of
inhomogeneous smoothness --- admit sparse expansions in a B-spline wavelet dictionary whose
elements ReLU subnetworks can realise at logarithmic cost, so that empirical risk
minimisation over the resulting class attains the spatially adaptive rate that no
\emph{linear} estimator can. \citet{imaizumi2019deep,imaizumi2022advantage} make the
analogous point for piecewise smooth targets with smooth discontinuity surfaces, and
\citet{suzuki2021deep} extend it to anisotropic smoothness. \citet{parhi2023near} show
that weight decay on a shallow ReLU network is, after the usual rescaling, exactly a
second-order Radon-domain total-variation penalty --- the same penalty that makes locally
adaptive regression splines locally adaptive --- and that estimation over the resulting ball
is near minimax optimal while every linear estimator suffers the curse of dimensionality.
\citet{zhang2023deep} generalise the algebra to depth, showing that $\ell_2$ weight decay on
an $L$-layer parallel ReLU network is equivalent to an $\ell_{2/L}$ penalty on the
coefficients of the learned dictionary, and obtain rates over Besov and bounded-variation
classes that approach the minimax rate exponentially fast in $L$. Most recently,
\citet{liu2025deep} define a function class through tree-based nonlinear approximation and
show that networks adapt to regularity that varies across location and scale, and to the
design distribution.
However, all of these guarantees are stated for an integrated $L_2$ risk. They say nothing about $\hat f_n(x)$ at a
particular rough $x$, which is the quantity of interest; an estimator can be minimax optimal in
integrated risk while being badly wrong on a small region, and a small region is exactly
where inhomogeneous smoothness lives. On the other hand, the adaptivity in these results is
\emph{approximation-theoretic and class-dependent}: there exists a network class, indexed
through its size, sparsity or regularisation strength by the smoothness parameters, over
which the idealised minimiser attains the rate. Ours is \emph{procedural}: $\hat k(x)$ is
computed from the data at each query point by \eqref{eq:lepski-rule}, and the smoothness
index appears nowhere in the algorithm. Moreover, none of these results is about the estimator
that training actually returns. \citet{liang2025stable} show that the minima of ReLU networks that gradient descent can stably occupy suffer a curse of dimensionality through a ``neural shattering'' phenomenon,
and \citet{jiang2026depth} construct a local packing establishing that the minimax risk over
the deep variation-norm class grows quadratically in depth, so that the depth deficit in
\citet{zhang2023deep} is not an artefact of a loose upper bound. The gap between the
idealised regularised minimiser and what is available in practice is real, and
Figure~\ref{fig:doppler} is an instance of it.

The classical analysis of $k$-NN regression \citep{Gyorfietal2002,biau2015lectures} balances
a variance $\sigma^2/k$ against a bias controlled by the $k$-NN radius
$R_k(x)\asymp(k/n)^{1/d}$. \citet{kpotufe2011knn} observed that under a local doubling
condition the radius is automatically $(k/n)^{1/d(x)}$ with $d(x)$ the local doubling
exponent, so that $k$-NN attains the rate $n^{-2/(2+d(x))}$ with nothing in the procedure
estimating $d(x)$. \citet{kpotufe2013adaptivity} added a Lepski-type bandwidth rule on top,
obtaining a kernel estimator that adapts simultaneously to the local H\"older exponent and
the local dimension. Variable-bandwidth local polynomial regression
\citep{fan1996local,fan1996study} is the classical route to the same goal through plug-in
selection, and forest methods \citep{athey2019generalized,friedberg2021local} learn an adaptive neighbourhood, then fit a flexible model in it. In all of this work the local model is a constant,
a kernel-weighted polynomial or a forest average. Ours is a neural network, which is what
gives the local fit approximation power beyond polynomials on the neighbourhood. The resulting estimator is not a weighted average of the responses, which is also what makes the analysis harder.%

The rule for choosing $k$ using the data goes back to \citet{lepskii1990problem,lepskii1991asymptotically}: order the
candidate bandwidths, and accept the largest one whose estimate is still statistically
indistinguishable from every smaller one. \citet{lepski1997optimal} applied it pointwise with
variable-bandwidth kernel estimates and obtained the spatial adaptivity of wavelet
thresholding by purely local means; \citet{lepski1997pointwise} established the matching
lower bound, showing that a logarithmic factor is \emph{unavoidable} for pointwise
adaptation. \citet{goldenshluger2011bandwidth} gave the modern form with genuine oracle
inequalities. The $\log$ factor in \eqref{eq:adaptive-rate} is the price
\citet{lepski1997pointwise} showed every pointwise adaptive estimator must pay to attain local adaptivity.
Here, the local parameter $k$ is a neighbourhood size, and the admissibility argument is handled by a stochastic scale derived from network entropy.

Wavelet shrinkage \citep{donoho1994ideal,donoho1995adapting,donoho1995asymptopia}, locally adaptive regression splines
\citep{mammen1997locally} and trend filtering \citep{tibshirani2014adaptive} are some of the popular spatially adaptive regression methods. Most recently \citet{chatterjee2024laser} propose variable-bandwidth local polynomial
regression with a split-based rule and a single global threshold, with pointwise bounds
adapting to the local H\"older exponent; \citet{kuzborskij2020locally} treat the
corresponding question in the sequential setting. These methods are used as 
competitors in Section~\ref{sec:simulation} in our simulation experiments. However, all of them are intrinsically one-dimensional and tied to an interval or dyadic structure --- wavelet shrinkage additionally requires an equispaced
design --- whereas our estimator makes no requirement on the design at all and is stated for covariates supported on a manifold.

\citet{nakada2020adaptive} obtain rates governed by the Minkowski dimension of the support,
and \citet{schmidthieber2019manifold} and \citet{chen2022nonparametric} give chart-based
network constructions for manifold-supported covariates with rates in the intrinsic
dimension $d^{*}$. The intrinsic-dimension adaptivity for
networks is not new, and we do not claim it as such. In
these constructions the chart structure is built into the approximating class, so the class
must be chosen knowing $d^{*}$; here, as in \citet{kpotufe2011knn}, the $k$-NN radius simply
is $(k/n)^{1/d^{*}}$ whatever the support looks like, and the procedure does not know
$d^{*}$, or even that the support is a manifold rather than a set of small Minkowski
dimension.

\subsection{Our Contributions}

\begin{table}[h!]
\centering
\caption{Where the present paper sits. \emph{Network local model:} the function fitted in
the region that determines $\hat f(x)$ is a neural network. \emph{Pointwise:} the guarantee
bounds $|\hat f(x)-f_0(x)|$ at a fixed $x$ rather than an integrated risk. \emph{Data-driven
adaptation:} the estimator computes its own tuning parameter from the data, as opposed to a
smoothness-dependent choice made inside the proof; $\prt$ marks a single global
regularisation parameter that is tuned but is not smoothness-indexed. \emph{Adapts to
$d^{*}$:} the rate involves an intrinsic rather than the ambient dimension.}
\label{tab:related}
\small
\begin{tabular}{lcccc}
 & network & & data-driven & adapts \\
reference & local model & pointwise & adaptation & to $d^{*}$ \\
\midrule
\citet{dnnregression}            & \yes & \no  & \no  & \no  \\
\citet{bauer2019deep}            & \yes & \no  & \no  & \no  \\
\citet{kohler2021rate}           & \yes & \no  & \no  & \no  \\
\citet{suzuki2019adaptivity}     & \yes & \no  & \no  & \no  \\
\citet{imaizumi2019deep}         & \yes & \no  & \no  & \no  \\
\citet{suzuki2021deep}           & \yes & \no  & \no  & \prt \\
\citet{parhi2023near}            & \yes & \no  & \prt & \no  \\
\citet{zhang2023deep}            & \yes & \no  & \prt & \no  \\
\citet{liu2025deep}              & \yes & \no  & \no  & \yes \\
\citet{chen2022nonparametric}    & \yes & \no  & \no  & \yes \\
\citet{nakada2020adaptive}       & \yes & \no  & \no  & \yes \\
\midrule
\citet{donoho1994ideal}          & \no  & \no  & \yes & \no  \\
\citet{mammen1997locally}        & \no  & \no  & \prt & \no  \\
\citet{tibshirani2014adaptive}   & \no  & \no  & \prt & \no  \\
\citet{lepski1997optimal}        & \no  & \yes & \yes & \no  \\
\citet{kpotufe2011knn}           & \no  & \yes & \no  & \yes \\
\citet{kpotufe2013adaptivity}    & \no  & \yes & \yes & \yes \\
\citet{chatterjee2024laser}      & \no  & \yes & \yes & \no  \\
\midrule
this paper                       & \yes & \yes & \yes & \yes \\
\bottomrule
\end{tabular}
\end{table}

We summarize the current literature and our contributions in Table~\ref{tab:related} on four broad points. Reading Table~\ref{tab:related} by rows, the following are established. \emph{(i)} Least squares over a suitably chosen sparse ReLU class attains minimax rates in integrated risk for composite H\"older targets \citep{dnnregression}, and with weight decay or a parallel architecture attains near-minimax rates over Besov and bounded-variation classes \citep{parhi2023near,zhang2023deep}. \emph{(ii)} Globally trained networks adapt to the intrinsic dimension of the support \citep{nakada2020adaptive,chen2022nonparametric} and, in an approximation-theoretic sense, to regularity that varies over the domain \citep{suzuki2019adaptivity,liu2025deep}. \emph{(iii)} Lepski-type selection of a locality parameter, with pointwise guarantees that adapt to the local smoothness exponent and to the local dimension, is available for kernel estimators \citep{lepski1997optimal,kpotufe2013adaptivity} and, very recently, for local polynomials \citep{chatterjee2024laser}. This paper contributes the following.
\begin{enumerate}
\item A \emph{pointwise} non-asymptotic risk bound for a \emph{network} estimator, \eqref{eq:pointwise-master}. The obstruction this has to overcome is that the observations entering the localized fit are selected by
their own covariates, so the fitted residuals are not the residuals of a fixed design and
the usual basic inequality is unavailable. Section~\ref{sec:rate} resolves this by
conditioning on the covariates and controlling the resulting empirical process at the
critical radius of the local Gaussian complexity of the difference class.
\item A Lepski-type selector for the \emph{neighbourhood size} whose stochastic scale
$\tau_k$ is derived from the metric entropy of the network class rather than from a kernel
variance, together with an oracle inequality over the grid and the resulting adaptation to
the local exponent $\alpha_x$ (Section~\ref{sec:adaptive}).
\item Adaptation to the intrinsic dimension obtained from the \emph{localisation} rather
than from the approximating class. The class of
supports covered is correspondingly wider than in chart-based constructions.
\end{enumerate}

It is important to note that we do not claim that globally trained networks fail to be
locally adaptive in integrated risk --- the work cited above shows otherwise --- but that
the guarantee available for them is of a different kind from the one given here. And we do
not claim that a network is the only sensible model to fit inside the neighbourhood. Whether it improves on a local polynomial fitted in the same neighbourhood is an empirical question, and Section~\ref{sec:simulation} answers it in both directions. The network wins where the local regression function has genuine non-polynomial shape at the resolvable scale, and loses to a local constant fit where it does not. Section~\ref{sec:simulation} also compares the ReLU activation, which the theory covers, with three other smooth activations (sigmoid, SiLU and GELU). All proofs are deferred to the appendices.

\section{The Localized Estimator}
\label{sec:framework}

Consider an i.i.d. sample $\{(X_i, Y_i)\}_{i=1}^n$ drawn from the underlying distribution of $(X, Y)$. For a given query point $x$, we reorder the data according to increasing values of $\|X_i - x\|$. The reordered sequence is denoted by
\[
(X_{(1,n)}(x), Y_{(1,n)}(x)), \ldots, (X_{(n,n)}(x), Y_{(n,n)}(x)),
\]
or simply by $(X_{(1,n)}, Y_{(1,n)}), \ldots, (X_{(n,n)}, Y_{(n,n)})$ if no confusion arises. The point $X_{(k,n)}(x)$ is called the $k$th nearest neighbour ($k$‑NN) of $x$. We keep the constrained class $\F_n=\F(L,\bm p,s,F,K)$ of
\eqref{eq:network-class} throughout. In the least squares approach the regression function $f_0$ satisfies

\[
\mathbb{E}\{(f_0(X) - Y)^2\} = \inf_f \mathbb{E}\{(f(X) - Y)^2\},
\]
where the infimum is taken over all measurable functions $f : \mathbb{R}^d \to \mathbb{R}$. Since the true distribution of $(X, Y)$ is unknown, one can approximate the $L_2$ risk by its empirical counterpart
\[
\frac{1}{n} \sum_{j=1}^{n} |f(X_j) - Y_j|^2.
\]
To avoid overfitting, we restrict the search to a suitable function class $\mathcal{F}_n$ and define the empirical risk minimizer
\[
\hat f_n(\cdot) \in \arg\min_{f \in \mathcal{F}_n} \frac{1}{n} \sum_{j=1}^{n} |f(X_j) - Y_j|^2,
\]
which satisfies, by definition,
\[
\hat f_n \in \mathcal{F}_n, \quad
\frac{1}{n} \sum_{j=1}^{n} |\hat f_n(X_j) - Y_j|^2
= \min_{f \in \mathcal{F}_n} \frac{1}{n} \sum_{j=1}^{n} |f(X_j) - Y_j|^2.
\]
It can be shown that $\hat f_n$ satisfies

\begin{align*}
\int |\hat f_n(x)-f_0(x)|^2\,\mu(dx)
&\le 2\sup_{f\in\mathcal F_n}\left|
\frac1n\sum_{j=1}^n|f(X_j)-Y_j|^2
-\E\{(f(X)-Y)^2\}\right|\\
&\quad+\inf_{f\in\mathcal F_n}\int|f(x)-f_0(x)|^2\,\mu(dx).
\end{align*}

Interested readers are referred to the monograph \cite{Gyorfietal2002} for a detailed discussion on the least squares regression problems for a wide variety of examples. The above equation states that the total estimation error of the empirical risk minimizer $\hat f_n$ decomposes into two parts:
\begin{enumerate}
    \item A \emph{stochastic error} term measuring how far the empirical $L_2$ risk deviates from the true (population) $L_2$ risk uniformly over the class $\mathcal{F}_n$.
    \item An \emph{approximation error} term quantifying how well the best function in $\mathcal{F}_n$ can approximate the true regression function $f_0$.
\end{enumerate}

The class determines both terms of this upper bound. Enlarging the class
can reduce the approximation error, while increasing the supremum in the
stochastic term. Choosing the class therefore balances approximation and
estimation. See \citet{Gyorfietal2002} for analyses based on this decomposition.

In contrast, our approach introduces a localized training mechanism inspired by the $k$-nearest neighbour ($k$-NN) principle. For each query point $x \in [0,1]^d$, we look at the subset $\mathcal{D}_{x} = \{(X_{(i,n)}(x), Y_{(i,n)}(x))\}_{i=1}^k$ and train a deep neural network on $\mathcal{D}_{x}$. Formally, let $N_k(x)$ contain exactly the first $k$ indices
when the observations are ordered by distance to $x$, breaking ties by
their original indices. With this convention, the localized estimator is defined as 
\[
\widehat{f}(x) = \widehat{g}_x(x) \quad\text{where}\quad \widehat{g}_x \in \arg\min_{f \in \mathcal{F}(L,\bm p,s)} \frac{1}{k}\sum_{i\in N_k(x)} \big(Y_i - f(X_i)\big)^2.
\]

As in other local regression methods, two features require attention. First, the
training sample $\{(X_i,Y_i):i\in N_k(x)\}$ is not an i.i.d.\ sample: membership of
$N_k(x)$ is determined by the whole covariate sample, so the summands are dependent and the
empirical measure is not an independent sample from a fixed neighbourhood distribution. Second, the
quantity of interest is the value $\hat f_{x,k}(x)$ of the fitted function at a single
point, not its $L_2(\mu)$ distance from $f_0$, so a uniform bound over $\F_n$ in
$L_2(\mu)$ is not what is required. We handle stochastic dependence by conditioning on
$X_1,\ldots,X_n$, and transfer training error to the query using the fixed
Lipschitz bound on the fitted class. In particular, given the covariates, the index set $N_k(x)$ is fixed and the noise
variables indexed by it are i.i.d.\ $N(0,\sigma^2)$ under the independent
Gaussian noise assumption used for our theory, so the empirical process in the basic
inequality is a Gaussian process and can be localised. The formal treatment of this modified risk structure and its implications for convergence will be presented in the next section.

\section{Rate Optimal Results: H\"{o}lder Smoothness}
\label{sec:rate}

Throughout the theoretical analysis,
\[
Y_i=f_0(X_i)+\sigma W_i,\qquad
W_1,\ldots,W_n\stackrel{\mathrm{i.i.d.}}{\sim}N(0,1),
\]
where $\sigma>0$ is fixed and the noise vector is independent of
$(X_1,\ldots,X_n)$. The pairs $(X_i,Y_i)$ are i.i.d.,
$X_i\in[0,1]^d$, and $\|f_0\|_\infty\le B$ on the support of $X$.
Use the class \eqref{eq:network-class} and the deterministic tie-breaking
rule of Section~\ref{sec:framework}.

Throughout, $B(x,r)$ denotes the closed Euclidean ball.
Fix a query point $x$ and a neighbourhood size $k$. Write
\[
\hat g=\hat f_{x,k}\in\arg\min_{f\in\mathcal F_n}
\frac1k\sum_{i\in N_k(x)}(Y_i-f(X_i))^2,
\qquad \hat f_{n,k}(x)=\hat g(x).
\]
All training values below are evaluations of this \emph{same} fitted
function $\hat g$, not values of the assembled predictor obtained by
refitting separately at each $X_i$. Define
\[
\|h\|_{2,k}^2=\frac1k\sum_{i\in N_k(x)}h(X_i)^2,
\qquad \partial\mathcal F_n=\{f-g:f,g\in\mathcal F_n\}.
\]
Scaling the output weights and output bias by $a\in[0,1]$ preserves every
constraint in \eqref{eq:network-class}. Thus $\mathcal F_n$ and its
difference class are star-shaped about zero.

\subsection{Basic Inequality and Related Results}
\label{sec:basic-facts}

The first step in our analysis is the standard consequence of $\hat g_x(\cdot)$ being a minimiser of the local sum of squares giving us:

\begin{lemma}[Basic inequality]
\label{lem:basic-inequality}
Assume $f_0\in\mathcal F(L,\bm p,s)$. Then
\begin{align}
    \label{eq:basic-inequality}
    \|\hat f_n - f_0\|_{2,k}^2 \leq \frac{2\sigma}{k}\sum_{i\in N_k(x)} W_i\big(\hat f_n(X_i)-f_0(X_i)\big).
\end{align}
\end{lemma}

This inequality is algebraic; it does not require independent selected
observations. Independence enters when we condition on the covariates.
Conditionally on $X_1,\ldots,X_n$ the index set $N_k(x)$ is fixed, since it is determined
by the covariates alone, so the right-hand side is $2\sigma Z(\hat g_x-f_0)$ with
$Z(g)=k^{-1}\sum_{i\in N_k(x)}W_i g(X_i)$ a centred Gaussian process indexed by the
difference class $\partial\mathcal F$. Because $\hat g_x-f_0$ is data dependent, $Z$ must be
controlled uniformly; and because the unknown error $\|\hat g_x-f_0\|_{2,k}$ appears on both
sides, it suffices to control $Z$ over the ball $\{g\in\partial\mathcal F:\|g\|_{2,k}\le\delta\}$
rather than over the whole class. Hence, given $\delta>0$, we define the local Gaussian complexity of $\mathcal F$ around $f_0$ at scale $\delta$ by
$$G_n(\delta; \partial\mathcal{F}) = \E\left[\sup_{g\in\partial\F:\|g\|_{2,k} \leq \delta} \Big|\frac{1}{k}\sum_{i\in N_k(x)} W_i g(X_i)\Big|\right].$$

The quantity $G_n(\delta;\partial\F)$ is the expected largest correlation between the noise
on the neighbourhood and a function in $\partial\F$ whose empirical size on that
neighbourhood is at most $\delta$. It is computed conditionally on $X_1,\ldots,X_n$, so it
depends on the query point $x$ and on $k$ through $N_k(x)$ and the norm $\|\cdot\|_{2,k}$,
and it measures the complexity of the class only where it matters for the fit at $x$:
on the $k$ design points of the neighbourhood, and only near the origin of $\partial\F$.
Two elementary bounds fix its scale. By the Cauchy--Schwarz inequality,
\[
G_n(\delta;\partial\F)\le \delta\,\E\Big[\Big(\frac1k\sum_{i\in N_k(x)}W_i^2\Big)^{1/2}\Big]\le\delta ,
\]
which holds for any class and gives no improvement with $k$. At the other extreme, if
$\partial\F$ were a $D$-dimensional linear space, $G_n(\delta;\partial\F)$ would be of order
$\delta\sqrt{D/k}$. The metric entropy of $\partial\F$ determines where between these two the
network class falls. The following result establishes a monotonicity property of $G_n(\delta;\partial\F)/\delta$ as a function in $\delta$ and defines the critical radius which helps us establish the upper bound on the training sample error rate. 
\begin{lemma}[Monotonicity and the critical radius]
\label{lem:critical-exists}
For $0<\delta\le t$,
\begin{align}
    \label{eq:star-monotone}
    \frac{\delta}{t} G_n(t;\partial\F) \leq G_n(\delta;\partial\F) \Rightarrow \frac{G_n(t;\partial\F)}{t} \leq \frac{G_n(\delta;\partial\F)}{\delta},
\end{align}
so that $\delta\mapsto G_n(\delta;\partial\F)/\delta$ is non-increasing. 
\end{lemma}

Hence, for every $c>0$ the inequality $G_n(\delta;\partial\F)\leq c\delta^2$ has a smallest
positive solution, which we denote by $\delta_n$ and call {\em the critical radius}.
The critical radius is the scale at which the complexity of the class, measured by
$G_n(\delta;\partial\F)$, is balanced by the squared radius $c\delta^2$. Note that the set of
solutions of $G_n(\delta;\partial\F)\le c\delta^2$ is non-empty, because the Cauchy--Schwarz bound
gives $G_n(\delta;\partial\F)\le\delta\le c\delta^2$ for all $\delta\ge1/c$. By
Lemma~\ref{lem:critical-exists} it is an interval unbounded above: writing the inequality as
$G_n(\delta;\partial\F)/\delta\le c\delta$, the left-hand side is non-increasing and the
right-hand side is increasing in $\delta$, so if $\delta$ is a solution then so is every
$\delta'\ge\delta$. That the smallest solution exists and is positive needs two further
facts. First, $\delta\mapsto G_n(\delta;\partial\F)$ is continuous on $(0,\infty)$: it is
non-decreasing, and for $\delta<t$ Lemma~\ref{lem:critical-exists} gives
$G_n(\delta;\partial\F)\le G_n(t;\partial\F)\le (t/\delta)G_n(\delta;\partial\F)$, so the left
and right limits coincide. Hence the solution set is closed and its infimum is attained.
Second, provided $\partial\F$ contains a function that does not vanish identically on
$\{X_i:i\in N_k(x)\}$ (otherwise the problem is trivial), $G_n(1;\partial\F)>0$, and for
$\delta\le1$ the monotonicity gives
$G_n(\delta;\partial\F)/\delta\ge G_n(1;\partial\F)>c\delta$ as soon as
$\delta<\min\{1,G_n(1;\partial\F)/c\}$. Hence no $\delta$ below this threshold solves the
inequality, and $\delta_n>0$. Finally, since every $\delta\ge1/c$ is a solution,
$\delta_n\le1/c$; with $c=1$ this gives $\delta_n\le1$. The critical radius depends
on $c$ and on the covariates; we take $c=1$ throughout, and every bound that involves
$\delta_n$ is conditional on $X_1,\ldots,X_n$ until the expectation over the design is taken
in Section~\ref{sec:train-to-test}. Only star-shapedness of $\partial\F$ is used in the lemma,
and this is the point at which the difference class, rather than $\mathcal F-f_0$, is required.

\subsection{Training Sample Error Rate}
\label{sec:train-error}

\begin{proposition}[Local training error, well-specified case]
\label{prop:training-error}
Suppose $f_0\in\mathcal F_n$ and put
$\kappa_\sigma=\max\{1,16\sigma^2\}$. Under the independent Gaussian
noise model above, for every $t\ge\delta_n$,
\[
\P\big(\|\hat g-f_0\|_{2,k}^2>\kappa_\sigma t\delta_n\mid X\big)
\le e^{-kt\delta_n/2}.
\]
\end{proposition}

Three features of the bound are worth noting. First, its scale. Taking $t=\delta_n$, the local
training error is at most $\kappa_\sigma\delta_n^2$ outside an event of probability
$\exp\{-k\delta_n^2/2\}$. For the network class,
one can show (Proposition~\ref{prop:critical-radius} below) that $k\delta_n^2\lesssim\mathcal C_n$, so this is the
familiar behaviour of least squares with $\mathcal C_n$ effective parameters fitted to $k$
observations: an error of at most order $\max\{1,\sigma^2\}\mathcal C_n/k$. The failure probability is
$e^{-k\delta_n^2/2}$, which is small precisely when $k\delta_n^2$ is large, that is, when the class is rich
relative to the local sample. Larger deviations,
$\kappa_\sigma m\delta_n^2$ with $m\ge1$, have probability at most
$\exp\{-mk\delta_n^2/2\}$ (take $t=m\delta_n$), and integrating this tail gives the following
Corollary.

\begin{corollary}[Expected local training error]
\label{cor:training-expectation}
Under the conditions of Proposition~\ref{prop:training-error},
\begin{align*}
\E[\|\hat g-f_0\|_{2,k}^2\mid X]
&\le\inf_{t\ge\delta_n}\left\{
\kappa_\sigma t\delta_n+\frac{2\kappa_\sigma}{k}e^{-kt\delta_n/2}\right\}\\
&\le\kappa_\sigma\delta_n^2+\frac{2\kappa_\sigma}{k}.
\end{align*}
\end{corollary}

Second, the relevant sample size is $k$, not $n$. The
proposition is the variance half of the bias--variance trade-off that determines $k$: a larger
neighbourhood reduces $\delta_n^2$, and the price, a larger local approximation error, appears
once the assumption $f_0\in\mathcal F$ is dropped in Section~\ref{sec:oracle}. Third, the
statement is conditional on the covariates. Both $\delta_n$ and the norm $\|\cdot\|_{2,k}$
depend on $X_1,\ldots,X_n$ through $N_k(x)$, but given the covariates the noise variables in
the neighbourhood are still i.i.d.\ $N(0,1)$, because the neighbourhood is selected without
reference to the responses. This is what allows a Gaussian argument for an estimator fitted
on a data-dependent subsample. The constant $\kappa_\sigma=\max\{1,16\sigma^2\}$ reflects the two
ways the event can fail to occur: an error already smaller than $\sqrt{t\delta_n}$, which needs
no argument, and a larger error ruled out by the concentration step, which costs the factor
$16\sigma^2$.

\subsection{The critical radius via metric entropy}
\label{sec:entropy-generic}

The critical inequality of Lemma~\ref{lem:critical-exists} is stated in terms of
$G_n(\delta;\partial\F)$, the expected supremum of a Gaussian process indexed by an
infinite-dimensional class, and cannot be evaluated directly. What can be controlled is the
\emph{size} of the class, through its covering numbers. Dudley's chaining argument connects
the two. It approximates each $g$ in the ball $\mathbb B(\delta;\partial\F)$ by a sequence of
progressively finer finite nets, applies a union bound at each resolution, and sums the
resulting maxima. The cost of resolution $t$ is $\sqrt{\log\mathcal N(t)}$, the typical size
of the maximum of $\mathcal N(t)$ Gaussian variables, scaled by $1/\sqrt k$ because the process
averages $k$ independent noise variables. The result converts the critical inequality into a
condition on the metric entropy of $\partial\F$, and for the network class the entropy admits
an explicit upper bound.

\begin{lemma}[Entropy bound for the local Gaussian complexity]
\label{lem:dudley}
Fix $X_1,\ldots,X_n$ and let $C_1>0$. For every $\delta\in(0,C_1]$,
\begin{align}
    \label{eq:dudley-Gn}
    G_n(\delta;\partial\F) \leq \frac{12}{\sqrt{k}} \int_{\delta^2/8C_1^2}^\delta \sqrt{\log \Big(2\mathcal{N}\big(t; \mathbb{B}(\delta;\partial\F), \|\cdot\|_{2,k}\big)\Big)}\mathrm d t + \frac{\delta^2}{2C_1^2},
\end{align}
where $\mathbb{B}(\delta;\partial\F)=\{h\in\partial\F:\|h\|_{2,k}\le\delta\}$ and
$\mathcal N(t;\mathbb B(\delta;\partial\F),\|\cdot\|_{2,k})$ denotes the \emph{internal}
covering number, i.e.\ the smallest number of points of $\mathbb B(\delta;\partial\F)$ whose
$t$-balls cover $\mathbb B(\delta;\partial\F)$. In
particular, with $C_1=1$, any $\delta\in(0,1]$ satisfying
$$\frac{12}{\sqrt{k}} \int_{\delta^2/8}^\delta \sqrt{\log \Big(2\mathcal{N}\big(t; \mathbb{B}(\delta;\partial\F), \|\cdot\|_{2,k}\big)\Big)}\mathrm d t \leq \frac{\delta^2}{2}$$
satisfies the critical inequality $G_n(\delta;\partial\F)\leq \delta^2$.
\end{lemma}

The two terms in \eqref{eq:dudley-Gn} correspond to the two stages of the argument. The
integral is the chaining bound, taken only over resolutions between (a constant multiple of)
$\delta^2/(2C_1^2)$ and
$\delta$. Resolutions coarser than $\delta$ are irrelevant, since every element of the ball
already lies within $\delta$ of the origin. Resolutions finer than $\delta^2/(2C_1^2)$ are not
needed, because the residual approximation error at that scale is controlled directly by the
Cauchy--Schwarz inequality, which contributes the second term. Truncating the chain there is
what makes the bound compatible with the critical inequality: the residual is itself of
order $\delta^2$, so it can be absorbed on the right-hand side of
$G_n(\delta;\partial\F)\le\delta^2$ at the cost of halving the room left for the integral. This is
the second statement of the lemma.

\subsection{Oracle inequality}
\label{sec:oracle}
The argument of Section~\ref{sec:train-error} used $f_0$ in two places: as the competitor in
the basic inequality, and as the centre of the class of possible errors. Neither use survives
when $f_0\notin\mathcal F(L,\bm p,s)$, but both can be transferred to an arbitrary
$\tilde f\in\mathcal F(L,\bm p,s)$. Comparing the local sum of squares of $\hat g_x$ with that
of $\tilde f$ instead of $f_0$ gives a basic inequality with an additional approximation term,
\[
\|\hat g_x-f_0\|_{2,k}^2\le\|\tilde f-f_0\|_{2,k}^2
+\frac{2\sigma}{k}\sum_{i\in N_k(x)}W_i\big(\hat g_x(X_i)-\tilde f(X_i)\big),
\]
and the noise is now paired with $\hat g_x-\tilde f$. Since both $\hat g_x$ and $\tilde f$ lie
in the class, this difference belongs to $\partial\F$. The empirical process is therefore
indexed by the same difference class as before, and the high-probability event constructed in
the proof of Proposition~\ref{prop:training-error} controls it without change. Because that
event does not depend on $\tilde f$, the resulting bound holds simultaneously for every
competitor, and the infimum over $\tilde f$ may be taken afterwards.

\begin{theorem}[Oracle inequality for the localized fit]
\label{thm:oracle}
Let $\delta_n$ be the critical radius of Lemma~\ref{lem:critical-exists} for
$\partial\mathcal F(L,\bm p,s)$, and let $\hat\Delta=\hat f_n-f_0$. Then, conditionally on
$X_1,\ldots,X_n$, for every $t\ge\delta_n$, with $Q_\sigma=\max\{2,36\sigma^2\}$,
\begin{align}
\label{eq:oracle-tail}
\|\hat \Delta\|_{2,k}^2 \leq 2 \inf_{f\in\mathcal{F}(L,\bm p, s)}\|f-f_0\|_{2,k}^2 + Q_\sigma t\delta_n
\end{align}
with probability at least $1-\exp\{-kt\delta_n/2\}$, and
\begin{align}
\label{eq:oracle-expectation}
\E\left[\frac{1}{k}\sum_{i\in N_k(x)} \big(\hat f_n(X_i) - f_0(X_i)\big)^2\ \Big|\ X_1,\ldots,X_n\right] \leq 2\inf_{f\in\mathcal{F}(L,\bm p, s)}\|f-f_0\|_{2,k}^2 + Q_\sigma\delta_n^2 + \frac{2Q_\sigma}{k}.
\end{align}
\end{theorem}

The theorem is a bias--variance decomposition on the neighbourhood. The first term is the
error of the best approximation to $f_0$ from the network class, measured only at the $k$
design points in $N_k(x)$. The second is the estimation error of Proposition~\ref{prop:training-error},
unchanged. Three features of the decomposition are worth recording.

First, the price of dropping the assumption $f_0\in\mathcal F(L,\bm p,s)$ is the approximation error. The stochastic term does not depend on $f_0$. The critical radius is defined through
$\partial\F$, a property of the class alone, so misspecification affects only the
approximation term. 

Second, the approximation term is \emph{local}. The class needs to approximate $f_0$ only on
the $k$ points near $x$, not over the whole support of $X$. This is where localisation pays
on the approximation side. The two terms also move in opposite directions as $k$ grows: the neighbourhood
widens, so the approximation error increases, while $\delta_n^2$ decreases. The balance
between them determines the choice of $k$ in Section~\ref{sec:train-to-test} and its adaptive choice
in Section~\ref{sec:adaptive}.

Third, both statements are conditional on $X_1,\ldots,X_n$. The critical radius, the norm
$\|\cdot\|_{2,k}$ and the index set $N_k(x)$ all depend on the covariates, and it is only in
Section~\ref{sec:train-to-test}, once $\delta_n$ has been bounded by a deterministic quantity,
that the expectation over the design can be taken.

\subsection{Expected local approximation error}
\label{sec:approx}

The oracle inequality leaves the approximation term
$\inf_{f\in\F_n}k^{-1}\sum_{i\in N_k(x)}\{f(X_i)-f_0(X_i)\}^2$ to be controlled. Because
$\F_n$ contains the constant functions bounded by $B$, this term is bounded by the
variation of $f_0$ over the neighbourhood, and hence by the neighbourhood radius
\[
R_{n,k}(x) = \max_{i\in N_k(x)} \|X_i-x\| = \|X_{(k,n)}(x)-x\| .
\]

\begin{proposition}[Local approximation error]
\label{prop:approx}
Suppose $f_0$ is locally $\alpha$-H\"older continuous at $x$ with exponent
$0<\alpha\le1$ and constant $L_0$ on a ball $B(x,r_0)$, that $\|f_0\|_\infty\le B$, that
$\F_n$ contains every constant function bounded by $B$, and that the assumptions of
Lemma~\ref{lem:local_mass} below hold with $k\le\frac12nc_0r_0^{d^*}$. Then
\begin{align*}
    \E\left[\inf_{f\in\mathcal F_n}\frac1k\sum_{i\in N_k(x)}\{f(X_i)-f_0(X_i)\}^2\right]
  & \le L_0^2\,\E\big[R_{n,k}(x)^{2\alpha}\big]+4B^2\,\P\big(R_{n,k}(x)>r_0\big)\\
  & \le L_0^2\,\E\big[R_{n,k}(x)^{2\alpha}\big]+4B^2\exp\Big\{-\frac18nc_0r_0^{d^*}\Big\},
\end{align*}
the second bound using Lemma~\ref{lem:knn_radius_tail} below. The last term is
exponentially small in $n$; in particular it is at most a constant times
$n^{-2\alpha/d^*}\le(k/n)^{2\alpha/d^*}$, which is the order of the first term established in
Lemma~\ref{lem:expected_knn_radius}.
\end{proposition}

It remains to bound $\E[R_{n,k}(x)^{2\alpha}]$, and it is here that the intrinsic
dimension enters. Let $M\subset \mathbb{R}^d$ be a compact $d^*$-dimensional smooth
Riemannian manifold without boundary. Suppose that the distribution of $X$ is supported on
$M$ and has density $p_X$ with respect to the $d^*$-dimensional Riemannian volume measure
on $M$, and that
\[
0<p_{X,0}\leq p_X(z)\leq p_{X,1}<\infty, \qquad z\in M.
\]
We also assume that $M$ has condition number $1/\tau$ for some $\tau>0$.

The following lemma gives a uniform lower bound on the probability mass of sufficiently
small Euclidean neighbourhoods of $x$.

\begin{lemma}[Local mass of the manifold]
\label{lem:local_mass}
There exist constants $r_0>0$ and $0<c_0\le c_1<\infty$, depending only on
$d^*$, $\tau$, $p_{X,0}$, $p_{X,1}$ and the geometry of $M$, such that
\[
c_1 r^{d^*}\ \geq\ \mathbb{P}\bigl(\|X-x\|\leq r\bigr)
\geq c_0 r^{d^*},
\qquad
0<r\leq r_0,
\]
for every $x\in M$. 
\end{lemma}

The lower bound is used for the upper-bound analysis, the upper bound
only in the proof of the minimax lower bound, Proposition~\ref{prop:lower}. We now use the preceding local mass bound to control the tail of the
$k$-nearest-neighbour radius.

\begin{lemma}[Tail bound for the $k$-nearest-neighbour radius]
\label{lem:knn_radius_tail}
Suppose the assumptions of Lemma~\ref{lem:local_mass} hold. 
Suppose that
\(
k\leq \frac12 n c_0r_0^{d^*}.
\)
Then, for every
\(
r\in [r_{n,k},r_0]
\) with \(
r_{n,k}
=
\left(\frac{2k}{c_0n}\right)^{1/d^*},
\)
we have
\[
\mathbb{P}\bigl(R_{n,k}(x)>r\bigr)
\leq
\exp\left\{-\frac{c_0}{8}nr^{d^*}\right\}.
\]
\end{lemma}

We next convert this tail bound into a bound on the expected
$2\alpha$-th power of the nearest-neighbour radius.

\begin{lemma}[Expected $k$-nearest-neighbour radius]
\label{lem:expected_knn_radius}
Suppose the assumptions of Lemma~\ref{lem:local_mass} hold and
$0<\alpha\leq 1$. Assume that
\(
k\leq \frac12 nc_0r_0^{d^*}.
\)
Then, uniformly for $x\in M$,
\[
\mathbb{E}\left[R_{n,k}(x)^{2\alpha}\right]
\lesssim
\left(\frac{k}{n}\right)^{2\alpha/d^*}.
\]
\end{lemma}

Proposition~\ref{prop:approx} and Lemma~\ref{lem:expected_knn_radius} combine at once
into the form in which the approximation error is used below.

\begin{corollary}[Approximation error at the intrinsic rate]
\label{cor:approx-rate}
Under the assumptions of Proposition~\ref{prop:approx} and
Lemma~\ref{lem:expected_knn_radius},
\[
\mathbb{E}\left[\inf_{f\in \F}\frac{1}{k}\sum_{i\in N_k(x)}\bigl\{f(X_i)-f_0(X_i)\bigr\}^2\right]
\lesssim \left(\frac{k}{n}\right)^{2\alpha/d^*}.
\]
\end{corollary}

The corollary shows that the approximation error of the localized
regression problem is governed by the radius of the $k$-nearest-neighbour ball. In
particular, the ambient dimension $d$ does not enter this bound directly; the effective
dimension is the intrinsic dimension $d^*$ of the support of $X$. Also, these upper bounds extend to other compact supports satisfying the lower
mass condition, including regular domains with boundary. For example,
uniform design on $[0,1]$ has the required lower mass bound at the endpoints
as well as in the interior. This covers the design in the motivating example
for the upper-bound argument, although it is not a boundaryless manifold.

\subsection{From Local Training Error to Pointwise Testing Error}
\label{sec:train-to-test}

We now transfer the error of the localized estimator on its training
neighbourhood to its error at the query point $x$. Define
\[
T_{n,k}(x)
=
\frac{1}{k}
\sum_{i\in N_k(x)}
\bigl\{
\widehat f_n(X_i)-f_0(X_i)
\bigr\}^2
\]
and recall that
\(
R_{n,k}(x)
=
\max_{i\in N_k(x)}\|X_i-x\|.
\)
We impose the following regularity on the fitted class itself.

\begin{assumption}[Smoothness of the network class]
\label{ass:class-holder}
For a fixed $0<\alpha\le1$, there are constants $L_\F,L_0<\infty$,
independent of $n$, such that for all $u,v\in M$,
\[
|f(u)-f(v)|\le L_\F\|u-v\|^\alpha\quad(f\in\mathcal F_n),\qquad
|f_0(u)-f_0(v)|\le L_0\|u-v\|^\alpha.
\]
For \eqref{eq:network-class}, the first condition holds simultaneously for
all $0<\alpha\le1$ with $L_\F=K\max\{1,\sqrt d\}$.
\end{assumption}

\begin{proposition}[Pointwise risk of the localized estimator]
\label{prop:pointwise}
Under
Assumption~\ref{ass:class-holder}, let $L_H=L_\F+L_0$ and $x\in M$. Then the deterministic
inequality
\[
\big|\widehat f_n(x)-f_0(x)\big|^2 \leq 2T_{n,k}(x) + 2L_H^2R_{n,k}(x)^{2\alpha}
\]
holds, and combining it with Theorem~\ref{thm:oracle} and
Corollary~\ref{cor:approx-rate} gives
\begin{align}
\label{eq:pointwise-master}
\mathbb{E}\big[\big|\widehat f_n(x)-f_0(x)\big|^2\big]
\lesssim \left(\frac{k}{n}\right)^{2\alpha/d^*} + Q_\sigma\,\E\big[\delta_n^2\big] + \frac{2Q_\sigma}{k},
\end{align}
where $Q_\sigma=\max\{2,36\sigma^2\}$ is the constant of Theorem~\ref{thm:oracle}. %
\end{proposition}

The two terms have the interpretation the construction suggests: the first is the
squared local approximation error, arising from the radius of the $k$-nearest
neighbourhood, and the second is the stochastic estimation error of fitting the network on
the local sample.
In particular, if the localized Gaussian complexity satisfies
\[
\delta_n^2\lesssim \frac{C_n}{k},
\]
where, by Section~\ref{sec:entropy-Fn} below, $C_n$ may be taken to be the
\emph{deterministic} quantity $(s+1)\log(1+n(L+1)V^2/(s+1))$ --- which is what permits the
expectation over the covariates taken in \eqref{eq:pointwise-master}, since $\delta_n$ is
itself covariate dependent ---
then
\[
\mathbb{E}\left[\bigl|\widehat f_n(x)-f_0(x)\bigr|^2\right]
\lesssim \left(\frac{k}{n}\right)^{2\alpha/d^*} + \frac{C_n+1}{k},
\]
and balancing the two terms by choosing
\[
k\asymp (C_n+1)^{d^*/(2\alpha+d^*)} n^{2\alpha/(2\alpha+d^*)}
\]
yields
\[
\mathbb{E}\left[\bigl|\widehat f_n(x)-f_0(x)\bigr|^2\right]
\lesssim \left(\frac{C_n+1}{n}\right)^{2\alpha/(2\alpha+d^*)} .
\]
The quantity $C_n$ cannot be bounded uniformly in $n$: since $s\ge0$, $L\ge0$ and $V\ge1$,
$C_n\ge\log(1+n)$. The smallest possible order is therefore $C_n\asymp\log n$, attained when
$s$, $L$ and $V$ do not grow with $n$, and the best rate this bound yields is
$(\log n/n)^{2\alpha/(2\alpha+d^*)}$. This logarithmic factor is the one that appears in
Theorem~\ref{thm:minimax} below.

\subsubsection{The critical radius for the network class}
\label{sec:entropy-Fn}

Here we evaluate $\delta_n$ for the sparse ReLU class. Combining the covering
number bound of \citet{schmidthieber2019manifold} for $\mathcal F_n$ with Lemma~\ref{lem:dudley} gives a
critical inequality that can be solved in closed form.

\begin{proposition}[Critical radius for the sparse ReLU class]
\label{prop:critical-radius}
Let $\mathcal F_n=\mathcal F(L,\bm p,s)$ with $\max_j\|W_j\|_\infty\le1$ and
$V=\prod_{\ell}(p_\ell+1)$. Then the critical inequality
$G_n(\delta;\partial\mathcal F_n)\le\delta^2$ is implied by
\begin{align}
\label{eq:critical-ineq}
\delta^2 \geq \frac{2048(s+1)}{k}\log\left(\frac{128(L+1)V^2}{\delta^2}\right),
\end{align}
whose smallest solution is
\begin{align}
\label{eq:delta-lambertW}
\bar\delta^2 = \frac{2048(s+1)}{k}\,W\!\left(\frac{k(L+1)V^2}{16(s+1)}\right),
\end{align}
with $W$ the Lambert $W$ function. The critical radius satisfies
$\delta_n\le\bar\delta$. Using $W(z)\le\log(1+z)$ for $z\ge0$,
\[
\delta_n^2\le\bar\delta^2 \lesssim \frac{s+1}{k}\log\left(1+\frac{k(L+1)V^2}{s+1}\right) =: \frac{\mathcal C_n}{k}.
\]
\end{proposition}

The quantity $\bar\delta$ is an upper bound for the critical radius, not the critical radius
itself: \eqref{eq:critical-ineq} is only a sufficient condition for $G_n(\delta;\partial\F_n)\le\delta^2$,
so its smallest solution can only be larger than the smallest solution of the critical
inequality. This is all that is needed, since every bound of Sections~\ref{sec:train-error}
and~\ref{sec:oracle} remains valid when $\delta_n$ is replaced by any larger radius at which the
critical inequality holds (see the proof of Lemma~\ref{lem:computable-radius}). Unlike $\delta_n$,
$\bar\delta$ does not depend on the covariates.

Substituting into Theorem~\ref{thm:oracle} and Proposition~\ref{prop:pointwise},
\[
\mathbb{E}\big[T_{n,k}(x)\big] \lesssim \left(\frac{k}{n}\right)^{2\alpha/d^*} + \frac{\mathcal C_n+1}{k},
\qquad
\mathbb{E}\left[\bigl|\widehat f_n(x)-f_0(x)\bigr|^2\right]
\lesssim \left(\frac{k}{n}\right)^{2\alpha/d^*} + \frac{\mathcal C_n+1}{k}.
\]

\begin{remark}[Removing the circularity in $\mathcal C_n$]
\label{rem:Cn}
The quantity $\mathcal C_n=(s+1)\log(1+k(L+1)V^2/(s+1))$ depends on $k$, while $k$ is
chosen below as a function of $\mathcal C_n$. Since $\mathcal C_n$ is increasing in $k$ and
$k\le n$, one may replace it throughout by the deterministic upper bound
$(s+1)\log(1+n(L+1)V^2/(s+1))$, which is of the same order whenever $k\asymp n^{\theta}$ for
a fixed $\theta\in(0,1)$. We do so tacitly in what follows.
\end{remark}

\subsection{Minimax lower bound}
\label{sec:lower}

For fixed $L_0,B>0$, use the bounded target class
\[
\mathcal H^\alpha(L_0,B)=\left\{f:M\to\mathbb R:
\|f\|_\infty\le B,\quad
|f(u)-f(v)|\le L_0\|u-v\|^\alpha\ \text{for all }u,v\in M\right\}.
\]
The amplitude restriction is part of the class over which we take the
supremum, not only an assumption on one chosen target.

\begin{proposition}[Pointwise minimax lower bound]
\label{prop:lower}
Fix $x\in M$, $0<\alpha\le1$, $L_0,B>0$ and $\sigma>0$, and assume the
geometric and density conditions above. Then, for sufficiently large $n$,
\[
\inf_{\tilde f_n}\sup_{f_0\in\mathcal H^\alpha(L_0,B)}
\E_{f_0}|\tilde f_n(x)-f_0(x)|^2
\gtrsim n^{-2\alpha/(2\alpha+d^*)}.
\]
\end{proposition}

The proof places a H\"older bump at $x$. Its height is of order $h^\alpha$
and its support has probability at most a constant times $h^{d^*}$, so
its Kullback--Leibler separation from zero is of order
$nh^{2\alpha+d^*}$. This gives the intrinsic-dimensional lower rate.

\begin{theorem}[Minimax rate-optimality up to logarithmic factors]
\label{thm:minimax}
Fix $F,K,B,L_0,\sigma$ as above, with $0<B\le1$, $F\ge B$ and $K,L_0>0$.
Use the class \eqref{eq:network-class} and the geometric design assumptions.
Suppose $0<\alpha\le1$ and $\mathcal C_n+1=O((\log n)^C)$ for some $C>0$;
for example, it suffices that $s,L,\log V$ grow polylogarithmically.
Choose
$k\asymp(\mathcal C_n+1)^{d^*/(2\alpha+d^*)}
n^{2\alpha/(2\alpha+d^*)}$.
Then, for sufficiently large $n$,
\[
\sup_{f_0\in\mathcal H^\alpha(L_0,B)}
\E_{f_0}|\hat f_{n,k}(x)-f_0(x)|^2
\lesssim n^{-2\alpha/(2\alpha+d^*)}
(\log n)^{2\alpha C/(2\alpha+d^*)}.
\]
Together with Proposition~\ref{prop:lower}, this proves pointwise minimax
rate-optimality over $\mathcal H^\alpha(L_0,B)$ up to logarithmic factors.
\end{theorem}

\begin{remark}
\label{rem:higher-smoothness}
The intrinsic dimension enters through the local mass condition and the
nearest-neighbour radius. Moreover, the analysis is restricted to $0<\alpha\le1$. For higher smoothness,
one needs derivative-based H\"older classes and an appropriate local
polynomial or network approximation result. 
\end{remark}

\subsection{How should $k$ be chosen in practice?}
\label{sec:question}

Theorem~\ref{thm:minimax} is an oracle statement. The neighbourhood size it prescribes,
\[
k\asymp (\mathcal C_n+1)^{d^*/(2\alpha+d^*)}\,n^{2\alpha/(2\alpha+d^*)},
\]
depends on the smoothness exponent $\alpha$ of $f_0$ and on the intrinsic dimension $d^*$
of the support of $X$, neither of which is known to the statistician. The dependence is not
a minor calibration issue: the exponent of $n$ varies from $0$ to $2/(2+d^*)$ as $\alpha$ ranges
over $(0,1]$, so a misspecified $\alpha$ produces a neighbourhood of the wrong
order and a suboptimal rate. Moreover, a correct \emph{global}
value of $\alpha$ need not help, because the difficulty of the problem varies across the
covariate space. 

For the Doppler function \eqref{eq:doppler} the local H\"older exponent is the same
everywhere, since $f_0$ is globally Lipschitz (see the discussion following \eqref{eq:doppler}), what changes is the
local constant. The derivative of $f_0$ reaches about $190$ in absolute value on $[0,0.15]$, is
at most about $55$ on $[0.15,1]$ and at most about $10$ on $[0.5,1]$. At the sample sizes of
interest the neighbourhoods near the origin therefore see a function that varies on a much
finer scale than elsewhere, and the bias term, which is proportional to the local constant
$L_x$ (see \eqref{eq:bias-bk} below), is correspondingly larger. The balance defining the oracle $k$ should therefore be struck separately at
each query point. Since the estimator is already computed one query point at a time, nothing
in the construction prevents this: the neighbourhood size may be chosen as a function of
$x$.

This is the problem taken up in the next section. The key structural observation is visible
in \eqref{eq:pointwise-master}: the approximation term $(k/n)^{2\alpha/d^*}$ is
increasing in $k$ while the stochastic term $\mathcal C_n/k$ is decreasing, and --- unlike
the first --- the second is \emph{computable}, because it is determined by the architecture
$(L,\bm p,s)$ and not by $f_0$. A selection rule may therefore compare estimates across
neighbourhood sizes on the scale of the computable term and stop when they begin to
disagree by more than that scale. This is Lepski's principle.

\section{Locally adaptive estimation}
\label{sec:adaptive}

The results in the preceding section assume that the smoothness of the
regression function is known and hence that the number of nearest
neighbours $k$ can be chosen according to the corresponding oracle
balance. In this section, we develop a data-driven choice of $k$ which
adapts to the local smoothness of the regression function.

The main idea is to construct, for each query point $x$, a collection of
localized neural-network estimators corresponding to different values of
$k$, and then select the largest neighbourhood whose estimate is
statistically indistinguishable from all estimators based on smaller
neighbourhoods. This is a Lepski-type procedure. The resulting estimator
has a pointwise risk which adapts to the unknown local Hölder exponent
$\alpha_x$.

The argument has four steps, and the subsections follow them. We first fix the notion
of local smoothness to which the procedure should adapt (Section~\ref{sec:local-smoothness})
and the finite collection of candidate estimators among which it chooses
(Section~\ref{sec:collection}). We then show that, with high probability and
\emph{simultaneously} over the collection, the error of each candidate splits into a bias
term $b_k(x)$ that increases with $k$ and a stochastic term $\tau_k$ that decreases with $k$
and, crucially, can be computed from the data (Section~\ref{sec:stochastic}). Lepski's rule
(Section~\ref{sec:lepski}) uses only the computable term. Section~\ref{sec:oracle-adaptive}
shows that it nevertheless selects a neighbourhood at least as large as the unknown optimal
one, and that the resulting estimator performs as well as the best candidate up to a
constant. Section~\ref{sec:adaptation} translates this into a rate. Proofs are in
Appendix~\ref{app:proofs-adaptive}.

Throughout this section, we keep the geometric assumptions of
Section~\ref{sec:approx}. The covariate $X$ is supported on
$\mathcal M\subset[0,1]^d$ (the set $M$ of Section~\ref{sec:approx}), a compact smooth embedded
$d^*$-dimensional manifold without boundary, with reach at least $\tau>0$
and $1\le d^*<d$, and $c_0,r_0$ are the constants of
Lemma~\ref{lem:local_mass}. The inclusion in $[0,1]^d$ is needed because
the constraints in \eqref{eq:network-class}, and hence the entropy bound and
the constant $L_\F$, are imposed on $[0,1]^d$. Only the lower mass bound is
used below, so, as noted after Corollary~\ref{cor:approx-rate}, compact
supports with a regular boundary are also covered. Let $p_X$ denote the density
of $X$ with respect to the Riemannian volume measure on $\mathcal M$ and
suppose that
\[
0<p_{X,0}
\le p_X(x)
\le p_{X,1}<\infty,
\qquad x\in\mathcal M.
\]
We consider the same regression model
\[
Y_i=f_0(X_i)+\sigma W_i,
\qquad
W_1,\ldots,W_n\stackrel{\mathrm{i.i.d.}}{\sim}N(0,1),
\]
with $(W_1,\ldots,W_n)$ independent of $(X_1,\ldots,X_n)$.
For every integer $k$, let
\[
N_k(x)
=
\left\{
i:X_i\text{ is among the $k$ nearest neighbours of }x
\right\}
\]
and define
\[
R_k(x)
=
\max_{i\in N_k(x)}\|X_i-x\|
=
\|X_{(k,n)}(x)-x\|.
\]
The localized neural-network estimator based on $k$ neighbours is defined
by
\[
\widehat f_k(x)
=
\widehat f_{x,k}(x),
\qquad
\widehat f_{x,k}
\in
\arg\min_{f\in\mathcal F_n}
\frac{1}{k}
\sum_{i\in N_k(x)}
\{Y_i-f(X_i)\}^2.
\]

\subsection{Local smoothness}
\label{sec:local-smoothness}

Theorem~\ref{thm:minimax} uses one exponent $\alpha$. We now allow the
available local H\"older exponent to vary with the query point. For example,
$f(z)=|z-x_0|^\beta$, $0<\beta<1$, has exponent $\beta$ at $x_0$ and is
locally Lipschitz at every $x\ne x_0$. Its local Lipschitz constant
$\beta|x-x_0|^{\beta-1}$ is unbounded as $x\to x_0$, so the uniform bounds
$L_x\le L_0$ and $r_x\ge r_{\min}$ in Definition~\ref{def:local-holder}
hold with $\alpha_x=1$ only on sets bounded away from $x_0$; with
$\alpha_x=\beta$ they hold everywhere, with $L_0=1$. The Doppler example has
different finite-scale variation but remains Lipschitz throughout.

\begin{definition}
\label{def:local-holder}
Let $x\in\mathcal M$. We say that $f_0$ is locally $\alpha_x$-Hölder
smooth at $x$ with constant $L_x$ if there exists $r_x>0$ such that
\[
|f_0(y)-f_0(x)|
\le
L_x\|y-x\|^{\alpha_x},
\qquad
y\in\mathcal M\cap B(x,r_x),
\]
for some $\alpha_x\in(0,1]$, with $L_x\leq L_0$ for some $L_0>0$. We assume in addition that $r_x\ge r_{\min}>0$ uniformly in $x$.
\end{definition}

Thus both the exponent $\alpha_x$ and the Hölder constant $L_x$ may
depend on the query point $x$.
We continue to assume that there exists $B>0$ such that
$\|f_0\|_\infty\le B$
and that the neural-network class $\mathcal F_n$ contains all constant
functions $f_a(\cdot)\equiv a$, $a\in[-B,B]$. The class remains \eqref{eq:network-class}, with the same fixed Lipschitz
bound $K$ at every candidate $k$. Hence its members are uniformly
$\alpha_x$-H\"older with the fixed constant $L_\F=K\max\{1,\sqrt d\}$.
We assume the noise level $\sigma$ is known for the selector analysed below.

\begin{remark}
\label{rem:absorb-event}
Let $r_* = \min\{r_0,r_{\min}\}$. Under the upper grid restriction below,
\[
\P\big(\exists k\in\mathcal K_n:R_k(x)>r_x\big)
\le\P(R_{k_J}(x)>r_*)\le e^{-nc_0r_*^{d^*}/8}.
\]
Proposition~\ref{prop:simultaneous} controls a stronger radius event directly,
so no separate exceptional event needs to be added there.
\end{remark}

\subsection{A collection of localized estimators}
\label{sec:collection}

Adaptation requires a set of candidates to choose among. Because the relevant
quantities change by a constant factor when $k$ doubles, a geometric grid is enough. Let
\[
\mathcal K_n
=
\{k_0,k_1,\ldots,k_J\},
\qquad
k_j=2^j k_0,
\quad
j=0,\ldots,J,
\]
be a geometric grid of candidate neighbourhood sizes. We assume that
\[
k_0\longrightarrow\infty,
\qquad
\frac{k_J}{n}\longrightarrow0,
\]
and, more precisely, that
\begin{equation}
k_0\ge4\log\Big(\frac{|\mathcal K_n|}{\eta}\Big),
\qquad
k_J\le\frac12\,c_0\,n\min\{r_0,r_{\min}\}^{d^*},
\label{eq:grid-conditions}
\end{equation}
where $c_0,r_0$ are the constants of Lemma~\ref{lem:local_mass} and $\eta\in(0,1)$ is a
confidence level fixed in advance. The first condition lets a union bound over the grid be
absorbed into the nearest-neighbour radius; the second keeps every neighbourhood inside the
region where Lemma~\ref{lem:knn_radius_tail} and Definition~\ref{def:local-holder} apply.
A concrete choice, for sufficiently large $n$, is
\begin{equation}
\label{eq:implementable-grid}
k_0=\lceil8\log(n/\eta)\rceil,\qquad
J=\left\lfloor\log_2\frac{n}{k_0\log n}\right\rfloor,
\qquad k_j=2^jk_0.
\end{equation}
This choice does not use $d^*$ or any geometric constant. It satisfies
\eqref{eq:grid-conditions} eventually for every fixed admissible design,
for fixed $\eta$ or for $\eta=n^{-1}$. Under the growth assumptions of
Corollary~\ref{cor:adaptive-OP}, the lower end is logarithmic and
$k_J\asymp n/\log n$. For every fixed $\alpha_x>0$,
$b_{k_0}/\tau_{k_0}\to0$ and $b_{k_J}/\tau_{k_J}\to\infty$;
also $k_0\le k_x^\circ\le k_J$ eventually, where $k_x^\circ$ is the balance point defined in Section~\ref{sec:adaptation}. These assertions follow from
the explicit upper bound for $\tau_k^2$ and the lower bound
$k\bar\delta_k^2\ge2048W(k/16)$ proved in Appendix~\ref{app:proofs-adaptive}.
Thus the window conditions below hold eventually at each fixed query with
a fixed positive exponent. This is not a uniform statement over exponents
that tend to zero with $n$.

The cardinality of the grid satisfies $|\mathcal K_n|=J+1\lesssim \log n$.
For every $k\in\mathcal K_n$, $\widehat f_k(x)=\widehat f_{x,k}(x)$ is defined as above.
The purpose of the geometric grid is only to ensure that the oracle
choice of $k$ can be approximated within a constant factor. In
particular, restricting the procedure to $\mathcal K_n$ does not affect
the polynomial convergence rate.

\subsection{Pointwise stochastic bounds for the localized estimators}
\label{sec:stochastic}

A selection rule compares estimators with one another, so it needs error bounds that
hold for all candidates at once, and it can act only on quantities it can compute. The
first requirement is met by a union bound over the grid, which costs a factor
$\log(|\mathcal K_n|/\eta)$. The second requires the critical radius of
Section~\ref{sec:rate}, which is defined through an expected supremum and is not itself
computable, to be replaced by an explicit upper bound. The following lemma shows that this
is harmless.

Recall that the relevant difference class is
$\partial\mathcal F_n=\{f-g:f,g\in\mathcal F_n\}$, and let $\delta_k$ denote the critical
radius associated with a local sample of size $k$, that is, the critical radius of
Lemma~\ref{lem:critical-exists} with $N_k(x)$ as the index set.

\begin{lemma}[A computable critical radius]
\label{lem:computable-radius}
Let $\mathcal F_n$ be the class in \eqref{eq:network-class}, put
$V=\prod_{\ell=0}^{L+1}(p_\ell+1)$, and define
\[
\bar\delta_k^2:=\frac{2048(s+1)}{k}\,W\!\left(\frac{k(L+1)V^2}{16(s+1)}\right),
\]
with $W$ the Lambert $W$ function. Then $\bar\delta_k\ge\delta_k$, and
Proposition~\ref{prop:training-error} and Theorem~\ref{thm:oracle} remain valid with
$\delta_n$ replaced by $\bar\delta_k$ throughout. Moreover
\begin{equation}
\label{eq:delta-k}
\bar\delta_k^2
\le
\frac{2048(s+1)}{k}
\log\left(
1+
\frac{k(L+1)V^2}{s+1}
\right).
\end{equation}
\end{lemma}

The point of the lemma is that $\bar\delta_k$ depends only on the architecture and on
$k$. Any radius at which the critical inequality holds can play the role of the critical
radius in the arguments of Section~\ref{sec:rate}; the smallest such radius gives the best
constant, but a slightly larger one that can be written down is what a procedure needs.

With $\bar\delta_k$ in hand, set
\[
 t_k=\max\Big\{\bar\delta_k,\ \frac{2}{k\bar\delta_k}\log\Big(\frac{|\mathcal K_n|}{\eta}\Big)\Big\},
\]
which is the smallest $t\ge\bar\delta_k$ for which the exceptional probability
$\exp(-tk\bar\delta_k/2)$ of Theorem~\ref{thm:oracle} is at most $\eta/|\mathcal K_n|$, and
define the stochastic scale
\begin{equation}
\label{eq:tau-scale}
\tau_k=6\sqrt{\kappa_\sigma\,t_k\bar\delta_k}
=6\sqrt{\kappa_\sigma}\,\max\Big\{\bar\delta_k^2,\ \frac2k\log\Big(\frac{|\mathcal K_n|}{\eta}\Big)\Big\}^{1/2},
\end{equation}
with $\kappa_\sigma=\max\{1,16\sigma^2\}$ as in Proposition~\ref{prop:training-error}.
It is the high-probability size of the stochastic error of $\hat f_k(x)$, uniformly over the
grid. The next result combines it with the local approximation error.

\begin{proposition}[Simultaneous pointwise bound over the grid]
\label{prop:simultaneous}
Let $f_0$ satisfy Definition~\ref{def:local-holder} at $x$, and let the grid satisfy
\eqref{eq:grid-conditions}. Define
\begin{equation}
\label{eq:bias-bk}
 b_k(x)=C_{\mathcal M}\,(L_\F+3L_x)\Big(\frac kn\Big)^{\alpha_x/d^*},
\qquad C_{\mathcal M}=\max\{1,(2/c_0)^{1/d^*}\}.
\end{equation}
Then there is an event $\mathcal A_x$ with $\P(\mathcal A_x)\ge1-2\eta$ on which,
simultaneously for every $k\in\mathcal K_n$,
\begin{equation}
\label{eq:simultaneous}
\big|\hat f_k(x)-f_0(x)\big|\le b_k(x)+\tau_k .
\end{equation}
\end{proposition}

The two terms in \eqref{eq:simultaneous} are the ones anticipated in
Section~\ref{sec:question}. The bias $b_k(x)$ grows with $k$ at a rate set by the unknown
pair $(\alpha_x,d^*)$ and is not available to the statistician. The stochastic term $\tau_k$
is available. Lepski's rule needs, in addition, that $\tau_k$ behaves regularly along the
grid: that it does not increase, and that it does not drop too sharply when $k$ doubles. For the scale \eqref{eq:tau-scale}, the next lemma proves both properties.

\begin{lemma}[Regularity of the stochastic scale]
\label{lem:tau-properties}
The sequence $(\tau_k)_{k\in\mathcal K_n}$ defined in \eqref{eq:tau-scale} satisfies,
for every $k$ with $k,2k\in\mathcal K_n$,
\[
\tau_{2k}\le\tau_k
\qquad\text{and}\qquad
\tau_{2k}\ge\tau_k/\sqrt2 .
\]
Moreover, with $C_1=36\cdot2048$,
\begin{equation}
\label{eq:tau-explicit}
\tau_k^2\le C_1\kappa_\sigma\left[\frac{s+1}{k}\log\Big(1+\frac{k(L+1)V^2}{s+1}\Big)
+\frac1k\log\Big(\frac{|\mathcal K_n|}{\eta}\Big)\right].
\end{equation}
\end{lemma}

No condition on $k_0$ is needed. Both properties follow from the fact that
$z\mapsto W(z)/z=e^{-W(z)}$ is decreasing, which is where the Lambert-$W$ form of
$\bar\delta_k$ pays off.

\subsection{Lepski's selection rule}
\label{sec:lepski}

Proposition~\ref{prop:simultaneous} bounds every candidate's error by its bias plus
its stochastic scale, and the bias is unknown. Lepski's idea is to detect the bias
indirectly. For small neighbourhoods the bias is negligible, so two such estimates can
differ only by stochastic fluctuation, of order $\tau$. Once $k$ is large enough for the
bias to matter, the estimate starts to drift away from the small-neighbourhood estimates by
more than the noise can explain. The rule stops just before this happens.

Define the adaptive neighbourhood size
\begin{equation}
\label{eq:lepski-rule}
\hat k(x)=\max\Big\{k\in\mathcal K_n:\big|\hat f_k(x)-\hat f_\ell(x)\big|\le A\tau_\ell
\ \text{for every}\ \ell\in\mathcal K_n,\ \ell\le k\Big\},
\end{equation}
where $A>0$ is a numerical constant, and set
$\hat f_{\mathrm{ad}}(x)=\hat f_{\hat k(x)}(x)$. The set in \eqref{eq:lepski-rule} always
contains $k_0$, so $\hat k(x)$ is well defined. Section~\ref{sec:oracle-adaptive} shows that
$A\ge4$ suffices.

The intuition behind the construction is straightforward. For a small
neighbourhood, the approximation error is small, but the stochastic
error is large. Increasing $k$ reduces the stochastic error because
more observations are used to train the network. However, the
neighbourhood radius increases with $k$, and therefore the local
approximation error increases. The selector chooses the largest
neighbourhood for which this increase in local approximation error
cannot yet be statistically distinguished from the stochastic
fluctuation.

\subsection{Oracle property of the adaptive estimator}
\label{sec:oracle-adaptive}

To analyse the rule we compare it with the choice that an oracle knowing $b_k(x)$
would make. For a fixed query point $x$, define the oracle neighbourhood
\begin{equation}
k_x^{*}=\max\{k\in\mathcal K_n:b_k(x)\le\tau_k\}.
\label{eq:oracle-index}
\end{equation}
Since $b_k(x)$ is increasing in $k$ while $\tau_k$ is  non-increasing
 (Lemma~\ref{lem:tau-properties}), the set in \eqref{eq:oracle-index} is an
initial segment of $\mathcal K_n$. It is the largest neighbourhood on which the bias is
still dominated by the stochastic error, and so the natural target. We need it to lie
strictly inside the grid.

\begin{assumption}[Adaptation window]
\label{ass:window}
$b_{k_0}(x)\le\tau_{k_0}$ and $b_{k_J}(x)>\tau_{k_J}$. Equivalently,
$k_0\le k_x^*\le k_J/2$.
\end{assumption}

The analysis rests on two facts. The first says that the rule never stops too early.

\begin{lemma}[The oracle neighbourhood is admissible]
\label{lem:admissible}
Let $A\ge4$ and suppose $b_{k_0}(x)\le\tau_{k_0}$. On the event $\mathcal A_x$ of
Proposition~\ref{prop:simultaneous},
\begin{equation}
\hat k(x)\ge k_x^{*}.
\label{eq:khat-lower}
\end{equation}
\end{lemma}

The reason is that on the good event every estimator up to $k_x^*$ is within
$b_\ell(x)+\tau_\ell\le2\tau_\ell$ of the truth, so any two of them are within $4\tau_\ell$ of
each other, and the test in \eqref{eq:lepski-rule} with $A\ge4$ cannot reject. For example, $A=4.5$ satisfies this constraint. The second fact says
that the oracle neighbourhood really balances the two error terms, rather than merely
satisfying $b\le\tau$ there.

\begin{lemma}[Balance at the oracle neighbourhood]
\label{lem:balance}
 Under Assumption~\ref{ass:window},
\begin{equation}
\label{eq:balance}
 \frac{1}{2\sqrt2}\,\tau_{k_x^{*}}<b_{k_x^{*}}(x)\le\tau_{k_x^{*}},
\qquad\text{so that}\qquad
(L_\F+3L_x)^2\Big(\frac{k_x^{*}}{n}\Big)^{2\alpha_x/d^*}\asymp\tau_{k_x^{*}}^2 .
\end{equation}
\end{lemma}

 This is where the upper end of the window is used: the lower bound comes from comparing
$k_x^*$ with the next grid point $2k_x^*$, at which the bias has overtaken the stochastic
term, and it needs $2k_x^*\in\mathcal K_n$. Together, the two lemmas give the oracle
inequality.

\begin{theorem}[Oracle inequality for the adaptive estimator]
\label{thm:adaptive-oracle}
 Let $A\ge4$ and let Assumption~\ref{ass:window} hold. On the event $\mathcal A_x$ of
Proposition~\ref{prop:simultaneous}, which has probability at least $1-2\eta$,
\begin{equation}
\big|\hat f_{\mathrm{ad}}(x)-f_0(x)\big|\le(A+2)\,\tau_{k_x^{*}},
\label{eq:adaptive-tau}
\end{equation}
 and consequently
\begin{equation}
\big|\hat f_{\mathrm{ad}}(x)-f_0(x)\big|^2\le 8(A+2)^2
\inf_{k\in\mathcal K_n}\left\{b_k(x)^2+\tau_k^2\right\}
\lesssim \inf_{k\in\mathcal K_n}\left\{(L_\F+3L_0)^2\Big(\frac kn\Big)^{2\alpha_x/d^*}+\tau_k^2\right\}.
\label{eq:oracle-adaptive}
\end{equation}
 In particular, by \eqref{eq:tau-explicit},
\begin{equation}
\begin{aligned}
\big|\hat f_{\mathrm{ad}}(x)-f_0(x)\big|^2\lesssim
\inf_{k\in\mathcal K_n}\Bigg\{&(L_\F+3L_0)^2\Big(\frac kn\Big)^{2\alpha_x/d^*}\\
&{}+\,\frac{s+1}{k}\log\Big(1+\frac{k(L+1)V^2}{s+1}\Big)
+\frac{1}{k}\log\Big(\frac{|\mathcal K_n|}{\eta}\Big)\Bigg\}.
\end{aligned}
\label{eq:oracle-explicit}
\end{equation}
\end{theorem}

 The right-hand side of \eqref{eq:oracle-adaptive} is, up to the constant $8(A+2)^2$,
the smallest error bound available to any single estimator in the collection, and in
particular to the one an oracle knowing $\alpha_x$, $L_x$ and $d^*$ would choose. The
selector achieves it without that knowledge. Two costs are hidden in the constant and the
logarithm. The constant reflects the fact that the rule controls differences between
estimators rather than errors directly. The term $\log(|\mathcal K_n|/\eta)$ is the price of
the union bound over the grid, and it is of order $\log\log n$ for fixed $\eta$.

\begin{remark}[Sequential form of the rule]
\label{rem:sequential-lepski}
Let $\tilde k(x)$ be the largest $k\in\mathcal K_n$ such that every
$k'\in\mathcal K_n$ with $k'\le k$ passes the test in
\eqref{eq:lepski-rule}. This is the rule obtained by scanning the grid
upwards and stopping at the first rejection, and $\tilde k(x)\le\hat k(x)$.
Lemma~\ref{lem:admissible} and Theorems~\ref{thm:adaptive-oracle}
and~\ref{thm:adaptive-rate} hold with $\tilde k(x)$ in place of $\hat k(x)$.
Indeed, on $\mathcal A_x$ the argument in the proof of
Lemma~\ref{lem:admissible} applies to every $k\le k_x^*$: for
$\ell\le k\le k_x^*$,
$|\hat f_k(x)-\hat f_\ell(x)|\le b_k(x)+\tau_k+b_\ell(x)+\tau_\ell\le4\tau_\ell$,
because $b_\ell(x)\le b_k(x)\le\tau_k\le\tau_\ell$. Hence every $k\le k_x^*$
passes, and $\tilde k(x)\ge k_x^*$. The proof of
Theorem~\ref{thm:adaptive-oracle} uses only this inequality and the test
between $\tilde k(x)$ and $\ell=k_x^*$, which $\tilde k(x)$ passes by
definition. The implementation in Section~\ref{sec:simulation} uses this
form.
\end{remark}

\subsection{Adaptation to local smoothness and intrinsic dimension}
\label{sec:adaptation}

 It remains to evaluate the infimum in \eqref{eq:oracle-explicit}. Ignoring logarithmic
factors and the network-complexity factors, the two dominant terms in
\eqref{eq:oracle-explicit} have the form
\begin{equation}
(L_\F+3L_0)^2\Big(\frac kn\Big)^{2\alpha_x/d^*}+\frac1k .
\label{eq:two-terms}
\end{equation}
The first term is increasing in $k$ and reflects the error incurred by using observations
farther from $x$; the second is decreasing in $k$ and reflects the stochastic error from
fitting a network on $k$ observations. Balancing only these two terms gives the heuristic scale
\begin{equation}
k_{\mathrm{heur}}(x)\asymp n^{\frac{2\alpha_x}{2\alpha_x+d^*}}
(L_\F+3L_0)^{-\frac{2d^*}{2\alpha_x+d^*}}.
\label{eq:kstar}
\end{equation}
 Keeping track of the complexity factors gives the following statement. Write
\[
 \mathcal D_n=\kappa_\sigma\Big[(s+1)\log\Big(1+\frac{n(L+1)V^2}{s+1}\Big)+\log\Big(\frac{|\mathcal K_n|}{\eta}\Big)\Big],
\]
 a deterministic upper bound for $k\tau_k^2/C_1$ over the grid, and let $k_x^\circ$ be the
value of $k$ at which $C_{\mathcal M}^2(L_\F+3L_0)^2(k/n)^{2\alpha_x/d^*}=C_1\mathcal D_n/k$.

\begin{theorem}[Adaptive pointwise rate]
\label{thm:adaptive-rate}
 Let $A\ge4$, let Assumption~\ref{ass:window} hold, and suppose that
$k_0\le k_x^\circ\le k_J$. Then, with probability at least $1-2\eta$,
\begin{equation}
\big|\hat f_{\mathrm{ad}}(x)-f_0(x)\big|^2\lesssim
 \Big(\frac{\mathcal D_n}{n}\Big)^{\frac{2\alpha_x}{2\alpha_x+d^*}}
(L_\F+3L_0)^{\frac{2d^*}{2\alpha_x+d^*}},
\label{eq:adaptive-rate}
\end{equation}
 where the implicit constant depends only on $A$ and $C_{\mathcal M}$.
\end{theorem}

Thus the adaptive estimator automatically attains the rate associated with the local smoothness
exponent $\alpha_x$, without requiring $\alpha_x$ to be known when selecting the neighbourhood
size.  Adaptation is to the exponent only: the H\"older constant enters
\eqref{eq:adaptive-rate} through the global bound $L_0$ rather than through the local
constant $L_x$. Moreover, the ambient feature dimension $d$ does not appear in the polynomial exponent of
the rate; the dimension entering the local approximation term is the intrinsic dimension
$d^*$ of the manifold support of $X$.  Nor does the procedure use $d^*$: it enters only
through the rate at which the neighbourhood radius shrinks, which the selector detects
indirectly through the bias. In particular, if $s$, $L$ and $\log V$ grow polylogarithmically, with the
same fixed $F,K,B,\sigma,L_0$,  we obtain the following.

\begin{corollary}[Adaptive rate up to logarithmic factors]
\label{cor:adaptive-OP}
Under the conditions of Theorem~\ref{thm:adaptive-rate}, suppose $\eta$ is
fixed, $F,K,B,\sigma,L_0$ are independent of $n$, and $s,L,\log V$ are
$O((\log n)^{C'})$. Then the single estimator defined with this fixed $\eta$
satisfies
\begin{equation}
\big|\hat f_{\mathrm{ad}}(x)-f_0(x)\big|^2
=O_{\P}\!\left(n^{-\frac{2\alpha_x}{2\alpha_x+d^*}}(\log n)^{C}\right)
\label{eq:adaptive-OP}
\end{equation}
for some constant $C>0$.
\end{corollary}

\begin{remark}[Confidence levels and expected risk]
\label{rem:adaptive-log}
For fixed $\eta$, the selection term $\log(|\mathcal K_n|/\eta)$ is
$O(\log\log n)$. The proof of Corollary~\ref{cor:adaptive-OP} establishes
a probability tending to one for that same fixed selector, rather than
changing its confidence parameter after it has been defined.
For expected risk, define the procedure with $\eta_n=n^{-1}$ and clip its
output to $[-B,B]$. If the grid and window conditions hold for this choice,
Theorem~\ref{thm:adaptive-rate} gives
\[
\E|\hat f_{\mathrm{ad}}^{\mathrm{clip}}(x)-f_0(x)|^2
\lesssim (\mathcal D_n/n)^{2\alpha_x/(2\alpha_x+d^*)}
(L_\F+3L_0)^{2d^*/(2\alpha_x+d^*)}+8B^2/n.
\]
Here $\mathcal D_n$ uses $\eta_n$, so its confidence term is $O(\log n)$.
Classical pointwise adaptation lower bounds \citep{lepski1997pointwise}
explain why logarithmic losses arise in adaptation, but we do not pursue that here.

\end{remark}

\subsection{The Doppler example revisited}
\label{sec:doppler-lepski}

\begin{figure}[htbp]
\centering
\includegraphics[width=3in, height = 2in]{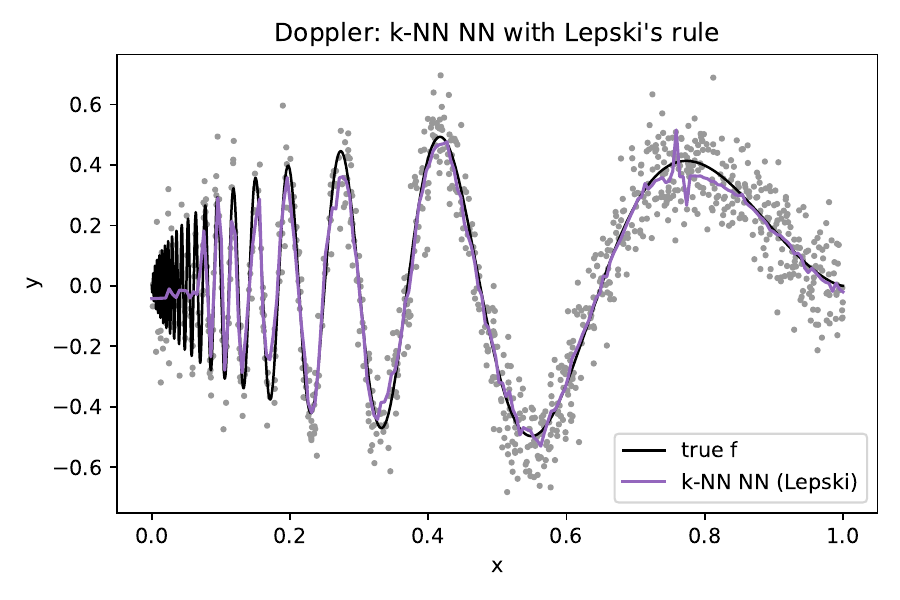}
\includegraphics[width=3in, height = 1.95in]{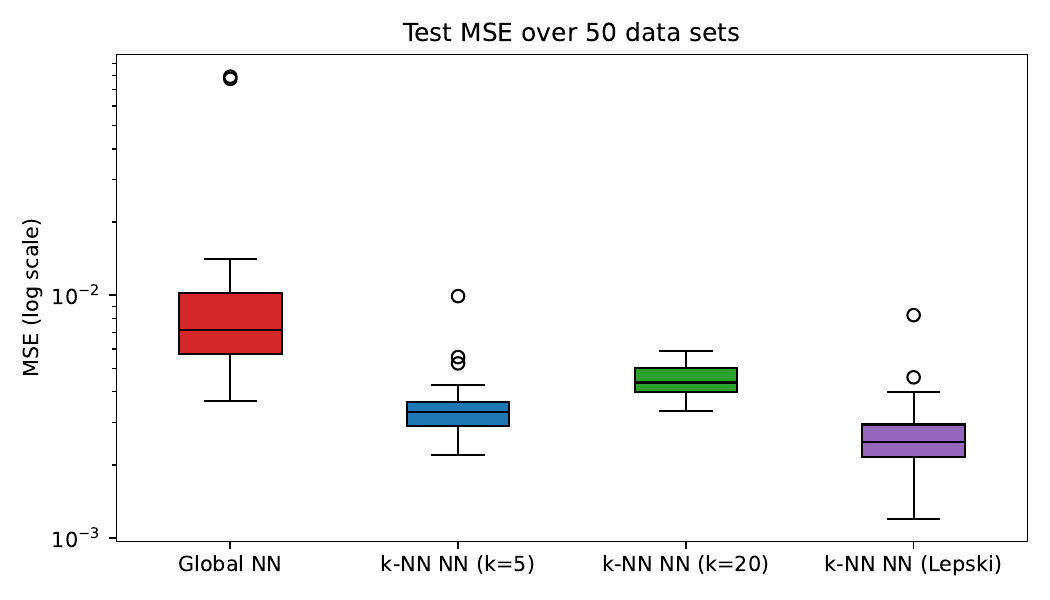}
\caption{The Doppler example of Figure~\ref{fig:doppler} with $k$ chosen by
Lepski's rule at each query point. (a) Data (grey), $f_0$ (black) and the
local fit with Lepski's choice of $k$ (purple). (b) Test MSE across the
$50$ data sets for the globally trained network, the local fits with $k=5$
and $k=20$, and the local fit with Lepski's choice of $k$.}
\label{fig:doppler-lepski}
\end{figure}

We return to the Doppler example of Section~\ref{sec:doppler}, which showed
that no single neighbourhood size suits the whole domain: $k=5$ follows the
oscillations near the origin but is noisy where $f_0$ varies slowly, and
$k=20$ does the reverse. Figure~\ref{fig:doppler-lepski} uses the same data
model, the same local architecture and the same $50$ data sets, and now
chooses $k$ at each query point by Lepski's rule. Panel (a) shows that the
selected fit follows the rapid oscillations near the origin and is smooth
where $f_0$ varies slowly, so that each part of the domain receives the
neighbourhood size it needs. Panel (b) shows that, in this example, the
Lepski fit has a smaller median test MSE than both fixed choices of $k$ and
the globally trained network. As in Section~\ref{sec:doppler}, these plots
concern this example and these architectures; a systematic comparison
follows in Section~\ref{sec:simulation}.

\section{Numerical experiments}
\label{sec:simulation}

In this section we study a practical version of the procedure on
simulated data. We address four questions.
\begin{enumerate}
\item \emph{Local adaptivity.} Does the localized network, with $k$ chosen
pointwise by Lepski's rule, compete with methods designed for spatially
inhomogeneous regression functions?
\item \emph{Activation.} Does the activation of the local network matter?
Our theory concerns ReLU networks; we also use the sigmoid, SiLU
and GELU $\big(\rho_2(x) = x\Phi(x)\big)$ activations.
\item \emph{Global versus local fitting.} Does localization improve on a
single deep network fitted to the whole sample?
\item \emph{Intrinsic dimension.} When the covariates lie on a
three-dimensional manifold in $\mathbb R^{30}$, does the local fit behave
as the intrinsic dimension, rather than the ambient dimension, suggests?
\end{enumerate}

In every design we observe $(X_i,Y_i)$, $i=1,\ldots,n$, with
\[
Y_i=f_0(X_i)+\sigma W_i,\qquad W_i\stackrel{\mathrm{i.i.d.}}{\sim}N(0,1),
\qquad\sigma=0.5,
\]
and the noise is independent of the design, as in Section~\ref{sec:rate}.
In the univariate designs the design is fixed and equispaced, $X_i=i/n$
with $n=10\,000$, and $f_0$ is one of the four test functions of
\citet{donoho1994ideal}:
\begin{align*}
\text{Blocks:}\quad & f_0(x)=\sum_{j=1}^{11}h_j\,\mathbf 1\{x\ge t_j\},\\
\text{Bumps:}\quad & f_0(x)=\sum_{j=1}^{11}a_j\,\bigl(1+|x-t_j|/w_j\bigr)^{-4},\\
\text{HeaviSine:}\quad & f_0(x)=4\sin(4\pi x)-\operatorname{sgn}(x-0.3)-\operatorname{sgn}(0.72-x),\\
\text{Doppler:}\quad & f_0(x)=\sqrt{x(1-x)}\,\sin\!\bigl(2.1\pi/(x+0.05)\bigr),
\end{align*}
where
\begin{align*}
t&=(0.10,0.13,0.15,0.23,0.25,0.40,0.44,0.65,0.76,0.78,0.81),\\
h&=(4,-5,3,-4,5,-4.2,2.1,4.3,-3.1,2.1,-4.2),\\
a&=(4,5,3,4,5,4.2,2.1,4.3,3.1,5.1,4.2),\\
w&=(0.005,0.005,0.006,0.01,0.01,0.03,0.01,0.01,0.005,0.008,0.005).
\end{align*}
Blocks is piecewise constant with eleven jumps, Bumps is smooth except at
eleven narrow peaks, HeaviSine is a smooth oscillation with two jumps, and
Doppler is the function \eqref{eq:doppler} of Section~\ref{sec:doppler}.

In the manifold designs a latent variable $z\sim\mathrm{Unif}[0,1]^3$ is
mapped to
\[
\psi(z)=\bigl(z_1,\,z_2,\,z_3,\,z_1^2,\,z_2z_3,\,\sin2\pi z_1,\,
\cos2\pi z_2,\,z_3^3\bigr)\in\mathbb R^8,
\]
which is padded with $22$ zeros and rotated into $\mathbb R^{30}$ by an
orthogonal matrix drawn once and used in every replication. The regression
function is
\[
f_0(X)=g(z_1)+0.8\sin(2\pi z_2z_3),\qquad
g\in\{\text{Blocks},\ \text{Bumps},\ \text{HeaviSine},\ \text{Doppler}\},
\]
which gives the designs CubicEmbed--Blocks, CubicEmbed--Bumps,
CubicEmbed--HeaviSine and CubicEmbed--Doppler, with $d^*=3$, $d=30$ and
$n=2000$. The support is the image of the cube $[0,1]^3$ under a bi-Lipschitz map,
a three-dimensional manifold with corners. Indeed, the first three
coordinates of $\psi(z)$ equal $z$ and the derivatives of $\psi$ are
bounded, so $\|z-z'\|\le\|\psi(z)-\psi(z')\|\le C\|z-z'\|$ for a finite
$C$, and padding and rotation preserve distances. Since the intersection of
$[0,1]^3$ with a ball of radius $r\le1$ centred in the cube has volume at
least a constant times $r^3$, the lower local mass condition used for the
upper bounds holds with $d^*=3$ (see the remark after
Corollary~\ref{cor:approx-rate}). The support is not contained in
$[0,1]^{30}$, and the inputs are not rescaled. All distance-based steps,
including every nearest-neighbour search, are the same as for the
eight-dimensional data $\psi(z)$. No method is
given $d^*$, the latent coordinates or the embedding.

For a query point $x$ and a neighbourhood size $k$, LDNN--Lepski (localized
deep neural network with Lepski's rule) fits a network with one hidden
layer of five units,
\[
x'\mapsto b_0+\sum_{j=1}^5b_j\,\phi\bigl(a_{j0}+a_j^\top x'\bigr),
\]
to the $k$ nearest neighbours of $x$, and evaluates the fit at
$x$. As in the theory, the fit is by least squares on the $k$
neighbours, with no penalty on the network parameters (no weight decay).
The criterion is minimised by the Broyden–Fletcher–Goldfarb–Shanno (BFGS)
algorithm \citep{dennis1977quasi} from a single start, with weights drawn from
$\mathrm{Unif}[-0.7,0.7]$, at most $200$ iterations and a gradient-norm
tolerance of $10^{-5}$; the inputs are used on their original scale. The activation $\phi$ is ReLU, the logistic sigmoid
$s(z)=1/(1+e^{-z})$, SiLU $z\,s(z)$, or GELU $z\,\Phi(z)$, where $\Phi$ is
the standard normal distribution function, and the four versions start
from the same weights. The network has $5d+11$ parameters: $16$ in the
univariate designs and $161$ in the manifold designs. The candidate
neighbourhood sizes form the geometric grids
$\mathcal K=\{10,20,\ldots,640\}$ (univariate) and $\{10,20,\ldots,320\}$
(manifold), and at each query point we use the following form of
\eqref{eq:lepski-rule}:
\begin{equation}
\label{eq:lepski-practical}
\hat k(x)=\max\Big\{k\in\mathcal K:\ \big|\hat f_{k'}(x)-\hat f_\ell(x)\big|
\le A\hat\tau_\ell\ \text{ for all }\ell,k'\in\mathcal K
\text{ with }\ell<k'\le k\Big\},
\qquad
\hat\tau_\ell=\hat\sigma\sqrt{\frac{\log n}{\ell}},
\end{equation}
with $A=1$ and
$\hat\sigma=\operatorname{median}_i|Y_i-Y_{\mathrm{nn}(i)}|/(0.6745\sqrt2)$,
where $\mathrm{nn}(i)$ indexes the nearest other design point to $X_i$; in
one dimension these are consecutive differences. The grid is thus scanned
upwards and the scan stops at the first rejection.

To separate the effect of the local network from that of the selection
rule, two linear local smoothers use \eqref{eq:lepski-practical} with the
same $A$ and $\hat\sigma$. Local linear--Lepski is weighted least squares
on $(1,X_i-x)$ with Epanechnikov weights and bandwidth equal to the
distance from $x$ to its $k$th nearest neighbour. In the manifold designs
the covariates lie in an eight-dimensional linear subspace of
$\mathbb R^{30}$, so the $31\times31$ weighted Gram matrix has rank at most
nine, and we use its Moore--Penrose pseudoinverse (relative tolerance
$10^{-10}$). Every weighted least-squares solution has the same intercept:
the query point is one of its own neighbours, so the null space of the
weighted design matrix consists of slope vectors orthogonal to that
subspace. Hence the fitted value $\hat f_k(x)=\sum_il_i(x)Y_i$ is well
defined. The plug-in scale $\hat\sigma\|l(x)\|_2$, the conditional standard
deviation of $\hat f_k(x)$ for fixed $k$ with $\sigma$ replaced by
$\hat\sigma$, replaces $\hat\sigma/\sqrt k$ in $\hat\tau_k$. $k$-NN--Lepski is the average of the
$k$ nearest responses. Their grids are $\{8,16,\ldots,2048\}$ and
$\{16,32,\ldots,512\}\cap[2(d+1),n)$ for local linear--Lepski, and
$\{4,8,\ldots,2048\}$ and $\{4,8,\ldots,512\}$ for $k$-NN--Lepski. Since
these grids differ from those of LDNN--Lepski, the comparison does not
control for the candidate neighbourhood sizes.
Trend filtering \citep{kim2009ell1,tibshirani2014adaptive} has order $r=0$
for Blocks and $r=2$ for the other functions; it is computed by the
algorithm of \citet{condat2013direct} for $r=0$ and the ADMM algorithm of
\citet{ramdas2016fast} for $r=2$ along a path of $30$ values of $\lambda$.
The tuning parameter minimises Stein's unbiased risk estimate
\citep{stein1981estimation} with $\hat\sigma$ in place of $\sigma$ and the
number of knots plus $r+1$ as degrees of freedom
\citep{tibshirani2012degrees}; with a plug-in noise level the criterion is
not exactly unbiased. In the manifold designs
we use order-$0$ trend filtering on the symmetrised $5$-nearest-neighbour
graph, the $k$-NN fused lasso of \citet{padilla2020adaptive}, computed by
ADMM over $25$ geometrically spaced values $\lambda\in\hat\sigma[10^{-2},10]$.
Its tuning parameter minimises the same plug-in criterion, with the number
of fused connected components of the graph, the groups of nodes joined by
edges whose estimated differences are exactly zero, as degrees of freedom
\citep{tibshirani2012degrees}. The cubic smoothing spline is tuned by generalised
cross-validation \citep{craven1979smoothing}. The global networks are
fully connected, with (depth, width) in $\{(4,64),(8,64)\}$, one for each of
the four activations; depths of this order, $\log n$, are used in the
theory of deep ReLU estimators \citep{dnnregression}. They are trained by
Adam \citep{kingma2015adam}, with mini-batches of $256$, learning rate
$10^{-3}$ and at most $500$ epochs, on $90\%$ of the sample; the remaining
$10\%$ is used for early stopping, for two learning-rate reductions and for
the choice of architecture, so the global networks are the only methods
tuned on independent data. The smoothing spline is used in the univariate
designs only: a thin-plate spline in $\mathbb R^{30}$ needs a penalty of
order $m\ge16$, whose null space has dimension
$\binom{m-1+d}{d}\approx3.4\times10^{11}$, far more than $n$.

In each replication all methods are applied to the same data set and
evaluated at a common set of design points $x_1,\ldots,x_m$: $m=300$
equally spaced points in the univariate designs and $m=100$ design points
in the manifold designs. The evaluation points belong to the training
sample, so the methods are compared as denoisers; trend filtering and the
smoothing spline are evaluated at their fitted values. The performance measure is
\[
\mathrm{MSE}=\frac1m\sum_{j=1}^m\{\hat f(x_j)-f_0(x_j)\}^2,
\]
the average of the squared pointwise errors studied in
Sections~\ref{sec:rate} and~\ref{sec:adaptive}. Each design has $200$
replications, except CubicEmbed--Doppler, for which $142$ replications
were completed. The random number generator is seeded by the pair (design,
replication), and the rotation matrix by a fixed seed.

\subsection{Results for the univariate designs}
\label{sec:sim-uni}

\begin{figure}[htbp]
\centering
\includegraphics[width=\textwidth]{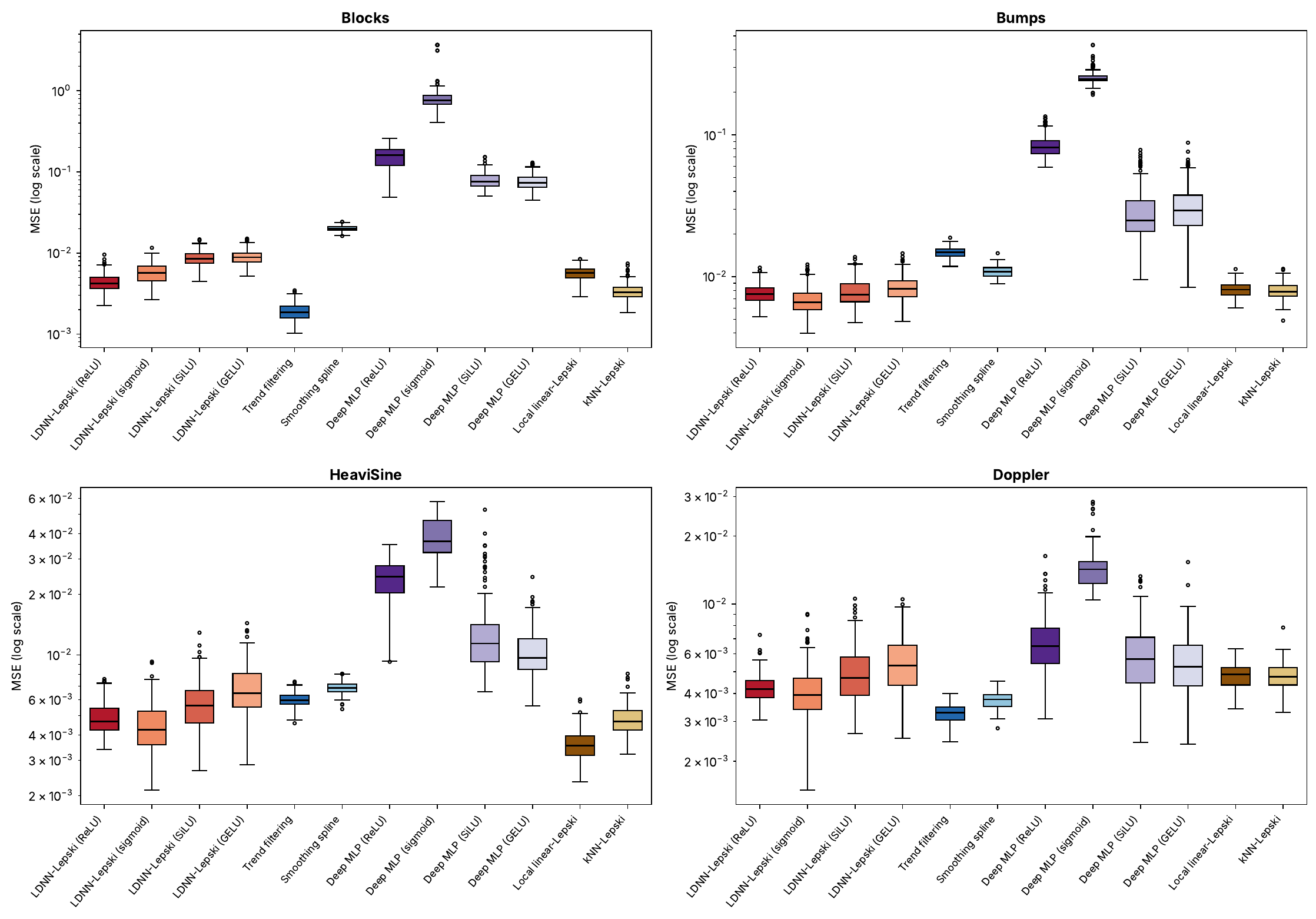}
\caption{MSE (log scale) over $200$ replications for the univariate designs.}
\label{fig:uni}
\end{figure}

Figure~\ref{fig:uni} shows the distribution of the MSE over the
replications. LDNN--Lepski with the sigmoid activation has the smallest median MSE for
Bumps. For the other three functions a competitor is best: order-$0$ trend
filtering for Blocks and order-$2$ trend filtering for Doppler, and local
linear--Lepski for HeaviSine. The sigmoid version of
LDNN--Lepski is second for HeaviSine, and for Doppler the sigmoid and ReLU
versions follow trend filtering and the smoothing spline; for Blocks the
ReLU version is third, after trend filtering and $k$-NN--Lepski. We report
the four versions separately rather than selecting one after seeing the
results. The sigmoid has smaller MSE than SiLU and GELU for all four
functions, and ReLU for Blocks, HeaviSine and Doppler; for Bumps ReLU and
SiLU are level. GELU has the largest median MSE of the four activations in
every case. Between ReLU and the sigmoid, ReLU is better for Blocks and the
sigmoid for Bumps and HeaviSine, while for Doppler they are close. These
are comparisons of the error averaged over the evaluation grid; they do not
show where in the domain each activation is more accurate. With the same
form of selection rule, the ReLU and sigmoid networks improve on local
linear--Lepski and $k$-NN--Lepski for Bumps and Doppler, in $61$--$84\%$ of
the replications. For HeaviSine local linear--Lepski is better than every
version of the network, and for Blocks $k$-NN--Lepski is. The ReLU and sigmoid versions have a smaller
median MSE than every fully connected global network for every function;
the margin is large for Blocks, Bumps and HeaviSine and smaller for
Doppler, where the SiLU and GELU global networks are within a factor of
about $0.7$--$0.8$ of them. Many global networks reached the epoch limit
for Blocks and Bumps. LDNN--Lepski selects the
largest candidate, $k=640$, at $52$--$78\%$ of the evaluation points for
Blocks and Bumps and at $77$--$88\%$ for HeaviSine and Doppler. The upper end of its grid is therefore often
binding, and a larger grid may change these comparisons.

\subsection{Results for the manifold designs}
\label{sec:sim-mv}

\begin{figure}[htbp]
\centering
\includegraphics[width=\textwidth]{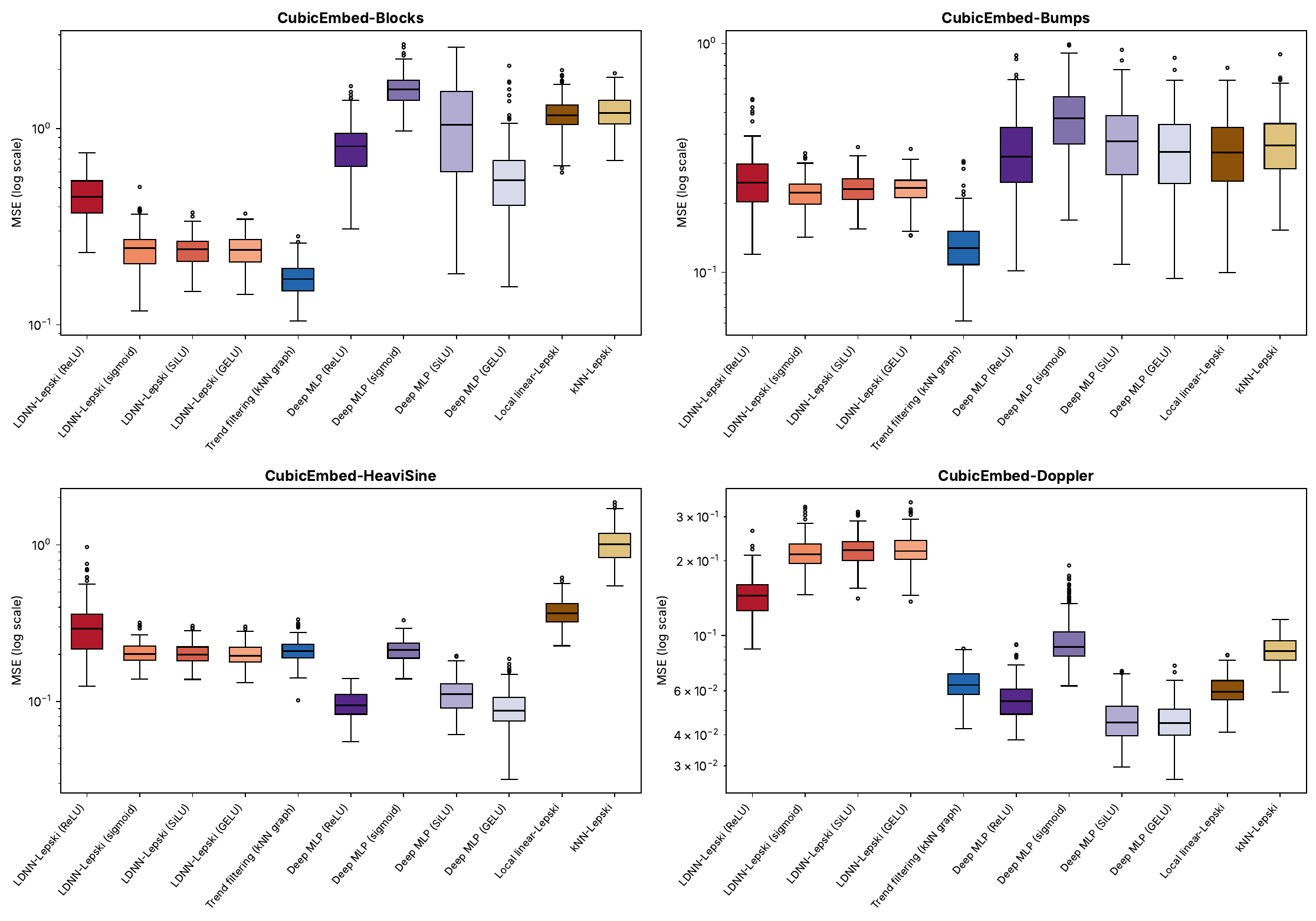}
\caption{MSE (log scale) over $200$ replications ($142$ for
CubicEmbed--Doppler) for the manifold designs.}
\label{fig:mv}
\end{figure}

Figure~\ref{fig:mv} shows the distribution of the MSE for the manifold
designs. Trend filtering on the $k$-NN graph has the smallest median MSE for
CubicEmbed--Blocks and CubicEmbed--Bumps, followed by the smooth
LDNN--Lepski versions. The GELU global network has the smallest median for
CubicEmbed--HeaviSine and, level with the SiLU global network, for
CubicEmbed--Doppler. For CubicEmbed--HeaviSine the
smooth LDNN--Lepski versions follow the ReLU, SiLU and GELU global networks
and are level with graph trend filtering, while for CubicEmbed--Doppler
every version of LDNN--Lepski has a larger median MSE than every
competitor. The smooth activations are better than
ReLU for the Blocks and HeaviSine variants, slightly better for
CubicEmbed--Bumps, and ReLU is better for CubicEmbed--Doppler, so the
ordering of the activations differs from that in one dimension; the three
smooth activations behave alike, with median MSE ratios within $5\%$ of one. With the same form of selection rule, the smooth
versions improve on local linear--Lepski and $k$-NN--Lepski for the Blocks,
Bumps and HeaviSine variants, but not for CubicEmbed--Doppler. The smooth local networks are better than every
global network for the Blocks and Bumps variants, and worse than the ReLU,
SiLU and GELU global networks for the HeaviSine and Doppler variants.
LDNN--Lepski selects the largest candidate, $k=320$, at no more than $14\%$
of the evaluation points, whereas local linear--Lepski selects its largest
candidate, $k=512$, at $90$--$93\%$ of the points for the Bumps and Doppler
variants. With ReLU, LDNN--Lepski selects the smallest candidate, $k=10$,
at $16$--$31\%$ of the points.

\subsection{Summary of the experiments}
\label{sec:sim-summary}

In one dimension, LDNN--Lepski with the sigmoid activation has the
smallest median MSE for Bumps, and its ReLU and sigmoid versions are close
behind the best method for HeaviSine and Doppler; trend filtering is best
for Blocks and Doppler and local linear--Lepski for HeaviSine. For every
function the ReLU and sigmoid versions have a smaller median MSE than each
of the global networks we trained. The best activation depends on the
function and the design: in one dimension the sigmoid, and for three of
the four functions also ReLU, do better than SiLU and GELU, whereas on the
manifold the
smooth activations do better for three of the four variants.
On the manifold, the smooth local networks improve on the local linear and
local constant fits for three of the four variants, and on every global
network for the Blocks and Bumps variants, for which graph trend filtering
is best; the global networks are best for the HeaviSine and Doppler
variants.

The experiments have limitations. They use a single manifold design, with
$n=2000$, $d^*=3$ and $d=30$, so they do not show how the error scales with
$n$, $d^*$ or $d$. LDNN--Lepski often selects the largest candidate in one
dimension, and the global networks often reached the epoch limit for
Blocks and Bumps; a larger grid, a longer training budget or restarts could
change the comparisons. The local networks are fitted without any
penalty and from a single start, so their accuracy also depends on the
optimiser. With $200$ replications a win percentage has a standard error of
at most $3.5$ points ($4.2$ points with the $142$ replications of
CubicEmbed--Doppler), so small differences in the paired comparisons should
not be read as established differences.

\bibliographystyle{apalike}
\bibliography{references}

\appendix

\appendix
\section{Proofs for Section~\ref{sec:rate}}
\label{app:proofs-rate}

This appendix collects the derivations behind the statements of
Section~\ref{sec:rate}. Throughout, $x$ is the fixed query point, $N_k(x)$ its
$k$-nearest-neighbour index set, and every bound stated conditionally on
$X_1,\ldots,X_n$ is flagged where the distinction matters.

\begin{proof}[Proof of Lemma~\ref{lem:basic-inequality}]
By the definition of $\hat f_n(x)$, we have
$$\frac{1}{k}\sum_{i\in N_k(x)}\big(Y_i - \hat f_n(X_i)\big)^2 \leq \frac{1}{k}\sum_{i\in N_k(x)}\big(Y_i - f_0(X_i)\big)^2.$$
This is where we use the assumption $f_0\in \mathcal{F}(L, \bm p, s)$; later we will relax this assumption. Then, we get
\begin{align*}
    & \frac{1}{k}\sum_{i\in N_k(x)}\big(Y_i - \hat f_n(X_i)\big)^2\\
    & = \frac{1}{k}\sum_{i\in N_k(x)}\big(Y_i - f_0(X_i) + f_0(X_i) - \hat f_n(X_i)\big)^2\\
    & = \frac{1}{k}\sum_{i\in N_k(x)}\big(Y_i - f_0(X_i)\big)^2 + \frac{1}{k}\sum_{i\in N_k(x)}\big(f_0(X_i) - \hat f_n(X_i)\big)^2 + \frac{2}{k}\sum_{i\in N_k(x)} \big(Y_i-f_0(X_i)\big)\big(f_0(X_i) - \hat f_n(X_i)\big)\\
    & = \frac{1}{k}\sum_{i\in N_k(x)}\big(Y_i - f_0(X_i)\big)^2 + \|\hat f_n - f_0\|_{2,k}^2 + \frac{2}{k}\sum_{i\in N_k(x)} \sigma W_i\big(f_0(X_i) - \hat f_n(X_i)\big) \leq \frac{1}{k}\sum_{i\in N_k(x)}\big(Y_i - f_0(X_i)\big)^2.
\end{align*}
Hence, we must have
\begin{align*}
    \|\hat f_n - f_0\|_{2,k}^2 &+ \frac{2}{k}\sum_{i\in N_k(x)} \sigma W_i\big(f_0(X_i) - \hat f_n(X_i)\big) \leq 0\\
    \|\hat f_n - f_0\|_{2,k}^2 &\leq \frac{2}{k}\sum_{i\in N_k(x)} \sigma W_i\big(\hat f_n(X_i)-f_0(X_i)\big).
\end{align*}
\end{proof}

\begin{proof}[Proof of Lemma~\ref{lem:critical-exists}]
Now for any $0<\delta\leq t$ and $\|h\|_{2,k} \leq t$ we have
\begin{align*}
    \frac{1}{k}\Big|\frac{\delta}{t}\sum_{i\in N_k(x)} W_i h(X_i)\Big| & = \frac{1}{k}\Big|\sum_{i\in N_k(x)} W_i \big(\frac{\delta}{t}h(X_i)\big)\Big|\leq \sup_{g\in\partial\F:\|g\|_{2,k}\leq \delta} \frac{1}{k}\Big|\sum_{i\in N_k(x)} W_i g(X_i)\Big|.
\end{align*}
This holds since $\frac{\delta}{t} h\in \partial\mathcal{F}$ by star-shapedness of $\partial\mathcal F$. Hence,
\begin{align}
    \frac{\delta}{t} G_n(t;\partial\F) \leq G_n(\delta;\partial\F) \Rightarrow \frac{G_n(t;\partial\F)}{t} \leq \frac{G_n(\delta;\partial\F)}{\delta}.
\end{align}
\end{proof}

\begin{proof}[Proof of Proposition~\ref{prop:training-error}]
For that, first let us consider $g\in\partial\F$ with $\|g\|_{2,k}\geq u$ such that
$$\Big|\frac{1}{k}\sum_{i\in N_k(x)} W_i g(X_i)\Big|\geq 2 \|g\|_{2,k} u.$$
Define $\Tilde{g} = \frac{u}{\|g\|_{2,k}} g,$ then $\|\Tilde{g}\|_{2,k} = u$. Since $g\in\partial\F$ which is star-shaped and $\frac{u}{\|g\|_{2,k}}\in [0,1]$, $\Tilde{g}\in\partial\F$. Also,
\begin{align*}
    \left|\frac{1}{k}\sum_{i\in N_k(x)} W_i \Tilde{g}(X_i)\right| & = \frac{u}{\|g\|_{2,k}} \left|\frac{1}{k}\sum_{i\in N_k(x)} W_i g(X_i)\right| \geq \frac{ 2 \|g\|_{2,k} u^2}{\|g\|_{2,k}} = 2u^2.
\end{align*}
Define,
$$\mathcal{A}(u) = \left\{\exists g\in \{h\in\partial\F: \|h\|_{2,k}\geq u\}: \Big|\frac{1}{k}\sum_{i\in N_k(x)} W_i g(X_i)\Big|\geq 2\|g\|_{2,k} u\right\}.$$
Therefore, whenever $\mathcal{A}(u)$ is true, there is $g\in\partial\F$ with $\|g\|_{2,k} = u$ such that
$$\left|\frac{1}{k}\sum_{i\in N_k(x)} W_i g(X_i)\right| \geq 2u^2 \Rightarrow \sup_{g\in\partial\F, \|g\|_{2,k}\leq u}\left|\frac{1}{k}\sum_{i\in N_k(x)} W_i g(X_i)\right| \geq 2u^2.$$
Therefore, 
$$\P[\mathcal{A}(u)] \leq \P\left[\sup_{g\in\partial\F, \|g\|_{2,k}\leq u}\left|\frac{1}{k}\sum_{i\in N_k(x)} W_i g(X_i)\right| \geq 2u^2\right].$$
Define
$$Z_{n,u}(W_1,\ldots, W_n) = \sup_{g\in\partial\F, \|g\|_{2,k}\leq u}\left|\frac{1}{k}\sum_{i\in N_k(x)} W_i g(X_i)\right|.$$
Conditioned on the feature vectors $X_1,\ldots, X_n$, $Z_{n,u}(W_1,\ldots, W_n)$ is a supremum of a Gaussian process. Now note that conditioned on $X_1,\ldots, X_n$,  
\begin{align*}
    & \Big|\sup_{g\in\partial\F, \|g\|_{2,k}\leq u}\big|\frac{1}{k}\sum_{i\in N_k(x)} W_i g(X_i)\big|-\sup_{g\in\partial\F, \|g\|_{2,k}\leq u}\big|\frac{1}{k}\sum_{i\in N_k(x)} W_i^\prime g(X_i)\big|\Big|\\
    & \leq \sup_{g\in\partial\F, \|g\|_{2,k}\leq u}\Big|\frac{1}{k}\sum_{i\in N_k(x)} (W_i -W_i^\prime) g(X_i)\Big|\\
    & \leq \sup_{g\in\partial\F, \|g\|_{2,k}\leq u} \|g\|_{2,k} \sqrt{\frac{1}{k} \sum_{i\in N_k(x)} (W_i -W_i^\prime)^2} \\
    & \leq\frac{u}{\sqrt{k}} \sqrt{\sum_{i=1}^n (W_i-W_i^\prime)^2}.
\end{align*}
Here the first inequality follows from the inequality,
$|\sup_g|a_g|-\sup_g|b_g||\le\sup_g|a_g-b_g|$, and the second follows from the Cauchy--Schwarz inequality
$$|k^{-1}\sum_{i\in N_k(x)}(W_i-W_i')g(X_i)|\le\{k^{-1}\sum_{i\in N_k(x)}(W_i-W_i')^2\}^{1/2}\|g\|_{2,k}.$$ Hence, $Z_{n,u}(W_1,\ldots, W_n)$ is a Lipschitz function of Gaussian random variables with Lipschitz constant $\frac{u}{\sqrt{k}}$. Hence, by the concentration of Lipschitz functions of Gaussian variables, we obtain
$$\P\Big[Z_{n,u}(W_1,\ldots, W_n) \geq \E[Z_{n,u}(W_1,\ldots, W_n)] + s\Big]\leq \exp\Big\{-\frac{ks^2}{2u^2}\Big\},$$
for all $s>0$. Choosing $s = u^2$ gives
$$\P\Big[Z_{n,u}(W_1,\ldots, W_n) \geq \E[Z_{n,u}(W_1,\ldots, W_n)] + u^2\Big]\leq \exp\Big\{-\frac{ku^2}{2}\Big\}.$$
Also, $\E[Z_{n,u}(W_1,\ldots, W_n)] = G_n(u;\partial\F)$. If $\delta_n$ is the smallest solution of $G_n(\delta; \partial\F)\leq c\delta^2$, then for any $u\geq \delta_n$ we get 
$$\frac{G_n(u;\partial\F)}{u} \leq \frac{G_n(\delta_n;\partial\F)}{\delta_n} \leq c \delta_n,$$
i.e., $\E[Z_{n,u}(W_1,\ldots, W_n)] \leq c \delta_n u$. Taking $c = 1$ and combining with the tail bound gives us
\begin{align*}
    & \P[\mathcal{A}(u)]\leq \P[Z_{n,u}(W_1,\ldots, W_n)\geq 2u^2]\\
    & \leq \P[Z_{n,u}(W_1,\ldots, W_n)\geq \delta_n u + u^2]\\
    & \leq \P[Z_{n,u}(W_1,\ldots, W_n)\geq \E[Z_{n,u}(W_1,\ldots, W_n)] + u^2] \leq \exp\Big\{-\frac{ku^2}{2}\Big\}.
\end{align*}
Here the second inequality uses $\delta_n u\le u^2$, which holds because $u\ge\delta_n$, and
the third uses $\E[Z_{n,u}]\le\delta_n u$. Taking $u = \sqrt{t\delta_n}$ for some $t\geq\delta_n$, so that $u\ge\delta_n$, gives us 
$$\P[\mathcal{A}(\sqrt{t\delta_n})^c] \geq 1-\exp\Big\{-\frac{tk\delta_n}{2}\Big\}.$$

Now suppose that $\mathcal{A}(u)^c$ holds with $u=\sqrt{t\delta_n}$. Since $\hat f_n$ and $f_0$
both lie in $\mathcal F_n$, the difference $\hat f_n-f_0$ belongs to $\partial\F$, and we distinguish
two cases according to its size.

\noindent\textbf{Case 1:} $\|\hat f_n-f_0\|_{2,k}<\sqrt{t\delta_n}$. The event $\mathcal A(u)^c$ says
nothing about functions of norm smaller than $u$, but none is needed, since in this case
directly
\[
\|\hat f_n-f_0\|_{2,k}^2< t\delta_n .
\]

\noindent\textbf{Case 2:} $\|\hat f_n-f_0\|_{2,k}\ge\sqrt{t\delta_n}$. Then $\hat f_n-f_0$ is one of the
functions over which $\mathcal A(u)$ quantifies, so on $\mathcal A(u)^c$
\begin{align*}
    \Big|\frac{1}{k}\sum_{i\in N_k(x)} W_i \big(\hat f_n(X_i)-f_0(X_i)\big)\Big| \leq 2 \|\hat f_n- f_0\|_{2,k} \sqrt{t\delta_n}.
\end{align*}
Consequently, the basic inequality \eqref{eq:basic-inequality} implies
\begin{align*}
    \|\hat f_n - f_0\|_{2,k}^2 \leq 4\sigma\|\hat f_n - f_0\|_{2,k} \sqrt{t\delta_n} \Rightarrow \|\hat f_n - f_0\|_{2,k}^2 \leq 16\sigma^2 t\delta_n.
\end{align*}
Combining the two cases, on $\mathcal A(\sqrt{t\delta_n})^c$ we have
$\|\hat f_n-f_0\|_{2,k}^2\le\max\{1,16\sigma^2\}\,t\delta_n=\kappa_\sigma t\delta_n$. Hence,
$$\P[\|\hat f_n - f_0\|_{2,k}^2 \leq \kappa_\sigma t\delta_n\mid X] \geq \P[\mathcal{A}(\sqrt{t\delta_n})^c\mid X] \geq 1 - \exp\Big\{-\frac{tk\delta_n}{2}\Big\}.$$
\end{proof}

\begin{proof}[Proof of Corollary~\ref{cor:training-expectation}]
All probabilities and expectations below are conditional on $X_1,\ldots,X_n$. Fix $u\ge\kappa_\sigma t\delta_n$
and put $t'=u/(\kappa_\sigma\delta_n)$. Then $t'\ge t\ge\delta_n$, so
Proposition~\ref{prop:training-error} applies with $t'$ in place of $t$ and gives
\[
\P\Big[\|\hat f_n-f_0\|_{2,k}^2>u\Big]=\P\Big[\|\hat f_n-f_0\|_{2,k}^2>\kappa_\sigma t'\delta_n\Big]
\le\exp\Big\{-\frac{kt'\delta_n}{2}\Big\}=\exp\Big\{-\frac{ku}{2\kappa_\sigma}\Big\}.
\]
This further gives us for $t\geq \delta_n$
\begin{align*}
    & \E\left[\frac{1}{k} \sum_{i\in N_k(x)} \big(\hat f_n(X_i) - f_0(X_i)\big)^2\right]\\
    & = \int_0^\infty \P\left[\frac{1}{k} \sum_{i\in N_k(x)} \big(\hat f_n(X_i) - f_0(X_i)\big)^2 \geq u\right] \mathrm d u\\
    & = \int_0^{\kappa_\sigma t\delta_n} \P\left[\frac{1}{k} \sum_{i\in N_k(x)} \big(\hat f_n(X_i) - f_0(X_i)\big)^2 \geq u\right] \mathrm d u + 
    \int_{\kappa_\sigma t\delta_n}^\infty \P\left[\frac{1}{k} \sum_{i\in N_k(x)} \big(\hat f_n(X_i) - f_0(X_i)\big)^2 \geq u\right] \mathrm d u\\
    & \leq \kappa_\sigma t\delta_n + \int_{\kappa_\sigma t\delta_n}^\infty \exp\Big\{-\frac{k\delta_n \frac{u}{\kappa_\sigma\delta_n}}{2}\Big\} \mathrm du\\
    & = \kappa_\sigma t\delta_n + \int_{\kappa_\sigma t\delta_n}^\infty \exp\Big\{-\frac{ku}{2\kappa_\sigma}\Big\} \mathrm du\\
    & = \kappa_\sigma t\delta_n + \frac{\exp\Big\{-\frac{ku}{2\kappa_\sigma}\}}{-\frac{k}{2\kappa_\sigma}}\Big|_{\kappa_\sigma t \delta_n}^\infty\\
    & = \kappa_\sigma t\delta_n + \frac{2\kappa_\sigma}{k} \exp\Big\{-\frac{kt \delta_n}{2}\Big\}
\end{align*}
This holds for all $t\geq\delta_n$. So, putting $t = \delta_n$ we get a crude upper bound
$$\E\left[\frac{1}{k} \sum_{i\in N_k(x)} \big(\hat f_n(X_i) - f_0(X_i)\big)^2\right] \leq \kappa_\sigma\delta_n^2 + \frac{2\kappa_\sigma}{k}.$$
However, one can get an even tighter upper bound 
$$\E\left[\frac{1}{k} \sum_{i\in N_k(x)} \big(\hat f_n(X_i) - f_0(X_i)\big)^2\right] \leq \inf_{t\geq \delta_n} \Big(\kappa_\sigma t\delta_n + \frac{2\kappa_\sigma}{k} \exp\Big\{-\frac{kt \delta_n}{2}\Big\}\Big).$$
\end{proof}

\begin{proof}[Proof of Lemma~\ref{lem:dudley}]
Write $\varepsilon=\delta^2/(2C_1^2)$ and $\mathcal N(t)=\mathcal N(t;\mathbb B(\delta;\partial\F),\|\cdot\|_{2,k})$
for the internal covering number. If $\varepsilon\ge\delta$ there is nothing to prove, since the
Cauchy--Schwarz bound gives $G_n(\delta;\partial\F)\le\delta\le\varepsilon$. Suppose therefore that
$\varepsilon<\delta$, which is the case for every $\delta\in(0,C_1]$ when $C_1\ge1$. So, we can construct a minimal $\delta^2/2C_1^2$-covering of $\mathbb{B}(\delta; \partial\F) = \{h\in \partial\F: \|h\|_{2,k}\leq \delta\}$, i.e., there exists $\{g_1,\ldots, g_M\}\subset \mathbb{B}(\delta;\partial\F)$ such that for all $g\in \mathbb{B}(\delta; \partial\F)$ there exists $j\in\{1,2,\ldots, M\}$ such that $\|g-g_j\|_{2,k}\leq \delta^2/2C_1^2$. Consequently,
\begin{align*}
    & \Big|\frac{1}{k}\sum_{i\in N_k(x)} W_i g(X_i)\Big| \\
    & \leq \Big|\frac{1}{k}\sum_{i\in N_k(x)} W_i g_j(X_i)\Big| + \Big|\frac{1}{k}\sum_{i\in N_k(x)} W_i \big(g(X_i)-g_j(X_i)\big)\Big|\\\
    & \leq \max_{j=1,2,\ldots, M} \Big|\frac{1}{k}\sum_{i\in N_k(x)} W_i g_j(X_i)\Big| + \sqrt{\frac{1}{k}\sum_{i\in N_k(x)} W_i^2} \sqrt{\frac{1}{k}\sum_{i\in N_k(x)} \big(g(X_i)-g_j(X_i)\big)^2}\\
    & \leq \max_{j=1,2,\ldots, M} \Big|\frac{1}{k}\sum_{i\in N_k(x)} W_i g_j(X_i)\Big| + \sqrt{\frac{1}{k}\sum_{i\in N_k(x)} W_i^2} ~~~\frac{\delta^2}{2C_1^2}.
\end{align*}
Hence, we obtain
$$G_n(\delta;\partial\F) \leq \E\left[\max_{j=1,2,\ldots, M} \Big|\frac{1}{k}\sum_{i\in N_k(x)} W_i g_j(X_i)\Big|\right] + \frac{\delta^2}{2C_1^2},$$
where we have used $\E[\sqrt{\frac{1}{k}\sum_{i\in N_k(x)} W_i^2}] \leq 1$. An upper bound to the first term can be obtained using the chaining method. Define,
$$Q_n(g_j) = \frac{1}{\sqrt{k}}\sum_{i\in N_k(x)} W_i g_j(X_i)\text{ for $j=1,2,\ldots, M$.}$$

These variables are mean zero Gaussian, conditionally on the covariates, and the associated semimetric is given by
$$\rho^2(g_j, g_{\ell}) = \text{Var}\big(Q_n(g_j)-Q_n(g_\ell)\big) = \|g_j-g_\ell\|_{2,k}^2.$$
Since $\|g\|_{2,k}\leq \delta$ for all $g\in\mathbb{B}(\delta;\partial\F)$, the coarsest resolution of the chaining can be set to $\delta$ and we can terminate it at the first dyadic resolution below $\varepsilon=\delta^2/2C_1^2$, where every member of our finite set is reconstructed exactly. We now carry out the chaining argument explicitly, with constants.

If $\xi_1,\ldots,\xi_N$ are centred Gaussian variables (not
necessarily independent) with $\mathrm{Var}(\xi_j)\le s^2$, then
\begin{equation}
\label{eq:gauss-max}
\E\max_{j\le N}|\xi_j|\le s\sqrt{2\log(2N)} .
\end{equation}
Indeed, $\max_j|\xi_j|$ is the maximum of the $2N$ variables $\pm\xi_j$, and for $\lambda>0$ Jensen's
inequality gives
$$\exp\{\lambda\E\max_j|\xi_j|\}\le\E\exp\{\lambda\max_j|\xi_j|\}\le\sum_{j}\E(e^{\lambda\xi_j}+e^{-\lambda\xi_j})\le2Ne^{\lambda^2s^2/2},$$
so $\E\max_j|\xi_j|\le\lambda^{-1}\log(2N)+\lambda s^2/2$. The choice $\lambda=\sqrt{2\log(2N)}/s$
gives \eqref{eq:gauss-max}.

Now, put $\varepsilon_m=2^{-m}\delta$ for $m\ge0$, and let $m^*\ge1$ be the smallest
integer with $\varepsilon_{m^*}\le\varepsilon$. For $1\le m\le m^*-1$ let $T_m\subset\mathbb B(\delta;\partial\F)$ be an
internal $\varepsilon_m$-cover of minimal cardinality $\mathcal N(\varepsilon_m)$, let $T_0=\{0\}$ (note
$0\in\partial\F$ and every element of the ball is within $\delta=\varepsilon_0$ of $0$), and let
$T_{m^*}=\{g_1,\ldots,g_M\}$. Because $\varepsilon_{m^*}\le\varepsilon$, an $\varepsilon_{m^*}$-cover is an
$\varepsilon$-cover, so by minimality $M=\mathcal N(\varepsilon)\le\mathcal N(\varepsilon_{m^*})$. For each $j$ define
$\pi_{m^*}(g_j)=g_j$ and, for $0\le m<m^*$, let $\pi_m(g_j)\in T_m$ be a point with
$\|g_j-\pi_m(g_j)\|_{2,k}\le\varepsilon_m$. Since $Q_n$ is linear and $Q_n(0)=0$,
\[
Q_n(g_j)=\sum_{m=1}^{m^*}\Big\{Q_n\big(\pi_m(g_j)\big)-Q_n\big(\pi_{m-1}(g_j)\big)\Big\}.
\]
For the $m$-th link,
$\|\pi_m(g_j)-\pi_{m-1}(g_j)\|_{2,k}\le\varepsilon_m+\varepsilon_{m-1}=3\varepsilon_m$, and as $j$ varies the link
takes at most $|T_m|\,|T_{m-1}|\le\mathcal N(\varepsilon_m)^2$ distinct values, because covering numbers are
non-increasing in the resolution. Applying \eqref{eq:gauss-max} to each level, with
$s=3\varepsilon_m$ and $N=\mathcal N(\varepsilon_m)^2$, and using $\log(2N^2)\le2\log(2N)$,
\[
\E\max_{j\le M}|Q_n(g_j)|
\le\sum_{m=1}^{m^*}3\varepsilon_m\sqrt{2\log\big(2\mathcal N(\varepsilon_m)^2\big)}
\le\sum_{m=1}^{m^*}6\varepsilon_m\sqrt{\log\big(2\mathcal N(\varepsilon_m)\big)} .
\]

Since $\varepsilon_m=2(\varepsilon_m-\varepsilon_{m+1})$ and
$t\mapsto\mathcal N(t)$ is non-increasing,
\[
\varepsilon_m\sqrt{\log\big(2\mathcal N(\varepsilon_m)\big)}
\le2\int_{\varepsilon_{m+1}}^{\varepsilon_m}\sqrt{\log\big(2\mathcal N(t)\big)}\,\mathrm dt .
\]
Summing over $1\le m\le m^*$,
\[
\E\max_{j\le M}|Q_n(g_j)|\le12\int_{\varepsilon_{m^*+1}}^{\delta/2}\sqrt{\log\big(2\mathcal N(t)\big)}\,\mathrm dt
\le12\int_{\varepsilon/4}^{\delta}\sqrt{\log\big(2\mathcal N(t)\big)}\,\mathrm dt,
\]
where the last step uses $\varepsilon_{m^*+1}\ge\varepsilon/4$: if $m^*\ge2$, minimality of $m^*$ gives
$\varepsilon_{m^*-1}>\varepsilon$, hence $\varepsilon_{m^*+1}=\varepsilon_{m^*-1}/4>\varepsilon/4$; if $m^*=1$, then
$\varepsilon_2=\delta/4>\varepsilon/4$. Dividing by $\sqrt k$, since
$k^{-1}\sum_{i\in N_k(x)}W_ig_j(X_i)=Q_n(g_j)/\sqrt k$, and recalling $\varepsilon/4=\delta^2/(8C_1^2)$, we find that
$$\E\left[\max_{j=1,2,\ldots, M} \Big|\frac{1}{k}\sum_{i\in N_k(x)} W_i g_j(X_i)\Big|\right] \leq \frac{12}{\sqrt{k}} \int_{\delta^2/8C_1^2}^\delta \sqrt{\log \Big(2\mathcal{N}\big(t; \mathbb{B}(\delta;\partial\F), \|\cdot\|_{2,k}\big)\Big)}\mathrm d t.$$
Combining this with the first display of the proof gives \eqref{eq:dudley-Gn}. For the
second statement, take $C_1=1$: if $\delta\in(0,1]$ satisfies the displayed condition, then
\eqref{eq:dudley-Gn} gives $G_n(\delta;\partial\F)\le\delta^2/2+\delta^2/2=\delta^2$.

\end{proof}

\begin{proof}[Proof of Theorem~\ref{thm:oracle}]
Recall the definition of $\mathcal{A}(u)$ using $\partial \mathcal{F}(L,\bm p, s)$, i.e.,
$$\mathcal{A}(u) = \left\{\exists g\in \{h\in\partial\mathcal{F}: \|h\|_{2,k}\geq u\}: \Big|\frac{1}{k}\sum_{i\in N_k(x)} W_i g(X_i)\Big|\geq 2\|g\|_{2,k} u\right\}.$$
The argument of Section~\ref{sec:train-error} applies verbatim and gives
$$\P[\mathcal{A}(\sqrt{t\delta_n})^c] \geq 1- \exp\Big\{-\frac{tk\delta_n}{2}\Big\},$$
for all $t\geq \delta_n$, where $\delta_n$ is the smallest positive solution of the critical inequality
$$G_n\big(\delta;\partial \mathcal{F}(L,\bm p, s)\big)=\E\left[\sup_{h\in \partial\mathcal{F}(L,\bm p, s), \|h\|_{2,k}\leq \delta}\Big|\frac{1}{k}\sum_{i\in N_k(x)} W_i h(X_i)\Big|\right] \leq \delta^2.$$
We will assume the event $\mathcal{A}(\sqrt{t\delta_n})^c$ is true. Then, given any $\Tilde{f}\in\mathcal{F}(L,\bm p, s)$ and $\hat f_n$ we have 
\begin{align}
\label{eq:krm-ineq}
    \frac{1}{2k}\sum_{i\in N_k(x)} \big(Y_i- \hat f_n(X_i)\big)^2 \leq \frac{1}{2k}\sum_{i\in N_k(x)} \big(Y_i- \Tilde{f}(X_i)\big)^2.
\end{align}
Also, note that
\begin{align*}
    & \frac{1}{k}\sum_{i\in N_k(x)} \big(Y_i- \hat f_n(X_i)\big)^2\\
    & = \frac{1}{k}\sum_{i\in N_k(x)} \big(Y_i-f_0(X_i) + f_0(X_i) - \hat f_n(X_i)\big)^2\\
    & = \frac{1}{k}\sum_{i\in N_k(x)} \big(Y_i-f_0(X_i)\big)^2 + \frac{1}{k}\sum_{i\in N_k(x)}\big(f_0(X_i) - \hat f_n(X_i)\big)^2 + \frac{2}{k}\sum_{i\in N_k(x)}\big(Y_i-f_0(X_i)\big)\big(f_0(X_i)-\hat f_n(X_i)\big).
\end{align*}
Similarly, we also get
\begin{align*}
    & \frac{1}{k}\sum_{i\in N_k(x)} \big(Y_i-\Tilde{f}(X_i)\big)^2 \\
    & = \frac{1}{k}\sum_{i\in N_k(x)} \big(Y_i-f_0(X_i)\big)^2 + \frac{1}{k}\sum_{i\in N_k(x)}\big(f_0(X_i) - \Tilde{f}(X_i)\big)^2 + \frac{2}{k}\sum_{i\in N_k(x)}\big(Y_i-f_0(X_i)\big)\big(f_0(X_i)- \Tilde{f}(X_i)\big).
\end{align*}
Now using \eqref{eq:krm-ineq} and the fact that $Y_i = f_0(X_i) + \sigma W_i$, we get
\begin{align*}
    \|\hat f_n - f_0\|_{2,k}^2 + \frac{2\sigma}{k}\sum_{i\in N_k(x)}W_i\big(f_0(X_i)-\hat f_n(X_i)\big) \leq \|\Tilde{f}-f_0\|_{2,k}^2 + \frac{2\sigma}{k}\sum_{i\in N_k(x)}W_i\big(f_0(X_i)-\Tilde{f}(X_i)\big),
\end{align*}
which gives us our basic inequality for this case as 
\begin{align}
    \label{eq:basic-ineq-2}
    \|\hat f_n - f_0\|_{2,k}^2 \leq \|\Tilde{f}-f_0\|_{2,k}^2 + \frac{2\sigma}{k}\sum_{i\in N_k(x)}W_i\big(\hat{f}_n(X_i)-\Tilde{f}(X_i)\big).
\end{align}
Clearly, to analyze the convergence rate of $\|\hat f_n-f_0\|_{2,k}^2$ we need to analyze the term 
$$\frac{2\sigma}{k}\sum_{i\in N_k(x)}W_i\big(\hat{f}_n(X_i)-\Tilde{f}(X_i)\big) = \frac{2\sigma}{k}\sum_{i\in N_k(x)}W_i\Tilde{\Delta}(X_i)$$
where $\Tilde{\Delta}(x) = \big(\hat f_n(x)-\Tilde{f}(x)\big)$. We also define $\hat \Delta(x) = \hat f_n(x) - f_0(x)$ for notational convenience. We break our analysis into two cases.

\noindent\textbf{Case 1:} First suppose that $\|\Tilde{\Delta}\|_{2,k}\leq \sqrt{t\delta_n}$. Then,
\begin{align*}
    \|\hat{\Delta}\|_{2,k}^2 = \|\hat f_n - f_0\|_{2,k}^2 & = \|(\Tilde{f}-f_0) + \Tilde{\Delta}\|_{2,k}^2\\
    & \stackrel{(a)}{\leq} \Big\{\|\Tilde{f}-f_0\|_{2,k} + \sqrt{t\delta_n}\Big\}^2\\
    & \stackrel{(b)}{\leq} 2\|\Tilde{f}-f_0\|_{2,k}^2 + 2t\delta_n,
\end{align*}
where in (a) we have used the triangle inequality and (b) we used the inequality $(a+b)^2\leq 2(a^2+b^2)$.

\noindent\textbf{Case 2:} Suppose now that $\|\Tilde{\Delta}\|_{2,k}>\sqrt{t\delta_n}$. Note that $\Tilde{\Delta}\in\partial\mathcal{F}$ and the event $\mathcal{A}(\sqrt{t\delta_n})^c$ is true. Hence,
$$\Big|\frac{1}{k}\sum_{i\in N_k(x)} W_i \Tilde{\Delta}(X_i)\Big|\leq 2\sqrt{t\delta_n} \|\Tilde{\Delta}\|_{2,k}.$$
Combining with \eqref{eq:basic-ineq-2}, we find that with probability at least $1-\exp\big\{-\frac{kt\delta_n}{2}\big\}$, the squared error is bounded as
\begin{align*}
    \|\hat \Delta\|_{2,k}^2 & \leq \|\Tilde{f}-f_0\|_{2,k}^2 + 4\sigma \sqrt{t\delta_n} \|\Tilde{\Delta}\|_{2,k}\\
    & \leq \|\Tilde{f}-f_0\|_{2,k}^2 + 4\sigma \sqrt{t\delta_n} \big\{\|\hat{\Delta}\|_{2,k} + \|\Tilde{f} - f_0\|_{2,k}\big\}\\
    & \leq \|\Tilde{f}-f_0\|_{2,k}^2 + 2\sigma\Big[\frac{t\delta_n}{\beta}+ \beta\|\hat{\Delta}\|_{2,k}^2\Big] + 2\sigma\Big[\frac{t\delta_n}{\beta} + \beta\|\Tilde{f}-f_0\|_{2,k}^2\Big],
\end{align*}
where the first inequality follows from the triangle inequality and the second inequality follows from the fact $2ab\leq \beta a^2 + b^2/\beta$ for $a,b\in\mathbb{R}$ and $\beta>0$. Taking $\beta = (6\sigma)^{-1}$, we have
\begin{align*}
    & \|\hat \Delta\|_{2,k}^2 \leq \|\Tilde{f}-f_0\|_{2,k}^2 + 2\Big[6\sigma^2t\delta_n+ \frac{1}{6}\|\hat{\Delta}\|_{2,k}^2\Big] + 2\Big[6\sigma^2t\delta_n+ \frac{1}{6}\|\Tilde{f}-f_0\|_{2,k}^2\Big]\\
    & \Rightarrow \frac{2}{3} \|\hat \Delta\|_{2,k}^2 \leq \frac{4}{3} \|\Tilde{f}-f_0\|_{2,k}^2 + 24 \sigma^2t\delta_n\\
    & \Rightarrow \|\hat \Delta\|_{2,k}^2 \leq 2 \|\Tilde{f}-f_0\|_{2,k}^2 + 36 \sigma^2 t\delta_n.
\end{align*}
Since $\Tilde{f}\in\mathcal{F}$ is arbitrary, in both cases we get
\begin{align}
\|\hat \Delta\|_{2,k}^2 \leq 2 \inf_{f\in\mathcal{F}(L,\bm p, s)}\|f-f_0\|_{2,k}^2 + Q_\sigma t\delta_n,
\end{align}
holds with probability at least $1-\exp\big\{-\frac{kt\delta_n}{2}\big\}$, where $Q_\sigma = \max\{2,36\sigma^2\}$. By integrating the tail bounds we can also show that
\begin{align}
\E\left[\frac{1}{k}\sum_{i\in N_k(x)} \big(\hat f_n(X_i) - f_0(X_i)\big)^2\ \Big|\ X_1,\ldots,X_n\right] \leq 2\inf_{f\in\mathcal{F}(L,\bm p, s)}\|f-f_0\|_{2,k}^2 + Q_\sigma\delta_n^2 + \frac{2Q_\sigma}{k}.
\end{align}
Both \eqref{eq:oracle-tail} and \eqref{eq:oracle-expectation} are statements
conditional on $X_1,\ldots,X_n$; note in particular that $\delta_n$ is a function of the
covariates, since the norm $\|\cdot\|_{2,k}$ and the index set $N_k(x)$ are.
\end{proof}

\begin{proof}[Proof of Proposition~\ref{prop:approx}]
We first control the approximation error induced by restricting the
estimation procedure to the $k$-nearest neighbourhood of a
fixed query point $x$. Let
\[
R_{n,k}(x)
=
\max_{i\in N_k(x)} \|X_i-x\|
=
\|X_{(k,n)}(x)-x\|.
\]
Recall that $\mathcal F_n=\mathcal{F}(L,\bm p,s)$ contains all constant functions
bounded by $B$. In particular, whenever $|f_0(x)|\leq B$, the constant
function
\[
f_x(\cdot)\equiv f_0(x)
\]
belongs to $\F_n$. Therefore, for every realization of
$X_1,\ldots,X_n$,
\[
\begin{aligned}
\inf_{f\in \F_n}
\frac{1}{k}
\sum_{i\in N_k(x)}
\bigl\{f(X_i)-f_0(X_i)\bigr\}^2
&\leq
\frac{1}{k}
\sum_{i\in N_k(x)}
\bigl\{f_0(x)-f_0(X_i)\bigr\}^2 .
\end{aligned}
\]
Assume that $f_0$ is locally $\alpha$-H\"older continuous at $x$ with $0<\alpha\leq 1$, that is, there exists an $r_0>0$ such that
\[
|f_0(u)-f_0(x)|
\leq
L_0\|u-x\|^\alpha,\text{ for all $u\in B(x,r_0)$}.
\]
Decreasing $r_0$ preserves the H\"older assumption, so we may and do assume that $r_0$ is
not larger than the constant $r_0$ of Lemma~\ref{lem:local_mass}; then the condition
$k\le\frac12nc_0r_0^{d^*}$ of Lemma~\ref{lem:knn_radius_tail} holds for this $r_0$ as well.
On the event $\mathcal R = \{R_{n,k}(x)\le r_0\}$ every $X_i$, $i\in N_k(x)$, lies in $B(x,r_0)$ and
\[
\inf_{f\in\mathcal F_n}\frac1k\sum_{i\in N_k(x)}\{f(X_i)-f_0(X_i)\}^2
\le\frac{L_0^2}{k}\sum_{i\in N_k(x)}\|X_i-x\|^{2\alpha}\le L_0^2R_{n,k}(x)^{2\alpha}.
\]
On $\mathcal R^c$ the left-hand side is at most $(2B)^2$, so
\[
\E\left[\inf_{f\in\mathcal F_n}\frac1k\sum_{i\in N_k(x)}\{f(X_i)-f_0(X_i)\}^2\right]
\le L_0^2\,\E\big[R_{n,k}(x)^{2\alpha}\big]+4B^2\,\P\big(R_{n,k}(x)>r_0\big),
\]
and by Lemma~\ref{lem:knn_radius_tail}, applied with $r=r_0$ (admissible because
$k\le\frac12nc_0r_0^{d^*}$ is equivalent to $r_{n,k}\le r_0$), the second term is at most
$4B^2\exp\{-\tfrac18nc_0r_0^{d^*}\}$, which is exponentially small in $n$. Hence
\[
\E\left[\inf_{f\in\mathcal F_n}\frac1k\sum_{i\in N_k(x)}\{f(X_i)-f_0(X_i)\}^2\right]
\le L_0^2\,\E\big[R_{n,k}(x)^{2\alpha}\big]+4B^2\exp\Big\{-\frac18nc_0r_0^{d^*}\Big\}.
\]
Since $\sup_{n}n^{2\alpha/d^*}\exp\{-\frac18nc_0r_0^{d^*}\}<\infty$ and $k\ge1$, the exponential term is
bounded by a constant multiple of $(k/n)^{2\alpha/d^*}$, which is the order of the first term by
Lemma~\ref{lem:expected_knn_radius}; in this sense the bound is
$\lesssim L_0^2\E[R_{n,k}(x)^{2\alpha}]+(k/n)^{2\alpha/d^*}$.
\end{proof}

\begin{proof}[Proof of Lemma~\ref{lem:local_mass}]
Let $\tau>0$ be the lower bound on the reach of $M$, and let $X$ have a
density $p_X$ with respect to the volume measure $\operatorname{Vol}_M$,
with $p_{X,0}\le p_X\le p_{X,1}$ on $M$. Since $M$ is a compact
$d^*$-dimensional submanifold with reach at least $\tau$, there are
constants $0<a_0\le a_1<\infty$, depending only on $d^*$, such that
\begin{equation}
\label{eq:vol-two-sided}
a_0\,v_{d^*}r^{d^*}\le\operatorname{Vol}_M\bigl(B(x,r)\cap M\bigr)
\le a_1\,v_{d^*}r^{d^*},
\qquad x\in M,\ 0<r\le\tau/4;
\end{equation}
see, for example, \citet[Lemma~5.3]{niyogi2008finding} for the lower
bound, which gives
$\operatorname{Vol}_M(B(x,r)\cap M)\ge(1-r^2/(4\tau^2))^{d^*/2}
v_{d^*}r^{d^*}\ge(15/16)^{d^*/2}v_{d^*}r^{d^*}$ for $r\le\tau/2$,
and \citet[Lemma~B.7]{aamari2019nonasymptotic} for the two-sided bound, which holds for
$0<r\le\tau/4$ with constants depending only on $d^*$. Put
$r_0=\min\{\tau/4,\operatorname{diam}(M)\}$. For $0<r\le r_0$,
\[
\P(\|X-x\|\le r)=\int_{B(x,r)\cap M}p_X\,d\operatorname{Vol}_M
\]
lies between $p_{X,0}\operatorname{Vol}_M(B(x,r)\cap M)$ and
$p_{X,1}\operatorname{Vol}_M(B(x,r)\cap M)$, and
\eqref{eq:vol-two-sided} gives the lemma with
$c_0=a_0p_{X,0}v_{d^*}$ and $c_1=a_1p_{X,1}v_{d^*}$, both the lower and the upper bound
stated in Lemma~\ref{lem:local_mass}.

\end{proof}

\begin{proof}[Proof of Lemma~\ref{lem:knn_radius_tail}]
For fixed $x$ and $r>0$,
\[
N_x(r)
=
\sum_{i=1}^n
\mathbf{1}\{\|X_i-x\|\leq r\}
\sim
\operatorname{Binomial}(n,q_x(r)),
\]
where
\[
q_x(r)
=
\mathbb{P}\bigl(\|X-x\|\leq r\bigr).
\]
Moreover,
\[
\{R_{n,k}(x)>r\}
=
\{N_x(r)<k\}.
\]

For $r\leq r_0$, Lemma~\ref{lem:local_mass} gives
\[
q_x(r)\geq c_0r^{d^*}.
\]
Hence
\[
\mu_x(r):=\mathbb{E}N_x(r)
=nq_x(r)
\geq nc_0r^{d^*}.
\]
If $r\geq r_{n,k}$, then
\[
nc_0r^{d^*}\geq 2k,
\]
and therefore
\[
k\leq \frac{\mu_x(r)}{2}.
\]
Consequently,
\[
\mathbb{P}\bigl(R_{n,k}(x)>r\bigr)
=
\mathbb{P}\bigl(N_x(r)<k\bigr)
\leq
\mathbb{P}\left(
N_x(r)\leq\frac{\mu_x(r)}{2}
\right).
\]

The standard multiplicative Chernoff bound for a binomial random
variable gives
\[
\mathbb{P}\left(
N_x(r)\leq\frac{\mu_x(r)}{2}
\right)
\leq
\exp\left(-\frac{\mu_x(r)}{8}\right).
\]
Since $\mu_x(r)\geq nc_0r^{d^*}$,
\[
\mathbb{P}\bigl(R_{n,k}(x)>r\bigr)
\leq
\exp\left\{-\frac{c_0}{8}nr^{d^*}\right\}.
\]
\end{proof}

\begin{proof}[Proof of Lemma~\ref{lem:expected_knn_radius}]
Since $M$ is compact, there exists a finite constant
\[
D_M=\sup_{u,v\in M}\|u-v\|<\infty.
\]
Consequently,
\[
0\leq R_{n,k}(x)\leq D_M.
\]
Using the tail-integral representation,
\[
\mathbb{E}\left[R_{n,k}(x)^{2\alpha}\right]
=
2\alpha
\int_0^{D_M}
r^{2\alpha-1}
\mathbb{P}\bigl(R_{n,k}(x)>r\bigr)\,dr.
\]
Put
\[
r_{n,k}
=
\left(\frac{2k}{c_0n}\right)^{1/d^*}.
\]
We divide the integral into three regions:
\[
[0,r_{n,k}],
\qquad
[r_{n,k},r_0],
\qquad
[r_0,D_M].
\]
Thus,
\[
\mathbb{E}\left[R_{n,k}(x)^{2\alpha}\right]
=I_1+I_2+I_3,
\]
where
\[
\begin{aligned}
I_1
&=
2\alpha\int_0^{r_{n,k}}
r^{2\alpha-1}
\mathbb{P}\bigl(R_{n,k}(x)>r\bigr)\,dr,\\
I_2
&=
2\alpha\int_{r_{n,k}}^{r_0}
r^{2\alpha-1}
\mathbb{P}\bigl(R_{n,k}(x)>r\bigr)\,dr,\\
I_3
&=
2\alpha\int_{r_0}^{D_M}
r^{2\alpha-1}
\mathbb{P}\bigl(R_{n,k}(x)>r\bigr)\,dr.
\end{aligned}
\]

For the first term we simply use
$\mathbb{P}(R_{n,k}(x)>r)\leq 1$:
\[
I_1
\leq
2\alpha\int_0^{r_{n,k}}r^{2\alpha-1}\,dr
=
r_{n,k}^{2\alpha}.
\]
Therefore,
\[
I_1
\leq
\left(\frac{2k}{c_0n}\right)^{2\alpha/d^*}.
\]

For the second term, Lemma~\ref{lem:knn_radius_tail} gives
\[
\begin{aligned}
I_2
&\leq
2\alpha
\int_{r_{n,k}}^{r_0}
r^{2\alpha-1}
\exp\left\{-\frac{c_0}{8}nr^{d^*}\right\}\,dr.
\end{aligned}
\]
Using the change of variables
\[
u=\frac{c_0}{8}nr^{d^*},
\]
we have
\[
r
=
\left(\frac{8u}{c_0n}\right)^{1/d^*},
\]
and hence
\[
dr
=
\frac1{d^*}
\left(\frac{8}{c_0n}\right)^{1/d^*}
u^{1/d^*-1}\,du.
\]
Therefore,
\[
\begin{aligned}
I_2
&\leq
C n^{-2\alpha/d^*}
\int_{(c_0/8)nr_{n,k}^{d^*}}^\infty
u^{2\alpha/d^*-1}e^{-u}\,du\\
&\leq
C^\prime n^{-2\alpha/d^*},
\end{aligned}
\]
where $C$ depends only on $\alpha,d^*$ and $c_0$.
Since $k\geq1$,
\[
n^{-2\alpha/d^*}
\leq
\left(\frac{k}{n}\right)^{2\alpha/d^*}.
\]
Thus,
\[
I_2
\lesssim
\left(\frac{k}{n}\right)^{2\alpha/d^*}.
\]

It remains to control $I_3$. Since $r_0$ is fixed and
\[
q_x(r_0)
\geq c_0r_0^{d^*},
\]
we have
\[
N_x(r_0)
\sim
\operatorname{Binomial}(n,q_x(r_0)),
\]
with mean satisfying
\[
\mathbb{E}N_x(r_0)
\geq nc_0r_0^{d^*}
\geq 2k.
\]
Consequently, another application of the Chernoff bound yields
\[
\mathbb{P}\bigl(R_{n,k}(x)>r_0\bigr)
\leq
\exp\left\{-\frac18nc_0r_0^{d^*}\right\}.
\]
Since the event $\{R_{n,k}(x)>r\}$ is decreasing in $r$,
\[
\mathbb{P}\bigl(R_{n,k}(x)>r\bigr)
\leq
\exp\left\{-\frac18nc_0r_0^{d^*}\right\},
\qquad r\geq r_0.
\]
Therefore,
\[
\begin{aligned}
I_3
&\leq
2\alpha
\exp\left\{-\frac18nc_0r_0^{d^*}\right\}
\int_{r_0}^{D_M}r^{2\alpha-1}\,dr\\
&\leq
D_M^{2\alpha}
\exp\left\{-\frac18nc_0r_0^{d^*}\right\}.
\end{aligned}
\]
The final quantity is exponentially small in $n$ and, in particular,
\[
I_3
\lesssim
\left(\frac{k}{n}\right)^{2\alpha/d^*}
\]
for all sufficiently large $n$.

Combining the three terms gives
\[
\mathbb{E}\left[R_{n,k}(x)^{2\alpha}\right]
\lesssim
\left(\frac{k}{n}\right)^{2\alpha/d^*}.
\]
The constants are uniform in $x\in M$.
\end{proof}

\begin{proof}[Proof of Proposition~\ref{prop:pointwise}]
Since $\widehat f_n\in \mathcal F_n$, it follows that
\[
\begin{aligned}
\big|
(\widehat f_n-f_0)(u)
-
(\widehat f_n-f_0)(v)
\big|
&\leq
|\widehat f_n(u)-\widehat f_n(v)|
+
|f_0(u)-f_0(v)|\\
&\leq
(L_\F+L_0)\|u-v\|^\alpha.
\end{aligned}
\]
Thus, setting
\[
L_H=L_\F+L_0,
\]
the function $\widehat f_n-f_0$ is $\alpha$-H\"older continuous with
constant $L_H$.
For every $i\in N_k(x)$, we therefore have
\[
\begin{aligned}
\big|
\widehat f_n(x)-f_0(x)
\big|
&\leq
\big|
\widehat f_n(X_i)-f_0(X_i)
\big|
+
\big|
(\widehat f_n-f_0)(x)
-
(\widehat f_n-f_0)(X_i)
\big|\\
&\leq
\big|
\widehat f_n(X_i)-f_0(X_i)
\big|
+
L_H\|X_i-x\|^\alpha\\
&\leq
\big|
\widehat f_n(X_i)-f_0(X_i)
\big|
+
L_H R_{n,k}(x)^\alpha.
\end{aligned}
\]

Since this holds for every $i\in N_k(x)$, averaging over the
neighbourhood gives
\[
\big|
\widehat f_n(x)-f_0(x)
\big|
\leq
\frac{1}{k}
\sum_{i\in N_k(x)}
\big|
\widehat f_n(X_i)-f_0(X_i)
\big|
+
L_HR_{n,k}(x)^\alpha.
\]
By the Cauchy--Schwarz inequality,
\[
\frac{1}{k}
\sum_{i\in N_k(x)}
\big|
\widehat f_n(X_i)-f_0(X_i)
\big|
\leq
\left[
\frac{1}{k}
\sum_{i\in N_k(x)}
\bigl\{
\widehat f_n(X_i)-f_0(X_i)
\bigr\}^2
\right]^{1/2}.
\]
Consequently,
\[
\big|
\widehat f_n(x)-f_0(x)
\big|
\leq
T_{n,k}(x)^{1/2}
+
L_HR_{n,k}(x)^\alpha.
\]
Squaring both sides and using
\[
(a+b)^2\leq 2a^2+2b^2,
\]
we obtain the deterministic inequality
\[
\big|
\widehat f_n(x)-f_0(x)
\big|^2
\leq
2T_{n,k}(x)
+
2L_H^2R_{n,k}(x)^{2\alpha}.
\]

This inequality shows explicitly how the pointwise testing error is
controlled by two quantities: the local training error and the radius
of the local neighbourhood.

Taking expectation on both sides yields
\[
\mathbb{E}
\big[
\big|
\widehat f_n(x)-f_0(x)
\big|^2
\big]
\leq
2\mathbb{E}\big[T_{n,k}(x)\big]
+
2L_H^2
\mathbb{E}\big[
R_{n,k}(x)^{2\alpha}
\big].
\]
By Lemma~\ref{lem:expected_knn_radius},
\[
\mathbb{E}\big[
R_{n,k}(x)^{2\alpha}
\big]
\lesssim
\left(\frac{k}{n}\right)^{2\alpha/d^*},
\]
and the bound on the expected approximation error,
\[
\mathbb{E}
\left[
\inf_{f\in \F}
\frac{1}{k}
\sum_{i\in N_k(x)}
\bigl\{f(X_i)-f_0(X_i)\bigr\}^2
\right]
\lesssim L_0^2 \E[R_{n,k}(x)^{2\alpha}] \lesssim
\left(\frac{k}{n}\right)^{2\alpha/d^*},
\]
we get
\begin{align}
    \mathbb{E}
\big[
\big|
\widehat f_n(x)-f_0(x)
\big|^2
\big]
& \lesssim
\mathbb{E}\big[T_{n,k}(x)\big]
+
\left(\frac{k}{n}\right)^{2\alpha/d^*}\\
&\lesssim
\mathbb{E}
\left[
\inf_{f\in \F}
\frac{1}{k}
\sum_{i\in N_k(x)}
\bigl\{f(X_i)-f_0(X_i)\bigr\}^2
\right]
+ Q_\sigma\E\big[\delta_n^2\big] + \frac{2Q_\sigma}{k}+
\left(\frac{k}{n}\right)^{2\alpha/d^*}\\
&\lesssim \left(\frac{k}{n}\right)^{2\alpha/d^*}
+
Q_\sigma\E\big[\delta_n^2\big]
+
\frac{2Q_\sigma}{k}.
\end{align}
In the second step we have taken the expectation over the covariates of the conditional
bound \eqref{eq:oracle-expectation}, using the tower property
$\E[T_{n,k}(x)]=\E\{\E[T_{n,k}(x)\mid X_1,\ldots,X_n]\}$. This turns the covariate-dependent
$\delta_n^2$ into $\E[\delta_n^2]$, and the factor $2$ in front of the approximation term of
\eqref{eq:oracle-expectation} is absorbed into the constant. The approximation term is
bounded by Corollary~\ref{cor:approx-rate}. By Proposition~\ref{prop:critical-radius},
$\E[\delta_n^2]\le\bar\delta^2$, a deterministic quantity.
\end{proof}

\begin{proof}[Proof of Proposition~\ref{prop:critical-radius}]
We first derive a covering-number bound for this class. For any
seminorm $\|\cdot\|$, the difference class satisfies
\[
\mathcal{N}
\left(
\epsilon;
\partial\mathcal{F}_n,
\|\cdot\|
\right)
\leq
\mathcal{N}
\left(
\frac{\epsilon}{2};
\mathcal{F}_n,
\|\cdot\|
\right)^2.
\]
Indeed, let
$\{f_1,\ldots,f_M\}$ be an $\epsilon/2$-cover of
$\mathcal{F}_n$. For any $f,g\in\mathcal{F}_n$, choose $f_j,f_\ell$
such that
\[
\|f-f_j\|\leq\frac{\epsilon}{2},
\qquad
\|g-f_\ell\|\leq\frac{\epsilon}{2}.
\]
Then
\[
\|(f-g)-(f_j-f_\ell)\|
\leq
\|f-f_j\|+\|g-f_\ell\|
\leq\epsilon.
\]
Thus, taking $\|\cdot\|=\|\cdot\|_{2,k}$,
\[
\mathcal{N}
\left(
\epsilon;
\partial\mathcal{F}_n,
\|\cdot\|_{2,k}
\right)
\leq
\mathcal{N}
\left(
\frac{\epsilon}{2};
\mathcal{F}_n,
\|\cdot\|_{2,k}
\right)^2.
\label{eq:diffclass-cover}
\]

Moreover,
\[
\|f-g\|_{2,k}
\leq
\|f-g\|_\infty,
\]
so every $\|\cdot\|_\infty$-cover is also a
$\|\cdot\|_{2,k}$-cover.  Lemma 5 of \cite{dnnregression} gives
\[
\log
\mathcal{N}
\left(
\epsilon;
\mathcal{F}_n,
\|\cdot\|_\infty
\right)
\leq
(s+1)
\log
\left(
\frac{2(L+1)V^2}{\epsilon}
\right),
\]
where
\[
V=\prod_{\ell=0}^{L+1}(p_\ell+1).
\]
Consequently,
\begin{equation}
\begin{aligned}
\log
\mathcal{N}
\left(
\epsilon;
\partial\mathcal{F}_n,
\|\cdot\|_{2,k}
\right)
&\leq
2\log
\mathcal{N}
\left(
\frac{\epsilon}{2};
\mathcal{F}_n,
\|\cdot\|_{2,k}
\right)\\
&\leq
2(s+1)
\log
\left(
\frac{4(L+1)V^2}{\epsilon}
\right).
\end{aligned}
\label{eq:sh-entropy}
\end{equation}

The covers in the last two displays are external: their centres $f_j-f_\ell$ need not lie
in $\partial\F_n$, and the centres of a $\|\cdot\|_\infty$-cover of $\mathcal F_n$ need not lie in
$\mathcal F_n$ (Lemma~5 of \cite{dnnregression} is stated for the larger class without the sup-norm
constraint, and an external cover of a class is an external cover of every subclass).
Lemma~\ref{lem:dudley}, however, uses the internal covering number of the ball
$\mathbb B(\delta;\partial\F_n)$. The two are related by the elementary inequality
\[
\mathcal N_{\mathrm{int}}(t;\mathcal T,\|\cdot\|)\le\mathcal N_{\mathrm{ext}}(t/2;\mathcal T,\|\cdot\|),
\]
valid for any set $\mathcal T$: given an external $t/2$-cover with centres $c_1,\ldots,c_N$, discard the
centres whose $t/2$-ball misses $\mathcal T$ and replace each remaining $c_j$ by a point
$h_j\in\mathcal T$ with $\|h_j-c_j\|\le t/2$; by the triangle inequality the $h_j$ form an internal
$t$-cover. Since $\mathbb B(\delta;\partial\F_n)\subseteq\partial\F_n$, every external cover of
$\partial\F_n$ is an external cover of the ball. Hence, by \eqref{eq:sh-entropy} with
$\epsilon=t/2$,
\[
\log\mathcal N\big(t;\mathbb B(\delta;\partial\F_n),\|\cdot\|_{2,k}\big)
\le\log\mathcal N_{\mathrm{ext}}\Big(\frac t2;\partial\F_n,\|\cdot\|_{2,k}\Big)
\le2(s+1)\log\Big(\frac{8(L+1)V^2}{t}\Big),
\]
and, since $\log2\le2(s+1)\log2$,
\begin{equation}
\label{eq:sh-entropy-int}
\log\Big(2\mathcal N\big(t;\mathbb B(\delta;\partial\F_n),\|\cdot\|_{2,k}\big)\Big)
\le2(s+1)\log\Big(\frac{16(L+1)V^2}{t}\Big),\qquad 0<t\le\delta\le1 .
\end{equation}
The right-hand side is positive, because $V\ge1$ and $t\le1$.

We now substitute this entropy bound into the Dudley bound \eqref{eq:dudley-Gn} for the
localized Gaussian complexity. Take first $\delta\in(0,1]$ and $C_1=1$. Then, by
\eqref{eq:sh-entropy-int},
\[
G_n(\delta;\partial\mathcal{F}_n)
\leq
\frac{12\sqrt{2(s+1)}}{\sqrt{k}}
\int_{\delta^2/8}^{\delta}
\sqrt{
\log
\left(
\frac{16(L+1)V^2}{t}
\right)
}
\,dt
+
\frac{\delta^2}{2}.
\]
Since $t\geq\delta^2/8$ throughout the integration interval,
\[
\log
\left(
\frac{16(L+1)V^2}{t}
\right)
\leq
\log
\left(
\frac{128(L+1)V^2}{\delta^2}
\right),
\]
and the interval has length at most $\delta$. Therefore,
\[
G_n(\delta;\partial\mathcal{F}_n)
\leq
\frac{12\sqrt{2(s+1)}}{\sqrt{k}}
\,\delta
\sqrt{
\log
\left(
\frac{128(L+1)V^2}{\delta^2}
\right)
}
+
\frac{\delta^2}{2}.
\]
It follows that a sufficient condition for the critical inequality
$G_n(\delta;\partial\mathcal{F}_n)\le\delta^2$, for $\delta\in(0,1]$, is
\[
\frac{12\sqrt{2(s+1)}}{\sqrt{k}}
\,\delta
\sqrt{
\log
\left(
\frac{128(L+1)V^2}{\delta^2}
\right)
}
\leq
\frac{\delta^2}{2},
\qquad\text{equivalently}\qquad
\delta^2
\geq
\frac{1152(s+1)}{k}
\log
\left(
\frac{128(L+1)V^2}{\delta^2}
\right).
\]
Since $1152\le2048$, this is implied by \eqref{eq:critical-ineq}; we keep the constant $2048$,
which is the one used in Section~\ref{sec:adaptive}. If $\delta>1$, the critical inequality holds
trivially, because $G_n(\delta;\partial\F_n)\le\delta<\delta^2$ by the Cauchy--Schwarz bound. Hence
every $\delta>0$ satisfying \eqref{eq:critical-ineq} satisfies the critical inequality.

Let
\[
A_n:=\frac{2048(s+1)}{k},
\qquad
D_n:=128(L+1)V^2,
\qquad
\phi(y)=y-A_n\log\Big(\frac{D_n}{y}\Big),\quad y>0 .
\]
The function $\phi$ is continuous and strictly increasing, with $\phi(y)\to-\infty$ as $y\downarrow0$
and $\phi(y)\to\infty$ as $y\to\infty$. Hence the set of $\delta$ satisfying \eqref{eq:critical-ineq}
is $\{\delta:\delta^2\ge\bar y\}$, where $\bar y$ is the unique root of $\phi$, i.e.\ of
\[
y=A_n\log\left(\frac{D_n}{y}\right).
\]
Exponentiating, $e^{y/A_n}=D_n/y$, that is,
\[
\frac{y}{A_n}
\exp\left(\frac{y}{A_n}\right)
=
\frac{D_n}{A_n},
\]
and therefore, by the definition of the Lambert $W$ function as the inverse of $w\mapsto we^w$ on
$[0,\infty)$,
\[
\bar\delta^2=\bar y
=
A_n
W\left(\frac{D_n}{A_n}\right)
=
\frac{2048(s+1)}{k}
W\left(
\frac{k(L+1)V^2}{16(s+1)}
\right).
\]
Since $\bar\delta$ satisfies \eqref{eq:critical-ineq}, it satisfies the critical inequality, and
because $\delta_n$ is the smallest radius at which the critical inequality holds,
$\delta_n\le\bar\delta$. Note that $\bar\delta$ is in general strictly larger than $\delta_n$, so
\eqref{eq:delta-lambertW} defines an upper bound, not the critical radius itself.

Finally, the bound
\[
W(z)\le\log(1+z),
\qquad z\ge0,
\]
holds because $w\mapsto we^w$ is increasing and, with $w=\log(1+z)$,
$we^w=(1+z)\log(1+z)\ge z$ (the function $(1+z)\log(1+z)-z$ vanishes at $0$ and has derivative
$\log(1+z)\ge0$). Hence
\[
\delta_n^2\le\bar\delta^2
\le
\frac{2048(s+1)}{k}
\log
\left(
1+
\frac{k(L+1)V^2}{16(s+1)}
\right)
\le
\frac{2048(s+1)}{k}
\log
\left(
1+
\frac{k(L+1)V^2}{s+1}
\right)
=\frac{2048\,\mathcal C_n}{k}.
\]
The right-hand side does not depend on the covariates, so the same bound holds for
$\E[\delta_n^2]$.
\end{proof}

\begin{proof}[Proof of Proposition~\ref{prop:lower}]
Fix a point $x\in\mathcal M$.
The construction places a bump of height of order $h^\alpha$ and width $h$ at $x$. One may
build it in a local coordinate chart $\phi_x$ of $\mathcal M$ around $x$, but then the H\"older
condition, which in $\mathcal H^\alpha(L_0,B)$ is stated in the ambient Euclidean distance on
$\mathcal M$, has to be transferred through the chart. It is simpler, and gives explicit
constants, to use the ambient distance directly. Let
\[
\psi(t)=(1-t)_+,\qquad t\ge0,
\]
so that $\psi(0)=1$, $\psi(t)=0$ for $t\ge1$, $0\le\psi\le1$, and $\psi$ is $1$-Lipschitz. Put
$a=\min\{L_0,B\}$ and, for $0<h\le1$, define
\[
f_0^{(0)}(z)=0,
\qquad
f_0^{(1)}(z)
=
a h^\alpha
\psi\left(\frac{\|z-x\|}{h}\right),
\qquad z\in\mathcal M .
\]

Both functions lie in $\mathcal H^\alpha(L_0,B)$. This is clear for $f_0^{(0)}$. For
$f_0^{(1)}$, first $\|f_0^{(1)}\|_\infty=ah^\alpha\le a\le B$. Next, let $u,v\in\mathcal M$. Since
$\psi$ is $1$-Lipschitz and $|\|u-x\|-\|v-x\||\le\|u-v\|$,
\[
|f_0^{(1)}(u)-f_0^{(1)}(v)|\le ah^{\alpha-1}\|u-v\| .
\]
If $\|u-v\|\le h$, then $h^{\alpha-1}\|u-v\|=\|u-v\|^\alpha(\|u-v\|/h)^{1-\alpha}\le\|u-v\|^\alpha$,
because $\alpha\le1$. If $\|u-v\|>h$, then simply
$|f_0^{(1)}(u)-f_0^{(1)}(v)|\le\|f_0^{(1)}\|_\infty\le ah^\alpha<a\|u-v\|^\alpha$. In both cases
$|f_0^{(1)}(u)-f_0^{(1)}(v)|\le a\|u-v\|^\alpha\le L_0\|u-v\|^\alpha$. Moreover,
\[
f_0^{(0)}(x)=0,
\qquad
f_0^{(1)}(x)=ah^\alpha.
\]

Under the model
\[
Y_i=f_0(X_i)+\sigma W_i,
\]
the conditional distribution of $Y_i$ given $X_i$ under $f_0^{(j)}$ is
\(
N(f_0^{(j)}(X_i),\sigma^2).
\) The corresponding joint distributions of $(X_i,Y_i)_{i\le n}$ are denoted by $P_1^n$ and $P_0^n$.
Since the covariate distribution is the same under both hypotheses, the chain rule for the
Kullback--Leibler divergence and $D_{\mathrm{KL}}(N(\mu_1,\sigma^2)\|N(\mu_0,\sigma^2))=(\mu_1-\mu_0)^2/(2\sigma^2)$ give
\[
D_{\mathrm{KL}}(P_1^n\|P_0^n)
=
\frac{n}{2\sigma^2}
\mathbb E
\left[
\{f_0^{(1)}(X)-f_0^{(0)}(X)\}^2
\right].
\]
The function $f_0^{(1)}$ vanishes outside $B(x,h)$ and is bounded by $ah^\alpha$. By the upper
bound in Lemma~\ref{lem:local_mass}, for $h\le r_0$,
$\P(\|X-x\|\le h)\le c_1h^{d^*}$, where $c_1$ is the constant of that lemma. Consequently,
\[
\mathbb E[f_0^{(1)}(X)^2]
\leq
{a}^2 h^{2\alpha}\,\P(\|X-x\|\le h)\le a^2c_1h^{2\alpha+d^*},
\]
and hence
\[
D_{\mathrm{KL}}(P_1^n\|P_0^n)
\le
\frac{a^2c_1}{2\sigma^2}\,nh^{2\alpha+d^*}.
\]
Choose
\[
h
=
\Big(\frac{\sigma^2}{a^2c_1\,n}\Big)^{1/(2\alpha+d^*)},
\]
which satisfies $h\le\min\{1,r_0\}$ for all sufficiently large $n$; then
\[
D_{\mathrm{KL}}(P_1^n\|P_0^n)
\le \frac12.
\]
By Pinsker's inequality,
\[
\|P_1^n-P_0^n\|_{\mathrm{TV}}
\le
\sqrt{
\frac12
D_{\mathrm{KL}}(P_1^n\|P_0^n)}\leq \frac12,
\]
and hence the total variation distance is bounded away from one.
We now derive Le Cam's two-point bound with its constant. Let $\theta_j=f_0^{(j)}(x)$ and
$\Delta=|\theta_1-\theta_0|=ah^\alpha$, and let $\widetilde f_n(x)$ be any estimator. Put
$E=\{|\widetilde f_n(x)-\theta_0|<\Delta/2\}$. On $E$, $|\widetilde f_n(x)-\theta_1|\ge\Delta/2$ by the
triangle inequality. Hence, by Markov's inequality applied to the squared error,
\[
\begin{aligned}
\max_{j=0,1}\E_{f_0^{(j)}}|\widetilde f_n(x)-\theta_j|^2
&\ge\frac12\Big\{\E_{f_0^{(0)}}|\widetilde f_n(x)-\theta_0|^2+\E_{f_0^{(1)}}|\widetilde f_n(x)-\theta_1|^2\Big\}\\
&\ge\frac{\Delta^2}{8}\Big\{P_0^n(E^c)+P_1^n(E)\Big\}
\ge\frac{\Delta^2}{8}\Big\{1-\|P_1^n-P_0^n\|_{\mathrm{TV}}\Big\},
\end{aligned}
\]
the last step because $P_1^n(E)\ge P_0^n(E)-\|P_1^n-P_0^n\|_{\mathrm{TV}}$. Now Le Cam's two-point lemma therefore implies that, for every estimator
$\widetilde f_n(x)$,
\[
\begin{aligned}
\sup_{f_0\in\mathcal H^\alpha(L_0,B)}
\mathbb E_{f_0}
\left[
|\widetilde f_n(x)-f_0(x)|^2
\right]
&\ge
\frac{1}{16}~
|f_0^{(1)}(x)-f_0^{(0)}(x)|^2
=\frac{a^2}{16}\Big(\frac{\sigma^2}{a^2c_1}\Big)^{\frac{2\alpha}{2\alpha+d^*}}n^{-\frac{2\alpha}{2\alpha+d^*}}
\asymp
n^{-\frac{2\alpha}{2\alpha+d^*}}.
\end{aligned}
\]
Thus,
\[
\inf_{\widetilde f_n}
\sup_{f_0\in\mathcal H^\alpha(L_0,B)}
\mathbb E_{f_0}
\left[
|\widetilde f_n(x)-f_0(x)|^2
\right]
\gtrsim
n^{-2\alpha/(2\alpha+d^*)}.
\]
\end{proof}

\section{Proofs for Section~\ref{sec:adaptive}}
\label{app:proofs-adaptive}

Throughout, $x\in\mathcal M$ is fixed, $f_0$ satisfies Definition~\ref{def:local-holder} at
$x$, and the grid $\mathcal K_n$ satisfies \eqref{eq:grid-conditions}. We write
$L'=L_\F+3L_0$.

\begin{proof}[Proof of Lemma~\ref{lem:computable-radius}]
By Proposition~\ref{prop:critical-radius}, $\bar\delta_k$ satisfies the sufficient
condition \eqref{eq:critical-ineq} and hence the critical inequality
$G_k(\bar\delta_k;\partial\mathcal F_n)\le\bar\delta_k^2$ (if $\bar\delta_k\le1$ through
Lemma~\ref{lem:dudley}, and if $\bar\delta_k>1$ directly from the Cauchy--Schwarz bound
$G_k(\bar\delta_k)\le\bar\delta_k<\bar\delta_k^2$). Indeed $\bar\delta_k$ is exactly the radius
\eqref{eq:delta-lambertW} with the local sample size $k$. Since $\delta_k$ is the smallest such
radius, $\bar\delta_k\ge\delta_k$.

For the second claim, observe that the proof of Proposition~\ref{prop:training-error} uses
the critical radius at exactly one point: to bound
$\E[Z_{n,u}]=G_k(u;\partial\mathcal F_n)\le u\delta_n$ for $u\ge\delta_n$. This follows from
the monotonicity of Lemma~\ref{lem:critical-exists} together with $G_k(\delta_n)\le\delta_n^2$.
Any radius $\bar\delta$ with $G_k(\bar\delta)\le\bar\delta^2$ has the same property, because
for $u\ge\bar\delta$,
\[
\frac{G_k(u;\partial\mathcal F_n)}{u}\le\frac{G_k(\bar\delta;\partial\mathcal F_n)}{\bar\delta}\le\bar\delta .
\]
The rest of the proof is unchanged, and so is the proof of Theorem~\ref{thm:oracle}, which
uses the same event. Finally, \eqref{eq:delta-k} follows from the elementary bound
$W(z)\le\log(1+z)$, valid for $z\ge0$, and $\log(1+z/16)\le\log(1+z)$.

\emph{A lower bound.} We also record the lower bound
\begin{equation}
\label{eq:delta-lower}
k\bar\delta_k^2\ge2048\,W\Big(\frac{k}{16}\Big),
\end{equation}
used in Section~\ref{sec:collection}. The number of non-zero parameters satisfies
$s\le\sum_{\ell=0}^{L}(p_\ell+1)p_{\ell+1}$, and each summand is at most
$(p_\ell+1)(p_{\ell+1}+1)\le V$, so $s\le(L+1)V$ and, as $V\ge2$,
$s+1\le(L+1)V+1\le(L+1)V^2$. Hence the argument of $W$ in the definition of $\bar\delta_k$ is at
least $k/16$, and since $W$ is increasing and $s+1\ge1$,
$k\bar\delta_k^2=2048(s+1)W\big(k(L+1)V^2/(16(s+1))\big)\ge2048\,W(k/16)$.
\end{proof}

\begin{proof}[Proof of Proposition~\ref{prop:simultaneous}]
 We construct $\mathcal A_x$ as the intersection of two events, one controlling the
training errors and one controlling the neighbourhood radii, each of probability at least
$1-\eta$.

\emph{ Step 1: training errors.} For each $k\in\mathcal K_n$, let
\[
T_k(x)
=
\frac{1}{k}
\sum_{i\in N_k(x)}
\bigl\{
\widehat f_{x,k}(X_i)-f_0(X_i)
\bigr\}^2.
\]
Because $f_0$ is not assumed to lie in $\mathcal F_n$ here, the relevant training-error bound
is the \emph{oracle} inequality of Section~\ref{sec:oracle}, equation~\eqref{eq:oracle-tail}.
 Applied with the computable radius $\bar\delta_k$ of Lemma~\ref{lem:computable-radius}
it gives, for every $t\ge\bar\delta_k$, conditionally on the covariates,
\begin{equation}
\label{eq:oracle-Tk}
\P\left(T_k(x)>2\inf_{f\in\mathcal F_n}\|f-f_0\|^2_{2,k}
+36\kappa_\sigma\,t\bar\delta_k\ \Big|\ X_1,\dots,X_n\right)\le\exp\Big(-\frac{tk\bar\delta_k}{2}\Big).
\end{equation}
 The constant $36\kappa_\sigma$ dominates $\max\{2,36\sigma^2\}$, and so covers both cases
in the proof of Theorem~\ref{thm:oracle}. With $t=t_k$ the right-hand side is at most
$\eta/|\mathcal K_n|$, and $36\kappa_\sigma t_k\bar\delta_k=\tau_k^2$ by \eqref{eq:tau-scale}. A
union bound over the grid gives an event $\mathcal E_1$, of conditional and hence
unconditional probability at least $1-\eta$, on which
\[
 T_k(x)\le2\inf_{f\in\mathcal F_n}\|f-f_0\|^2_{2,k}+\tau_k^2
\qquad\text{for every }k\in\mathcal K_n .
\]

\emph{ Step 2: neighbourhood radii.}  Let $r_{n,k}=(2k/(c_0n))^{1/d^*}$. By the
second condition in \eqref{eq:grid-conditions}, $r_{n,k}\le\min\{r_0,r_{\min}\}$ and
$k\le\frac12c_0nr_0^{d^*}$ for every $k\in\mathcal K_n$, so Lemma~\ref{lem:knn_radius_tail}
applies with $r=r_{n,k}$ and gives
\[
 \P\big(R_k(x)>r_{n,k}\big)\le\exp\Big(-\frac{c_0}{8}nr_{n,k}^{d^*}\Big)=e^{-k/4}
\le e^{-k_0/4}\le\frac{\eta}{|\mathcal K_n|},
\]
 the last step by the first condition in \eqref{eq:grid-conditions}. A union bound gives
an event $\mathcal E_2$ of probability at least $1-\eta$ on which
$R_k(x)\le r_{n,k}\le r_{\min}\le r_x$ for every $k\in\mathcal K_n$. Set
$\mathcal A_x=\mathcal E_1\cap\mathcal E_2$, so that $\P(\mathcal A_x)\ge1-2\eta$.

\emph{ Step 3: approximation error.}  On $\mathcal E_2$ every $X_i$ with
$i\in N_k(x)$ lies in $B(x,r_x)$. Taking the constant function $f\equiv f_0(x)\in\mathcal F_n$,
\[
\inf_{f\in\mathcal F_n}\|f-f_0\|^2_{2,k}\le
\frac1k\sum_{i\in N_k(x)}\{f_0(x)-f_0(X_i)\}^2\le L_x^2R_k(x)^{2\alpha_x},
\]
 so that on $\mathcal A_x$, $T_k(x)\le2L_x^2R_k(x)^{2\alpha_x}+\tau_k^2$ for every $k$.

\emph{ Step 4: from the neighbourhood to the query point.}  Let
$\Delta_k=\widehat f_{x,k}-f_0$. For $i\in N_k(x)$, since $\widehat f_{x,k}\in\mathcal F_n$ is
$\alpha_x$-H\"older with constant $L_\F$ and $X_i\in B(x,r_x)$,
$$|\Delta_k(x)-\Delta_k(X_i)|\le(L_\F+L_x)\|X_i-x\|^{\alpha_x}\le(L_\F+L_x)R_k(x)^{\alpha_x}.$$
 Averaging over $i\in N_k(x)$ and applying the Cauchy--Schwarz inequality, exactly as in
the proof of Proposition~\ref{prop:pointwise}, gives on $\mathcal A_x$
\begin{equation}
\begin{aligned}
\big|\hat f_k(x)-f_0(x)\big|&\le T_k(x)^{1/2}+(L_\F+L_x)R_k(x)^{\alpha_x}\\
&\le  \sqrt2\,L_xR_k(x)^{\alpha_x}+\tau_k+(L_\F+L_x)R_k(x)^{\alpha_x}
 {\ \le\ (L_\F+3L_x)R_k(x)^{\alpha_x}+\tau_k.}
\end{aligned}
\label{eq:train-to-test-k}
\end{equation}

\emph{ Step 5: the rate of the radius.}  On $\mathcal E_2$,
$R_k(x)^{\alpha_x}\le(2/c_0)^{\alpha_x/d^*}(k/n)^{\alpha_x/d^*}\le C_{\mathcal M}(k/n)^{\alpha_x/d^*}$,
because $0<\alpha_x\le1$. Substituting into \eqref{eq:train-to-test-k} gives
\eqref{eq:simultaneous}.
\end{proof}

\begin{proof}[Proof of Lemma~\ref{lem:tau-properties}]
 Write $\tau_k^2=36\kappa_\sigma\max\{a_k,\beta_k\}$ with $a_k=\bar\delta_k^2$ and
$\beta_k=2k^{-1}\log(|\mathcal K_n|/\eta)$. Put $c=(L+1)V^2/(16(s+1))$, so that
\[
a_k=2048(s+1)\,c\,\frac{W(ck)}{ck}.
\]
The Lambert function satisfies $W(z)e^{W(z)}=z$, hence $W(z)/z=e^{-W(z)}$, and $W$ is
increasing on $(0,\infty)$. Therefore $k\mapsto a_k$ is non-increasing, and so is
$k\mapsto\beta_k$. The maximum of two non-increasing sequences is non-increasing, which
gives $\tau_{2k}\le\tau_k$.

For the lower bound, $a_{2k}/a_k=\exp\{W(ck)-W(2ck)\}$. Taking logarithms in
$W(z)e^{W(z)}=z$ gives $W(z)=\log z-\log W(z)$, so
\[
W(2z)-W(z)=\log2-\log\frac{W(2z)}{W(z)}\le\log 2 ,
\]
because $W(2z)\ge W(z)$. Hence $a_{2k}\ge a_k/2$. Since also $\beta_{2k}=\beta_k/2$,
\[
\max\{a_{2k},\beta_{2k}\}\ge\tfrac12\max\{a_k,\beta_k\},
\]
that is, $\tau_{2k}^2\ge\tau_k^2/2$.

Finally, $\max\{a_k,\beta_k\}\le a_k+\beta_k$. Bounding $a_k$ by \eqref{eq:delta-k} with its
explicit constant $2048$, and using $2\le2048$,
\[
\tau_k^2\le36\kappa_\sigma\Big[\frac{2048(s+1)}{k}\log\Big(1+\frac{k(L+1)V^2}{s+1}\Big)
+\frac2k\log\Big(\frac{|\mathcal K_n|}{\eta}\Big)\Big],
\]
which is \eqref{eq:tau-explicit} with $C_1=36\cdot2048$.
\end{proof}

\begin{proof}[Proof of Lemma~\ref{lem:admissible}]
 Because $b_{k_0}(x)\le\tau_{k_0}$, the set in \eqref{eq:oracle-index} is non-empty and
$k_x^*$ is well defined. Suppose that the event $\mathcal A_x$ in \eqref{eq:simultaneous}
holds and let $\ell\le k_x^{*}$. Then
\[
\big|\hat f_{k_x^{*}}(x)-\hat f_\ell(x)\big|
\le\big|\hat f_{k_x^{*}}(x)-f_0(x)\big|+\big|\hat f_\ell(x)-f_0(x)\big|
\le b_{k_x^{*}}(x)+\tau_{k_x^{*}}+b_\ell(x)+\tau_\ell .
\]
By the definition of $k_x^{*}$, $b_{k_x^{*}}(x)\le\tau_{k_x^{*}}$, and because $\ell\le k_x^{*}$,
$b_\ell(x)\le b_{k_x^{*}}(x)\le\tau_{k_x^{*}}\le\tau_\ell$  (the last step by
Lemma~\ref{lem:tau-properties}). Consequently
$|\hat f_{k_x^{*}}(x)-\hat f_\ell(x)|\le4\tau_\ell$. Therefore, if $A\ge4$, the oracle
neighbourhood $k_x^{*}$ is admissible in \eqref{eq:lepski-rule}, and hence
$\hat k(x)\ge k_x^{*}$.
\end{proof}

\begin{proof}[Proof of Lemma~\ref{lem:balance}]
 The upper bound is the definition of $k_x^*$. For the lower bound, Assumption~\ref{ass:window}
gives $k_x^*\le k_J/2$, so $2k_x^*\in\mathcal K_n$, and by maximality
$b_{2k_x^{*}}(x)>\tau_{2k_x^{*}}$.  From \eqref{eq:bias-bk},
$b_{2k}(x)=2^{\alpha_x/d^*}b_k(x)\le2b_k(x)$, since $\alpha_x\le1\le d^*$; and
$\tau_{2k}\ge\tau_k/\sqrt2$ by Lemma~\ref{lem:tau-properties}. Hence
\[
 b_{k_x^*}(x)\ge\tfrac12 b_{2k_x^*}(x)>\tfrac12\tau_{2k_x^*}\ge\frac{1}{2\sqrt2}\,\tau_{k_x^*}.
\]
 The equivalence in \eqref{eq:balance} follows by substituting \eqref{eq:bias-bk}.
This step is where $k_x^{*}<k_J$ is needed: if the oracle index sits at the top of the grid
the comparison $b_{2k_x^{*}}>\tau_{2k_x^{*}}$ is unavailable and only ``$\lesssim$'' survives.
\end{proof}

\begin{proof}[Proof of Theorem~\ref{thm:adaptive-oracle}]
 By Lemma~\ref{lem:admissible}, $k_x^{*}\le\hat k(x)$,  so the definition of
$\hat k(x)$  with $\ell=k_x^*$ gives
$|\hat f_{\hat k(x)}(x)-\hat f_{k_x^{*}}(x)|\le A\tau_{k_x^{*}}$. Combining with \eqref{eq:simultaneous},
\[
\big|\hat f_{\mathrm{ad}}(x)-f_0(x)\big|
\le A\tau_{k_x^{*}}+b_{k_x^{*}}(x)+\tau_{k_x^{*}}\le(A+2)\tau_{k_x^{*}},
\]
which is \eqref{eq:adaptive-tau}.  For \eqref{eq:oracle-adaptive} it suffices to show
that $\tau_{k_x^*}^2\le8\inf_{k\in\mathcal K_n}\{b_k(x)^2+\tau_k^2\}$. For $k\le k_x^{*}$
we have $\tau_k\ge\tau_{k_x^{*}}$, so $b_k^2+\tau_k^2\ge\tau_{k_x^{*}}^2$; and for $k>k_x^{*}$ we
have $b_k(x)\ge b_{k_x^{*}}(x) >\tau_{k_x^*}/(2\sqrt2)$ by \eqref{eq:balance}, so
$b_k^2+\tau_k^2 \ge\tau_{k_x^{*}}^2/8$ as well. Taking the infimum over $k\in\mathcal K_n$
and combining with \eqref{eq:adaptive-tau} gives the first inequality in
\eqref{eq:oracle-adaptive}; the second uses $L_x\le L_0$. Substituting
\eqref{eq:tau-explicit} gives \eqref{eq:oracle-explicit}.
\end{proof}

\begin{proof}[Proof of Theorem~\ref{thm:adaptive-rate}]
 Since $k\le n$ on the grid, \eqref{eq:tau-explicit} gives $\tau_k^2\le C_1\mathcal D_n/k$.
By Theorem~\ref{thm:adaptive-oracle} and $L_x\le L_0$, on $\mathcal A_x$,
\[
\big|\hat f_{\mathrm{ad}}(x)-f_0(x)\big|^2\le8(A+2)^2\inf_{k\in\mathcal K_n}\varphi(k),
\qquad
\varphi(k)=C_{\mathcal M}^2L'^2\Big(\frac kn\Big)^{2\alpha_x/d^*}+\frac{C_1\mathcal D_n}{k}.
\]
The two terms of $\varphi$ are equal at
\[
k_x^\circ=(C_1\mathcal D_n)^{\frac{d^*}{2\alpha_x+d^*}}(C_{\mathcal M}L')^{-\frac{2d^*}{2\alpha_x+d^*}}n^{\frac{2\alpha_x}{2\alpha_x+d^*}},
\qquad
\frac{C_1\mathcal D_n}{k_x^\circ}=(C_{\mathcal M}L')^{\frac{2d^*}{2\alpha_x+d^*}}\Big(\frac{C_1\mathcal D_n}{n}\Big)^{\frac{2\alpha_x}{2\alpha_x+d^*}} .
\]
Because $k_0\le k_x^\circ\le k_J$ and the grid is geometric with ratio $2$, there is
$k\in\mathcal K_n$ with $k_x^\circ\le k<2k_x^\circ$. For this $k$ the first term of $\varphi$
is at most $2^{2\alpha_x/d^*}\le4$ times its value at $k_x^\circ$ and the second is at most
its value at $k_x^\circ$, so $\varphi(k)\le5C_1\mathcal D_n/k_x^\circ$. Hence
\[
\big|\hat f_{\mathrm{ad}}(x)-f_0(x)\big|^2\le40(A+2)^2(C_{\mathcal M}L')^{\frac{2d^*}{2\alpha_x+d^*}}\Big(\frac{C_1\mathcal D_n}{n}\Big)^{\frac{2\alpha_x}{2\alpha_x+d^*}},
\]
which is \eqref{eq:adaptive-rate}, with a constant depending only on $A$ and
$C_{\mathcal M}$, since $C_1$ is absolute.
\end{proof}

\begin{proof}[Proof of Corollary~\ref{cor:adaptive-OP}]
 Since $|\mathcal K_n|\lesssim\log n$ and $\eta$ is fixed,
\[
\mathcal D_n\le\kappa_\sigma\Big[(s+1)\big\{\log(1+n)+\log(L+1)+2\log V\big\}+\log\log n+O(1)\Big]
=O\big((\log n)^{C''}\big)
\]
under the stated growth conditions. Substituting into \eqref{eq:adaptive-rate} gives a
deterministic sequence $\rho_n=C\,n^{-2\alpha_x/(2\alpha_x+d^*)}(\log n)^{C}$ with
$|\hat f_{\mathrm{ad}}(x)-f_0(x)|^2\le\rho_n$ on the event $\mathcal A_x$ of
Proposition~\ref{prop:simultaneous}.

Now for $\mathcal E_1$, the union bound in Step~1 of that proof gives, since $t_k\ge\bar\delta_k$,
\[
\P(\mathcal E_1^c)\le\sum_{k\in\mathcal K_n}\exp\Big(-\frac{t_kk\bar\delta_k}{2}\Big)
\le\sum_{k\in\mathcal K_n}\exp\Big(-\frac{k\bar\delta_k^2}{2}\Big)
\le|\mathcal K_n|\exp\big\{-1024\,W(k_0/16)\big\},
\]
by \eqref{eq:delta-lower} and the monotonicity of $W$. For $z\ge e$ we have $W(z)\le\log z$, hence
$W(z)=\log z-\log W(z)\ge\log z-\log\log z$, so $e^{-W(z)}\le\log z/z$. With $z=k_0/16$, which tends
to infinity,
\[
\P(\mathcal E_1^c)\le|\mathcal K_n|\Big(\frac{16\log(k_0/16)}{k_0}\Big)^{1024}\to0,
\]
because $|\mathcal K_n|\lesssim\log n$ while $k_0\ge8\log(n/\eta)$.

For $\mathcal E_2$, Step~2 of the same proof gives, since $k_j=2^jk_0\ge(j+1)k_0$,
\[
\P(\mathcal E_2^c)\le\sum_{j=0}^{J}e^{-k_j/4}\le\sum_{j\ge0}e^{-(j+1)k_0/4}
=\frac{e^{-k_0/4}}{1-e^{-k_0/4}}\to0 .
\]
Hence $\P(\mathcal A_x)\to1$. Now, given $\varepsilon>0$, choose $n_\varepsilon$ such that
$\P(\mathcal A_x^c)\le\varepsilon$ and the window conditions hold for $n\ge n_\varepsilon$; then
$\P\big(|\hat f_{\mathrm{ad}}(x)-f_0(x)|^2>\rho_n\big)\le\varepsilon$ for all $n\ge n_\varepsilon$, which is
\eqref{eq:adaptive-OP}.

\end{proof}

\end{document}